\documentclass[11pt,reqno]{amsart}
\usepackage[T1]{fontenc}
\usepackage{amsmath,amssymb,mathtools}
\usepackage{mathrsfs}
\usepackage[margin=1.05in]{geometry}
\usepackage{tikz}
\usepackage{xcolor}
\usepackage[colorlinks=true,linkcolor=blue!45!black,citecolor=green!35!black,urlcolor=blue!45!black]{hyperref}
\hypersetup{bookmarksdepth=2}
\numberwithin{equation}{section}
\newtheorem{theorem}{Theorem}[section]
\newtheorem{proposition}[theorem]{Proposition}
\newtheorem{lemma}[theorem]{Lemma}
\newtheorem{corollary}[theorem]{Corollary}
\newtheorem{mtheorem}{Theorem}
\newtheorem{mcorollary}[mtheorem]{Corollary}

\theoremstyle{definition}

\theoremstyle{remark}
\newtheorem{remark}[theorem]{Remark}
\newcommand{\CC}{\mathbb C}
\newcommand{\PP}{\mathbb P}
\newcommand{\RR}{\mathbb R}
\newcommand{\HH}{\mathcal H}
\newcommand{\HHH}{\mathscr H}

\newcommand{\MM}{\mathcal M}
\newcommand{\M}{\MM}
\newcommand{\E}{\mathcal E}
\newcommand{\ddc}{dd^c}
\DeclareMathOperator{\Ent}{Ent}
\DeclareMathOperator{\Ric}{Ric}
\DeclareMathOperator{\tr}{tr}
\DeclareMathOperator{\supp}{supp}
\DeclareMathOperator{\PSH}{PSH}

\title[Uniqueness of scalar-flat ALE K\"ahler metrics]
{Uniqueness of scalar-flat ALE K\"ahler metrics}
\author{Qi Yao}
\address{IASM, Zhejiang University, Hangzhou, Zhejiang, China}
\email{qiyao11@zju.edu.cn}
\date{}

\begin{document}

\begin{abstract}
We study uniqueness of scalar-flat K\"ahler metrics in a fixed
K\"ahler class on an asymptotically locally Euclidean (ALE)
K\"ahler manifold, assuming only differentiated $o(1)$ decay to
a fixed Euclidean metric on the end. In complex dimension $n\geq2$,
under suitable polynomial decay of the reference metric, we prove
that the relative Mabuchi functional is finite, continuous, and
convex along weak geodesics whose endpoint potentials satisfy
$u=O(r^{-\gamma})$, $\gamma>2n-4$, together with the corresponding
derivative estimates. Combining this convexity with the analysis
of the Lichnerowicz operator in weighted spaces, we obtain
uniqueness by a perturbation argument. Decay improvement and
comparison of the asymptotic charts produce a global potential with
the required decay for the difference of the two cohomologous
scalar-flat K\"ahler forms. The argument also
extends to asymptotically conical K\"ahler metrics approaching a
fixed Ricci-flat K\"ahler cone under prescribed decay assumptions.
\end{abstract}

\maketitle

\setcounter{tocdepth}{1}
\tableofcontents

\section{Introduction}\label{sec:introduction}

Let $(X,J,g)$ be a connected complete noncompact K\"ahler manifold of
complex dimension $n\geq2$ wit+h the reference K\"ahler metric $g$. We assume that $g$ is asymptotically
locally Euclidean (ALE): there exists a compact subset $K\subset X$,
a radius $R_0>0$, and a diffeomorphism
\[
 \Upsilon:X\setminus K\longrightarrow
       (\CC^n\setminus\overline B_{R_0})/\Gamma,
\]
where $\Gamma\subset U(n)$ is a finite group acting freely on the
unit sphere, such that
\begin{equation}\label{eq:ale-decay}
 \left|\nabla_{\mathrm{euc}}^j
 \big((\Upsilon^{-1})^*g-g_{\mathrm{euc}}\big)\right|_{g_{\mathrm{euc}}}
 =O(r^{-\delta-j}),\qquad 0\leq j\leq k,\qquad \delta>n-1.
\end{equation}
We call $X\setminus K \cong (\CC^n \setminus \overline{B}_{R_0})/\Gamma$ an \textit{end} of $X$.   
Here $r$ is the Euclidean radius from the end extending smoothly across the compact set $K$, $\nabla_{\mathrm{euc}}$ is the
Euclidean connection, and $k$ is a sufficiently large integer. We also assume that the scalar
curvature belongs to $L^1(X,\omega^n)$. We fix the ALE chart
$\Upsilon$ and the smooth extension of $r$ throughout the paper,
and use $\omega(\cdot, \cdot) = g(J \cdot, \cdot) $ as a reference K\"ahler form.

The choice $\delta>n-1$ is motivated by the Arnowitt--Deser--Misner
(ADM) mass in general relativity \cite{ADM61}. In real dimension $2n$, the
threshold is $(2n-2)/2=n-1$: the strict inequality $\delta>n-1$ is
the weakest fall-off condition which guarantees that the
mass is finite and independent of the admissible coordinates at
infinity. The coordinate-invariance results are due to Bartnik \cite{Bartnik86}
and Chru\'sciel \cite{Chrusciel86}; see also \cite[Section~1]{YaoMass} for their use in
the ALE setting.

The search for canonical representatives of K\"ahler classes is
a central problem in K\"ahler geometry. One natural condition is
that the metric have constant scalar curvature. For an ALE metric,
the scalar curvature tends to zero at infinity, so any constant
scalar curvature ALE K\"ahler metric is necessarily scalar-flat.
The study of scalar-flat ALE K\"ahler metrics thus forms a
noncompact counterpart of the constant scalar curvature problem
on compact K\"ahler manifolds.

In this paper, we study the uniqueness of these metrics within a
fixed K\"ahler class. More precisely, let $\omega_0$ and $\omega_1$
be scalar-flat ALE K\"ahler forms on the fixed complex manifold
$(X,J)$, with $[\omega_0]=[\omega_1]$. Assuming that their metrics
approach the Euclidean metric in the fixed ALE coordinates
and satisfy the prescribed decay conditions, we ask whether
$\omega_0=\omega_1$. Our main result establishes this uniqueness
under the assumptions above.

\subsection{Uniqueness in the compact case}
Calabi's study of extremal K\"ahler metrics \cite{Calabi82,Calabi85}
provides a general setting for the existence and uniqueness problem.
For K\"ahler--Einstein metrics, Aubin and Yau \cite{Aubin78,Yau78}
settled the negative first Chern class case, while Yau proved existence
and uniqueness of Ricci-flat metrics in each K\"ahler class when
$c_1(X)=0$. Bando--Mabuchi \cite{BandoMabuchi} proved uniqueness
of positive K\"ahler--Einstein metrics modulo the identity component
of the holomorphic automorphism group. The work of Matsushima,
Lichnerowicz, and Calabi \cite{Matsushima57,Lichnerowicz57,Calabi85}
and Futaki's invariant \cite{Futaki83} describe obstructions arising
from holomorphic vector fields.

Mabuchi \cite{Mabuchi85,Mabuchi86} introduced the $K$-energy,
whose first variation on a compact K\"ahler manifold is
\[
 d\MM_{\rm cpt}|_u(v)
       =-\int_Xv(S_{\omega_u}-\overline S)\omega_u^n,
 \qquad
 \omega_u=\omega+\ddc u,
\]
where $\overline S$ is the average scalar curvature of the class.
Its critical points are therefore precisely the cscK metrics.
Mabuchi \cite{Mabuchi87} also introduced the $L^2$ metric on the
space of K\"ahler potentials and established convexity of the
$K$-energy along smooth geodesics, with degeneracy determined by
holomorphic vector fields.

Semmes \cite{Semmes} and Donaldson \cite{Donaldson99} related
geodesics to the homogeneous complex Monge--Amp\`ere equation,
whose weak formulation uses the pluripotential theory of
Bedford--Taylor \cite{BT76,BT}. Chen \cite{Chen00} constructed
weak geodesics between smooth K\"ahler potentials with bounded
complex Hessian. B{\l}ocki obtained gradient estimates
\cite{Blocki09} and developed the geodesic regularity theory
\cite{Blocki12}; Chu--Tosatti--Weinkove \cite{CTW} subsequently
proved the joint $C^{1,1}$ bound. This regularity is optimal in
general: the examples of Lempert--Vivas and Darvas--Lempert
\cite{LempertVivas,DarvasLempert} show that smooth endpoint
potentials need not be joined by a $C^2$ geodesic. Moreover, in joint work with Xia Xiao
\cite[Theorem~A]{XiaoYao}, we showed that the spatial K\"ahler form
can degenerate on an open subset during an interval of interior
times.

In polarized classes, Donaldson \cite{Donaldson01} approximated
cscK metrics by balanced metrics when the polarized automorphism
group is discrete, and Mabuchi \cite{Mabuchi04} extended this
approach to extremal metrics in integral K\"ahler classes.
Phong--Sturm \cite{PhongSturm06,PhongSturm07} further connected
projective approximation with weak geodesics.
Berman--Berndtsson \cite{BB} proved convexity of the Mabuchi
functional along weak geodesics in arbitrary compact K\"ahler
classes, and uniqueness of cscK and extremal metrics modulo
holomorphic automorphisms. Their convexity proof combines
Berndtsson's plurisubharmonic variation of Bergman kernels
\cite{Ber} with Bergman approximation of Monge--Amp\`ere measures;
their uniqueness argument uses a strictly convex perturbation
and the Lichnerowicz equation. Chen--Li--P\u{a}un \cite{CLP}
provided another proof of convexity.

\subsection{Earlier noncompact results}
Prescribed asymptotics are part of any noncompact uniqueness problem. 
Auvray \cite{Auvray17} proved the $\mathcal C^{1,1}$ regularity of weak geodesics for Poincar\'e-type metrics on the complement of a divisor;
the uniqueness theorem there includes an ampleness hypothesis on the
logarithmic canonical bundle. 
Xu--Zheng \cite[Theorem~1.2]{XuZheng} proved uniqueness of
Poincar\'e-type cscK metrics in a fixed K\"ahler class on the
complement of a divisor with simple normal crossings, modulo
automorphisms of the compactification preserving the divisor. 

In the ALE setting, Aleyasin \cite{Aleyasin} considered a more restrictive decay class and proved scalar-flat uniqueness under a weighted nonpositive first Chern class hypothesis; however, the restrictive assumption ruled out many interesting cases.
In our earlier work \cite{YaoGeo}, we established joint $C^{1,1}$ regularity for ALE weak geodesics and studied weighted estimates for smooth $\varepsilon$-geodesics. Those weighted constants can depend on $\varepsilon^{-1}$. The convexity result there uses a nonpositive Ricci reference and is not a general theorem for arbitrary ALE weak geodesics without that sign condition.

The geometry at infinity also has an independent development.
Hein--LeBrun \cite{HeinLeBrun} expressed the mass of an ALE K\"ahler metric in terms of scalar curvature and cohomological data. 
In \cite{YaoMass}, we proved the global-potential and expansion results used to identify the decay of the difference of cohomologous scalar-flat metrics. 
Conlon--Hein \cite{ConlonHein} established existence and uniqueness results for asymptotically conical Calabi--Yau metrics with prescribed asymptotic behavior.
Han--Viaclovsky \cite{HanViaclovsky} developed existence and compactness theory for scalar-flat ALE K\"ahler surfaces. 
None of these results establishes uniqueness for scalar-flat
ALE K\"ahler metrics under the general assumptions considered
in this paper.

\subsection{The difficulties caused by infinite volume}
The compact theory does not directly establish finiteness of the Mabuchi functional
along weak geodesics on an ALE manifold. Relative entropy removes the
contribution of the infinite background volume, but the decay of smooth
endpoint potentials need not persist for the complex Hessian of the
geodesic, as the examples in \cite{XiaoYao} show. Thus the
available regularity theory does not ensure integrability of the
entropy. We approach the geodesic by weak geodesics joining compactly
supported endpoint potentials and use uniform estimates at infinity
to establish finiteness of the limiting functional.

The principal difficulty in proving convexity is convergence of the
entropy at interior times. Uniform convergence of potentials implies
only weak convergence of their volume measures, which does not imply
convergence of entropy; on a noncompact manifold one must also exclude
loss of entropy at infinity. Neither issue is resolved by the compact
convexity theorem. We use the envelope description of weak geodesics
to identify the limiting volume densities on compact subsets, and
localize the convexity argument to control the contribution from
infinity. Together these arguments prove convergence of the entropy
along the approximating geodesics, and hence convexity and continuity
of the relative Mabuchi functional.

Passing from convexity to uniqueness requires a further extension of
the compact perturbation argument of Berman--Berndtsson \cite{BB}.
On an ALE manifold, a solution of the linearized scalar curvature
equation need not have the decay required of the K\"ahler potentials,
even when the source is compactly supported. We solve the
Lichnerowicz equation with controlled behavior at infinity and then
construct compactly supported perturbation terms. The essential point
is to control the resulting error in the first variation as the
perturbation tends to zero. This makes it possible to use strict
convexity on each compact subset and obtain uniqueness without an
additional decay assumption on the Hessians of weak geodesics.

\subsection{Main Theorems}
Since we discuss scalar-flat metrics, the relevant first variation
is the compact formula with zero target scalar curvature,
$\overline S=0$: $d\MM|_u(v)=-\int_Xv\,S_{\omega_u}\,\omega_u^n$
for a compactly supported test function $v$, where
$S_\eta=\operatorname{tr}_\eta\Ric(\eta)$ is our scalar curvature
convention. The critical points of this functional are therefore the scalar-flat metrics. 
On an ALE manifold, however, the ordinary entropy integral may diverge because the volume is infinite. 
We use the relative entropy, obtained by subtracting the linear term of the entropy integrand at the reference density.
Fix $\gamma=2n-4+\varepsilon'$ with $\varepsilon'>0$, and let
$\HH_{-\gamma}(X,\omega)$ consist of smooth real K\"ahler potentials
satisfying $|\nabla^ju|=O(r^{-\gamma-j})$ for $0\leq j\leq k$.
For a smooth potential $u\in\HH_{-\gamma}(X,\omega)$, write
$\omega_u=\omega+\ddc u$ and set
\[
 f_u=\frac{\omega_u^n}{\omega^n},\qquad
 \Ent_{\mathrm{rel}}(u)=\int_X\bigl(f_u\log f_u-f_u+1\bigr)\,\omega^n.
\]
The relative Mabuchi functional is defined by
\[
 \MM(u)=\Ent_{\mathrm{rel}}(u)-nE_{\Ric(\omega)}(u),
 \qquad
 E_\rho(u)=\frac1n\sum_{j=0}^{n-1}
 \int_Xu\,\rho\wedge\omega_u^{n-1-j}\wedge\omega^j,
\]
where $\rho$ is a smooth closed real $(1,1)$-form.

For these smooth potentials, $f_u-1=O(r^{-\gamma-2})$, so the
relative entropy density is $O(r^{-2\gamma-4})$ and is integrable
since $2\gamma+4>2n$. If the reference metric $\omega$ is
scalar-flat, then $\Ric(\omega)\wedge\omega^{n-1}=0$.
Consequently, every nonzero term in the expansion of
$E_{\Ric(\omega)}(u)$ contains a factor
$\ddc u$. Its density relative to $\omega^n$
is $O(r^{-2\gamma-\delta-4})$, which is integrable under our decay
assumptions. Thus $E_{\Ric(\omega)}(u)$ is finite as well.
For a general reference metric, the additional term linear in
$u$ is also integrable under the standing assumption
$S_\omega\in L^1(X,\omega^n)$.
Our first result establishes the finiteness of this functional
along weak geodesics with endpoint potentials in $\HH_{-\gamma}(X,\omega)$.

\begin{mtheorem}\label{thm:main-a}
Let $(X,J,\omega)$ be an ALE K\"ahler manifold satisfying the
preceding geometric assumptions, and let
$u_t$, $t\in[0,1]$, be a weak geodesic whose endpoint potentials satisfy
\[
 u_0,u_1\in\HH_{-\gamma}(X,\omega).
\]
Then the Mabuchi functional $\MM(u_t)$ is well-defined for every
$t\in[0,1]$.
\end{mtheorem}

By Theorem~\ref{thm:main-a}, the relative Mabuchi functional defines a
finite real-valued function along every weak geodesic with endpoint
potentials in $\HH_{-\gamma}(X,\omega)$. We next establish its
continuity and convexity, the variational properties underlying the
uniqueness problem.

\begin{mtheorem}\label{thm:main-b}

Under the preceding geometric assumptions, the Mabuchi functional
is continuous and convex along the weak geodesic $u_t$ connecting potentials
$u_0,u_1\in\HH_{-\gamma}(X,\omega)$.
\end{mtheorem}

For two cohomologous scalar-flat ALE K\"ahler forms satisfying
\eqref{eq:ale-decay} in a common chart, the global-potential theorem
\cite[Theorem~D(i)]{YaoMass} provides a global difference potential in
$\HH_{-\gamma}$ for some $\gamma>2n-4$. Thus the prescribed metric
decay places the geometric uniqueness problem in the space of
potentials considered above. Our main theorem answers this problem.

\begin{mtheorem}[Uniqueness of scalar-flat ALE K\"ahler metrics]
\label{thm:uniqueness}
Let $\omega_0$ and $\omega_1$ be scalar-flat ALE K\"ahler
forms on a complex manifold $(X,J)$ of complex dimension
$n\geq2$. Suppose that their associated metrics $g_0,g_1$ are
asymptotic to the fixed Euclidean metric $g_{\mathrm{euc}}$ with
decay $-\delta$ as in \eqref{eq:ale-decay}, where $\delta>n-1$
and $k$ is sufficiently large. If $[\omega_0]=[\omega_1]$, then
$\omega_0=\omega_1$.
\end{mtheorem}

The decay assumption in Theorem~\ref{thm:uniqueness} can be weakened
in the ALE case. The following statement uses the decay improvement
of Ach\'e--Viaclovsky \cite{AcheViaclovsky}.

\begin{mcorollary}[Uniqueness under differentiated $o(1)$ asymptotics]
\label{cor:order-zero}
Let $\omega_0,\omega_1$ be complete scalar-flat K\"ahler forms on a
connected noncompact complex manifold $(X,J)$ of complex dimension $n\geq2$.
Suppose that there is a compact subset $K\subset X$ such that
$X\setminus K$ is diffeomorphic to
$(\CC^n\setminus\overline B_{R_0})/\Gamma$, and that their associated
metrics $g_0,g_1$ are asymptotic to a fixed Euclidean metric
$g_{\mathrm{euc}}$ on each end, in the sense that
\begin{equation}\label{eq:order-zero-pair}
 \left|\nabla_{\mathrm{euc}}^j
 (g_i-g_{\mathrm{euc}})\right|_{g_{\mathrm{euc}}}
 =o(r^{-j}),\qquad i=0,1,\quad 0\leq j\leq k,
\end{equation}
for sufficiently large $k$, where $r$ and $\nabla_{\mathrm{euc}}$
denote the Euclidean radius and connection on each end.
If $[\omega_0]=[\omega_1]$, then $\omega_0=\omega_1$.
\end{mcorollary}

We also consider an asymptotically conical (AC) K\"ahler manifold.

Fix a K\"ahler cone $(C(L),J_C,g_C)$ with smooth compact
link $L$, where $g_C=dr^2+r^2g_L$. We assume that $g_C$ is
Ricci-flat. Fix an identification
$\Phi:(R_0,\infty)\times L\longrightarrow X\setminus K$.
The AC assumptions below are
\begin{equation}\label{eq:ac-introduction}
 \Phi^*J\longrightarrow J_C,\qquad
 |\nabla_C^j(\Phi^*g-g_C)|_{g_C}=O(r^{-\delta-j}),\quad
 0\leq j\leq k,\qquad \delta>n-1,
\end{equation}
where $\nabla_C$ is the cone connection and $k$ is sufficiently
large. We extend $r$ smoothly over $X$, with $r\geq1$, and define
$\HH_{-\gamma}$ using this radius and the same differentiated
potential decay as above.

\begin{mcorollary}[The AC case]\label{cor:ac}
Let $(X,J,\omega)$ be a connected complete AC K\"ahler manifold of
complex dimension $n\geq2$, with one end satisfying
\eqref{eq:ac-introduction}, and assume
$S_\omega\in L^1(X,\omega^n)$. For every $\gamma>2n-4$, the relative
Mabuchi functional is finite, continuous, and convex along the weak
geodesic joining any $u_0,u_1\in\HH_{-\gamma}(X,\omega)$.
Moreover, two cohomologous scalar-flat K\"ahler forms on $(X,J)$
whose metrics satisfy \eqref{eq:ac-introduction} relative to the
same cone and the same identification $\Phi$ coincide.
\end{mcorollary}

\begin{remark}\label{rem:ac-decay}
The decay improvement used in Corollary~\ref{cor:order-zero} has no
general AC analogue. Our line-bundle constructions in
\cite{YaoInvariant}, and the Calabi--Yau examples of van Coevering
\cite{VanCoevering}, include metrics with a nonzero term of order
$r^{-2}$ relative to their cone. A nonzero restriction of the
K\"ahler class to the link prevents faster decay. These examples
do not contradict uniqueness with a fixed K\"ahler class and cone;
they show why differentiated $o(1)$ asymptotics alone cannot be
upgraded automatically to the decay assumed in
Corollary~\ref{cor:ac}. The optimal assumptions for AC uniqueness
are not addressed here.
\end{remark}

\subsection{Outline of the proof of convexity}\label{subsec:proof-outline}
To prove Theorem~\ref{thm:main-b}, we cut off the endpoint
potentials and join them by weak geodesics $u_{t,l}$. Each
approximating geodesic is supported in a compact subset independent
of time, so the local argument of Berman--Berndtsson \cite{BB}
implies convexity of $t\mapsto\MM(u_{t,l})$. The potentials
$u_{t,l}$ converge uniformly to $u_t$, and by the weighted estimates
we have $E_{\Ric(\omega)}(u_{t,l})\to E_{\Ric(\omega)}(u_t)$
uniformly in $t$. Since the endpoint Mabuchi energies converge,
convexity bounds the approximating energies from above and hence
bounds their relative entropies uniformly. Lower semicontinuity
then implies finiteness of the entropy for $u_t$, and the entropy estimate for
the complex Hessian verifies finiteness of the integrals defining
$\MM(u_t)$. To pass convexity to the limit, it remains to prove
convergence of the entropy at every interior time.

By the uniform energy bounds and convexity, we have locally uniform
subsequential limits $\MM(u_{t,l})\to G(t)$ on $(0,1)$.
Set $D(t)=G(t)-\MM(u_t)\geq0$, where nonnegativity follows from
lower semicontinuity. The rooftop representation uses the exact
agreement of the endpoint potentials on expanding compact sets
to prove convergence of entropy on every compact subset
(Lemma~\ref{lem:compact-density-convergence}). Thus $D(t)$
measures only the possible loss of entropy at infinity.
We apply the local convexity argument outside a large compact
set; the weighted estimates and the integral control of the
complex Hessian bound the cutoff errors by a quantity tending
to zero uniformly in $l$ (Lemma~\ref{lem:uniform-entropy-tails}).
Passing first to the subsequential limit and then letting the
cutoff radius tend to infinity, we obtain
\[
 D(t)\leq\frac{b-t}{b-a}D(a)+\frac{t-a}{b-a}D(b),
 \qquad 0<a<t<b<1.
\]
Convergence of the endpoint energies and lower semicontinuity
at the endpoints also imply $D(t)\to0$ as $t\downarrow0$
or $t\uparrow1$.
Hence $D\equiv0$: a nonnegative convex function with zero boundary
limits must vanish. Every subsequential limit therefore equals
$\MM(u_t)$, so the Mabuchi energies converge and their convexity
passes to the limit. Convexity and lower semicontinuity also imply
continuity at the endpoints, completing Theorem~\ref{thm:main-b}.

\subsection*{Acknowledgements}
The author thanks Hans-Joachim Hein, Bianca Santoro, and Song Sun for their support. Special thanks are
due to Song Sun for drawing the author's attention to the work of
Ach\'e--Viaclovsky \cite{AcheViaclovsky}, which allows the ALE
uniqueness theorem to hold under differentiated $o(1)$ asymptotics.

\subsection*{Statement on AI assistance}
The mathematical ideas and the formulation of the statements are
due to the author. AI assistance contributed in two ways. First,
it helped identify the failure of the expected decay of the complex
Hessian along weak geodesics, leading the author to abandon a proof
strategy based on that decay. Second, it accelerated the development
and checking of proof details and assisted with polishing the
language. The author takes full responsibility for the mathematical content
and the final presentation.

\section{Preliminaries}\label{sec:preliminary}

We retain the complete ALE K\"ahler manifold $(X,J,\omega)$ of
complex dimension $n\geq2$, the fixed ALE chart, and the radius
$r\geq1$ from the introduction. Write $g$ for the associated metric
and retain its decay rate $\delta$ from \eqref{eq:ale-decay}. Thus
\begin{equation}\label{eq:ale-setting}
 S_\omega\in L^1(X,\omega^n),\qquad
 |\nabla^j(g-g_{\mathrm{euc}})|_{g_{\mathrm{euc}}}
 =O(r^{-\delta-j}),\quad 0\leq j\leq k,\qquad \delta>n-1,
\end{equation}
where $k$ is sufficiently large for the estimates below. We use
$dd^c=\sqrt{-1}\partial\bar\partial$, $\omega_u=\omega+dd^cu$, and
$S_\omega=\tr_\omega\Ric(\omega)$; the Riemannian scalar curvature
is $2S_\omega$. Derivatives on the end may be taken with respect to
$g$ or $g_{\mathrm{euc}}$. In particular,
\begin{equation}\label{eq:ricci-decay}
 |\Ric(\omega)|_\omega\leq Cr^{-\delta-2},\qquad
 \operatorname{Vol}_\omega\{r<L\}\leq CL^{2n}.
\end{equation}

Fix
\begin{equation}\label{eq:gamma}
 \gamma=2n-4+\varepsilon',\qquad \varepsilon'>0,
\end{equation}
and let $\HH_{-\gamma}(X,\omega)$ consist of smooth real potentials
with $\omega_u>0$ and
\begin{equation}\label{eq:endpoint-decay}
 |\nabla^j u|_\omega=O(r^{-\gamma-j}),\qquad 0\leq j\leq k.
\end{equation}
The zeroth-order condition fixes the additive constant. A growth
weight $\lambda$ means $|\nabla^ju|=O(r^{\lambda-j})$; weighted
H\"older norms use the corresponding norms on rescaled annuli.
The finite-order condition above does not impose decay of derivatives
of every order. For cohomologous scalar-flat ALE forms approaching
the same Euclidean metric, \cite[Theorem~D(i)]{YaoMass} provides a
global difference potential of order $2-2\widetilde\delta$, with
$\widetilde\delta>n-1$. Taking $\gamma=2\widetilde\delta-2$, we obtain
\eqref{eq:gamma}. Rescaling the potential accounts for the different
$dd^c$ convention in that paper.

Choose a smooth profile $0\leq\chi\leq1$, equal to one on $[0,1]$
and zero on $[2,\infty)$, and set $\chi_R=\chi(r/R)$. On
$A_R=\{R<r<2R\}$, the estimates used below are
\[
 |\nabla^j\chi_R|\leq C_jR^{-j},\qquad
 \operatorname{Vol}_\omega(A_R)\leq CR^{2n},\qquad
 \int_{\{r>R\}}r^{-a}\omega^n\leq C_aR^{2n-a}\quad(a>2n).
\]
The last estimate follows by summing over dyadic annuli. The
constants depend on $n$, $(X,\omega)$, the fixed radius and cutoff
profile, and the indicated derivative order or exponent, but not on
$R$. We will use $\gamma>n-2$, which follows from \eqref{eq:gamma}.

\subsection{Weak geodesics and rooftop envelopes}
An upper semicontinuous function $u$ is $\omega$-psh if
$\omega+dd^cu\geq0$ as a current. In a coordinate ball with
$\omega=dd^c\rho$, this means that $\rho+u$ is plurisubharmonic.
For locally bounded plurisubharmonic potentials, Bedford--Taylor
\cite{BT76,BT} define mixed products inductively by
$dd^cv\wedge T=dd^c(vT)$. The products are positive, closed, symmetric,
and continuous under locally uniform convergence of the potentials.
For $C^{1,1}$ potentials they agree with the almost-everywhere products
of the complex Hessians. The continuity assertion is in terms of measures;
uniform convergence of potentials alone does not imply $L^1$
convergence of their Monge--Amp\`ere densities.

Let $A=\{\tau\in\CC:1<|\tau|<e\}$, $p:X\times A\longrightarrow X$, and $\pi:X\times A\longrightarrow A$.
Write $\Omega=p^*\omega$, $U(x,e^{t+\sqrt{-1}s})=u_t(x)$ for $ t = \log |\tau|$, and
$\Omega_U=\Omega+dd^cU$. A weak geodesic solves
\begin{equation}\label{eq:hcma}
 \begin{cases}
 \Omega_U^{n+1}=0,\qquad &\text{on }X\times A,\\
 \Omega_U \ge 0, \qquad & \text{on } X \times A,\\
 U(x,\tau)=u_j(x)&\text{when }|\tau|=e^j,\quad j=0,1,
 \end{cases}
\end{equation}
in the sense of Bedford--Taylor. The bounded Perron envelope is
\begin{equation}\label{eq:perron}
 U=\sup_{V\in\mathcal B(u_0,u_1)}V.
\end{equation}
Here $\mathcal B(u_0,u_1)$ consists of bounded, $S^1$-invariant
$\Omega$-psh functions whose joint boundary upper limits are at most
$u_j$ on $|\tau|=e^j$. The defining family in
\cite[Section~3]{YaoAsymptotic} has the same boundedness and boundary
conditions, without the $S^1$-invariance requirement.
The endpoint potentials are $\omega$-psh, bounded, and uniformly
continuous by \eqref{eq:endpoint-decay}, and the annulus has strongly
regular boundary. Thus \cite[Proposition~3.3]{YaoAsymptotic} ensures
continuity of the full upper envelope on $X\times\overline A$ and
attainment of the boundary values. Rotations of the annulus preserve
the full defining family, so its envelope is $S^1$-invariant and
belongs to $\mathcal B(u_0,u_1)$. Consequently, restricting to
$\mathcal B(u_0,u_1)$ leaves the supremum unchanged, and no upper
semicontinuous regularization is needed in \eqref{eq:perron}.
By \cite[Theorem~3.4]{YaoAsymptotic}, $U$ is the unique bounded
continuous solution of \eqref{eq:hcma}. It therefore agrees with
the solution of \cite[Theorem~A]{YaoGeo}, whose regularity we record
below.

\begin{proposition}[Uniform $C^{1,1}$ estimate]
\label{prop:joint-regularity}
For $u_0,u_1\in\HH_{-\gamma}(X,\omega)$, the weak geodesic $U$
is unique and continuous on $X\times\overline A$ and has joint $C^{1,1}$
regularity. For $\widehat g=g+dt^2+ds^2$, where $\tau=e^{t+\sqrt{-1}s}$,
\begin{equation}\label{eq:joint-c11-bound}
 \|U\|_{L^\infty}
 +\|\widehat\nabla U\|_{L^\infty}
 +\|\widehat\nabla^2U\|_{L^\infty}\leq C.
\end{equation}
The norms are on $X\times\overline A$, and the Hessian is understood
almost everywhere. The constant depends on $(X,\omega)$, a common
positive lower bound for $\omega_{u_0},\omega_{u_1}$, and a bound
for $\|u_0\|_{C^3(X,g)}+\|u_1\|_{C^3(X,g)}$.
\end{proposition}

In particular, for a constant $\Lambda$ with the same dependencies,
\begin{equation}\label{eq:spatial-bound}
 0\leq\omega_{u_t}\leq\Lambda\omega.
\end{equation}
The estimate is unweighted: it implies neither Hessian decay at infinity
nor continuity of the path in the strong $C^{1,1}$ norm.

Stability follows immediately from the Perron definition. If $U,V$ have
endpoint potentials $u_j,v_j$ and
$a=\max_j\|u_j-v_j\|_\infty$, subtraction of $a$ sends every
function in $\mathcal B(u_0,u_1)$ into $\mathcal B(v_0,v_1)$.
Exchanging the pairs and taking envelopes, we have
\begin{equation}\label{eq:perron-stability}
 \|U-V\|_{L^\infty(X\times\overline A)}
 \leq\max_{j=0,1}\|u_j-v_j\|_{L^\infty(X)}.
\end{equation}

We next introduce rooftop envelopes, which are closely related to weak geodesics.
For bounded continuous real functions $u_0,u_1$, define
\begin{equation}\label{eq:rooftop-definition}
 P_\omega(u_0,u_1)
 =\sup\{v\in\PSH(X,\omega)\cap L^\infty(X):
                          v\leq u_0,\ v\leq u_1\}.
\end{equation}
The functions $u_0,u_1$ are called the obstacles. Their rooftop
envelope is the largest bounded $\omega$-psh function below both.
A sufficiently negative constant makes the defining family nonempty.
The representation below and its Lipschitz dependence on the parameter
are proved in \cite[Lemmas~2.3--2.4]{XiaoYao} for bounded endpoint
potentials on any K\"ahler manifold, using the local minimum principle
of Darvas--Rubinstein \cite[Theorem~3.1]{DarvasRubinstein}.
We also record the local regularity and contact identities
needed in Section~\ref{sec:ale}.

\begin{proposition}[Rooftop representation]
\label{prop:rooftop-representation}
Let $u_0,u_1$ be bounded continuous $\omega$-psh functions, and let
$U$ be their bounded Perron envelope \eqref{eq:perron}.
For $\lambda\in\RR$, put $P_\lambda=P_\omega(u_0,u_1-\lambda)$. Then
\begin{equation}\label{eq:rooftop-representation}
 P_\lambda(x)=\inf_{0\leq s\leq1}\{U(x,s)-\lambda s\},
 \qquad
 U(x,t)=\sup_{\lambda\in\RR}\{P_\lambda(x)+\lambda t\},
\end{equation}
and $\|P_\lambda-P_{\lambda'}\|_\infty\leq|\lambda-\lambda'|$.
\end{proposition}

Next, we record the $C^{1,1}$ regularity of the rooftop envelope, $P_\lambda$.

\begin{lemma}\label{lem:rooftop-c11}
In the setting of Proposition~\ref{prop:rooftop-representation}, suppose
that the slices $u_s=U(\cdot,s)$ satisfy
\[
 \sup_{0\leq s\leq1}\|u_s\|_{C^{1,1}(W')}<\infty
 \qquad\text{for a coordinate ball }W'\Subset X.
\]
Then, for coordinate balls $W\Subset W'\Subset X$,
\begin{equation}\label{eq:rooftop-c11}
 \sup_{\lambda\in\RR}
 \|P_\lambda+\max\{\lambda,0\}\|_{C^{1,1}(W)}\leq C_{W'}.
\end{equation}
The constant depends on $n$, the coefficients of $\omega$ in the
fixed coordinates on $W'$, and the uniform spatial $C^{1,1}(W')$
bound above, but not on $W$ or $\lambda$.
\end{lemma}

\begin{proof}
Fix coordinates on $W'$, and let $C$ bound
$\sup_{s\in[0,1]}\|u_s\|_{C^{1,1}(W')}$.
Fix $x\in W$ and let $s_x\in[0,1]$ attain the infimum in
\eqref{eq:rooftop-representation} such that $P_\lambda(x)=U(x,s_x)-\lambda s_x$. For
$|h|<\operatorname{dist}(W,\partial W')$, we have $x\pm h\in W'$, and
\[
 \begin{aligned}
 P_\lambda(x+h)+P_\lambda(x-h)-2P_\lambda(x)
 &\leq U(x+h,s_x)+U(x-h,s_x)-2U(x,s_x)\\
 &\leq C|h|^2.
 \end{aligned}
\]
Here distance and derivatives are Euclidean in the chosen coordinates.

The second difference estimate means that $D^2P_\lambda\leq CI$
in the sense of distributions. Since $W\Subset W'$ is arbitrary,
the inequality holds throughout $W'$. Write $\omega=dd^c\rho$ on
$W'$. By the plurisubharmonicity of $\rho+P_\lambda$, we have
\[
 \Delta_{\mathrm{euc}}P_\lambda\geq-\Delta_{\mathrm{euc}}\rho
 \geq-C_\omega,
 \qquad C_\omega=\sup_{W'}\Delta_{\mathrm{euc}}\rho,
\]
where $C_\omega$ depends only on the coefficients of $\omega$ in
these coordinates. For any real orthonormal basis
$e_1,\ldots,e_{2n}$, we therefore have, distributionally,
\[
 \begin{aligned}
 \partial_{e_1e_1}P_\lambda
 =\Delta_{\mathrm{euc}}P_\lambda
       -\sum_{a=2}^{2n}\partial_{e_ae_a}P_\lambda \geq-C_\omega-(2n-1)C.
 \end{aligned}
\]
The upper and lower Hessian bounds also hold for local convolutions
with sufficiently small radius. Letting the radius tend to zero
shows that $P_\lambda\in C^{1,1}(W)$, with a Hessian bound depending
only on $n$, $C$, and $C_\omega$.

To bound the first derivatives, take $x,y\in W$ and choose a minimizing
time $s_y$ for $P_\lambda(y)$. Then
\[
 \begin{aligned}
 P_\lambda(x)-P_\lambda(y)
 &\leq U(x,s_y)-\lambda s_y
       -\bigl(U(y,s_y)-\lambda s_y\bigr)\\
 &=U(x,s_y)-U(y,s_y)\leq C|x-y|.
 \end{aligned}
\]
Exchanging $x$ and $y$ proves that $P_\lambda$ is Lipschitz with
constant $C$, and hence
$\|\nabla_{\mathrm{euc}}P_\lambda\|_{L^\infty(W)}\leq C$.
Together with 
\[
 |P_\lambda+\max\{\lambda,0\}|\leq\sup_s|U(\cdot,s)|\leq C,
\]
we obtain \eqref{eq:rooftop-c11}. 
\end{proof}

Next, we discuss the Monge--Amp\`ere measure of the rooftop envelope $P_\lambda$.

\begin{lemma}\label{lem:rooftop-contact-measures}
Let $u_0,u_1$ be bounded real functions in $C^{1,1}_{\mathrm{loc}}(X)$.
For a fixed $\lambda\in\RR$, suppose that
$P_\lambda=P_\omega(u_0,u_1-\lambda)$ belongs to
$C^{1,1}_{\mathrm{loc}}(X)$. Then
\begin{equation}\label{eq:rooftop-contact-measures}
 \begin{aligned}
 &\mathbf1_{\{P_\lambda<\min(u_0,u_1-\lambda)\}}
                    \omega_{P_\lambda}^n=0,\\
 &\mathbf1_{\{P_\lambda=u_j-j\lambda\}}\omega_{P_\lambda}^n
 =\mathbf1_{\{P_\lambda=u_j-j\lambda\}}\omega_{u_j}^n, \quad a.e.
                    \qquad j=0,1.
 \end{aligned}
\end{equation}
\end{lemma}

\begin{proof}
We first prove that the Monge--Amp\`ere measure vanishes where
$P_\lambda<\min(u_0,u_1-\lambda)$. Fix a coordinate ball
$B\Subset\{P_\lambda<\min(u_0,u_1-\lambda)\}$. By continuity, we have
\[
 c:=\min_{\overline B}
 \bigl(\min(u_0,u_1-\lambda)-P_\lambda\bigr)>0.
\]
To prove local maximality, let $v$ be a bounded $\omega$-psh
function on $B$ with
\[
 \limsup_{x\to y,\,x\in B}v(x)\leq P_\lambda(y),
 \qquad y\in\partial B.
\]
We claim that $v\leq P_\lambda$ on $B$. Otherwise,
$M:=\sup_B(v-P_\lambda)_+>0$. Choose $0<\varepsilon<1$ with
$\varepsilon M<c$. The function
$(1-\varepsilon)P_\lambda+\varepsilon v$ is $\omega$-psh, since
\[
 \omega+dd^c\bigl((1-\varepsilon)P_\lambda+\varepsilon v\bigr)
 =(1-\varepsilon)\omega_{P_\lambda}+\varepsilon\omega_v\geq0.
\]
Moreover, on $B$,
\[
 (1-\varepsilon)P_\lambda+\varepsilon v
 \leq P_\lambda+\varepsilon M
 <P_\lambda+c\leq\min(u_0,u_1-\lambda).
\]
The boundary inequality permits us to glue by a maximum:
\[
 \widetilde P=
 \begin{cases}
 \max\{P_\lambda,(1-\varepsilon)P_\lambda+\varepsilon v\},
       &\text{on }B,\\
 P_\lambda,&\text{on }X\setminus B.
 \end{cases}
\]
The resulting function is bounded and $\omega$-psh on $X$ and
lies below both obstacles. Hence $\widetilde P\leq P_\lambda$
by the definition of the rooftop envelope. At any point where
$v>P_\lambda$, however, we have by definition
$\widetilde P>P_\lambda$, a contradiction. Thus $P_\lambda$
is locally maximal. Writing $\omega=dd^c\rho$ in the coordinate
ball, $\rho+P_\lambda$ is a locally bounded maximal
plurisubharmonic function. By the Bedford--Taylor maximality criterion
\cite{BT76,BT}, we have
\[
 \omega_{P_\lambda}^n=(dd^c(\rho+P_\lambda))^n=0
 \quad\text{on }B.
\]
Since $B$ was arbitrary, the first identity follows.

We next prove the Hessian agreement on the contact sets. We use
the elementary fact that a locally Lipschitz function $F$ has
zero derivative almost everywhere on $\{F=0\}$. In one dimension,
$F$ is differentiable almost everywhere. At a nonisolated zero
$s$ where $F'(s)$ exists, choose distinct zeros $s_k\to s$; then
\[
 F'(s)=\lim_{k\to\infty}\frac{F(s_k)-F(s)}{s_k-s}=0.
\]
The isolated zeros form a countable set. Applying the
one-dimensional conclusion on lines parallel to each coordinate
axis and using Fubini's theorem proves the assertion in higher
dimensions; the derivatives along those lines agree almost
everywhere with the weak partial derivatives.

For $j=0,1$, set $w=P_\lambda-u_j+j\lambda\in C^{1,1}_{\mathrm{loc}}(X)$.
Apply the preceding fact first
to $w$. Each $\partial_{x_a}w$ vanishes almost everywhere on
$\{w=0\}$. Since
$\partial_{x_a}w$ is locally Lipschitz, we may apply the same fact again to obtain
\[
 \partial_{x_b}\partial_{x_a}w=0
 \quad\text{almost everywhere on }\{w=0\},
 \qquad a,b=1,\ldots,2n.
\]
Thus two $C^{1,1}$ functions agreeing on any measurable set have
equal Hessians almost everywhere on that set. In particular,
\[
 dd^cP_\lambda=dd^cu_j
 \quad\text{almost everywhere on }
 \{P_\lambda=u_j-j\lambda\},\qquad j=0,1.
\]
For $C^{1,1}$ potentials the wedge products have locally bounded
determinant densities, which agree with the Bedford--Taylor
measures when the potentials are $\omega$-psh. Equality of the
complex Hessians therefore implies equality of the restricted
measures, proving the second identity in
\eqref{eq:rooftop-contact-measures}.
\end{proof}

Bedford--Taylor continuity \cite{BT76,BT} ensures convergence of
Monge--Amp\`ere measures against continuous compactly supported test
functions when the potentials converge locally uniformly. Convergence
of the masses on arbitrary fixed measurable sets is stronger and does
not hold in general. The following lemma obtains that conclusion
under common local $C^{1,1}$ bounds; we will apply it to the contact
sets in Section~\ref{subsec:compact-density-convergence}.

\begin{lemma}\label{lem:local-density-weak}
Let $v_k,v$ be $\omega$-psh functions in
$C^{1,1}_{\mathrm{loc}}(X)$. Suppose that $v_k\to v$
locally uniformly and that the $v_k$ have common local
$C^{1,1}$ bounds. Then, for every fixed relatively compact
measurable set $E\subset X$,
\[
 \int_E\omega_{v_k}^n\longrightarrow\int_E\omega_v^n.
\]
\end{lemma}

\begin{proof}
Fix a relatively compact measurable set $E$ and choose a relatively
compact open set $V$ containing $\overline E$. Write
$f_k=\omega_{v_k}^n/\omega^n$ and $f=\omega_v^n/\omega^n$.
By the common local $C^{1,1}$ bounds, we have $0\leq f_k,f\leq C$ on $V$,
where $C$ depends only on $n$, the reference metric on $\overline V$,
and the common $C^{1,1}$ bound there, independently of $k$.
For $\varphi\in C^\infty_{\mathrm c}(V)$, integration by parts gives
\[
 \int_V\varphi(f_k-f)\,\omega^n
 =\int_V(v_k-v)\ddc\varphi\wedge
       \sum_{j=0}^{n-1}\omega_{v_k}^j\wedge\omega_v^{n-1-j}
 \longrightarrow0.
\]
Indeed, the common Hessian bound controls the mixed forms on
$\supp\varphi$, while $v_k-v\to0$ uniformly there.

To pass to the set $E$, use
\[
 \left|\int_E(f_k-f)\,\omega^n\right|
 \leq\left|\int_V\varphi(f_k-f)\,\omega^n\right|
       +2C\|\mathbf1_E-\varphi\|_{L^1(V,\omega^n)}.
\]
First let $k\to\infty$ with $\varphi$ fixed. Since
$\mathbf1_E\in L^1(V,\omega^n)$ can be approximated by functions
in $C^\infty_{\mathrm c}(V)$, the last term can be made arbitrarily
small. This proves the asserted convergence.
\end{proof}

\subsection{The relative Mabuchi functional and its variations}
For a smooth closed real $(1,1)$-form $\eta$, define
\begin{equation}\label{eq:twisted-energy}
 E_\eta(u)=\frac1n\sum_{j=0}^{n-1}
 \int_Xu\,\eta\wedge\omega_u^{n-1-j}\wedge\omega^j,
\end{equation}
on the subset of the space of potentials for which the integrals converge.
Set $\vartheta=\Ric(\omega)$, $ f_u=\frac{\omega_u^n}{\omega^n}$, and $h(s)=s\log s-s+1$.
We define the relative entropy of $u$ with respect to $\omega$ by
\begin{equation}\label{eq:entropy}
 \Ent_{\mathrm{rel}}(u)=\int_Xh(f_u)\,\omega^n.
\end{equation}
The relative Mabuchi functional is
\begin{equation}\label{eq:mabuchi}
 \MM(u)=\Ent_{\mathrm{rel}}(u)-nE_\vartheta(u).
\end{equation}
In proofs we may write $H(u)=\Ent_{\mathrm{rel}}(u)$ and
$\mathcal E(u)=nE_\vartheta(u)$. When the reference K\"ahler form
must be specified, we use $\Ent_{\mathrm{rel},\eta}$, or $H_\eta$
in proofs.

The subtraction in $h$ removes the linear part of the entropy
density at $f=1$, which need not be integrable on an ALE end.
For $u\in\HH_{-\gamma}(X,\omega)$,
\begin{equation}\label{eq:smooth-entropy-tail}
 f_u-1=O(r^{-\gamma-2}),\qquad h(f_u)=O(r^{-2\gamma-4}),
\end{equation}
so $\Ent_{\mathrm{rel}}(u)<\infty$. Convergence of $E_\vartheta$
for weak geodesic slices is proved in
Section~\ref{sec:energy-finiteness}.

We next establish lower semicontinuity of the relative entropy on the
whole manifold under local weak convergence of the Monge--Amp\`ere
measures. The result holds on a general noncompact K\"ahler manifold.
We will apply it to weak geodesics with compactly supported endpoint
potentials whose supports increase to $X$.

\begin{lemma}
\label{lem:entropy-lsc}
Let $(X,\omega)$ be a noncompact K\"ahler manifold of complex
dimension $n$. Let $u_k,u$ be locally bounded $\omega$-psh potentials
whose Monge--Amp\`ere measures have locally integrable densities
with respect to $\omega^n$. Suppose that
\[
 \int_X\varphi\,\omega_{u_k}^n
 \longrightarrow\int_X\varphi\,\omega_u^n
 \qquad\text{for every }\varphi\in C_c(X,\RR).
\]
Then
\begin{equation}\label{eq:entropy-lsc}
 \Ent_{\mathrm{rel}}(u)\leq\liminf_{k\to\infty}
                         \Ent_{\mathrm{rel}}(u_k).
\end{equation}
The entropy integrals are taken over the whole manifold $X$ and
may equal $+\infty$.
\end{lemma}

\begin{proof}
For every $s\geq0$, by the Legendre transform of $h$, we have
\[
 h(s)=\sup_{b\in\RR}\{bs-e^b+1\}.
\]
For $s>0$, the supremum is attained at $b=\log s$; for $s=0$,
it is approached as $b\to-\infty$. Consequently,
\[
 \sup_{a\in C_c(X,\RR)}
 \left\{\int_Xa\,\omega_u^n-\int_X(e^a-1)\,\omega^n\right\}
 \leq\Ent_{\mathrm{rel}}(u).
\]
Both integrals in braces are finite because $a$ has compact
support and $f_u$ is locally integrable. To prove equality, we
approximate the pointwise maximizing function $\log f_u$ by
continuous compactly supported functions.

The reverse inequality is new in the ALE setting. Choose a compact exhaustion
$K_m\subset\operatorname{int}K_{m+1}$, $\bigcup_mK_m=X$, and fix
$m,N\geq1$. Since $\log f_u$ may be unbounded and is $-\infty$
where $f_u=0$, set
\[
 a_N=\max\{-N,\min\{N,\log f_u\}\},\qquad
 a=\mathbf1_{K_m}a_N,
\]
with $a_N=-N$ where $f_u=0$. The function $a$ is bounded by $N$
and supported in $K_m$, and
\[
 \int_X\bigl(af_u-e^a+1\bigr)\,\omega^n
 =\int_{K_m}\bigl(a_Nf_u-e^{a_N}+1\bigr)\,\omega^n,
\]
because the integrand vanishes outside $K_m$. However, $a$ is
only measurable in general, so it cannot be used directly in the
supremum over continuous compactly supported functions. We therefore
approximate $a$ by functions
$a_j\in C_c(\operatorname{int}K_{m+1},\RR)$, with $|a_j|\leq N$,
such that
\[
 \int_{K_{m+1}}|a_j-a|(f_u+1)\,\omega^n\longrightarrow0.
\]
The approximation is possible because $(f_u+1)\omega^n$ is a finite
Radon measure on $K_{m+1}$. Indeed, approximate $a$ by simple
functions and use inner and outer regularity to approximate their
indicator functions by continuous cutoffs supported in
$\operatorname{int}K_{m+1}$.

The weight $f_u+1$ controls both the linear and exponential terms.
Since $|a_j|,|a|\leq N$, we have by the mean value theorem
$|e^{a_j}-e^a|\leq e^N|a_j-a|$. Hence
\[
 \begin{aligned}
 &\left|\int_X\bigl(a_jf_u-e^{a_j}+1\bigr)\,\omega^n
       -\int_{K_m}\bigl(a_Nf_u-e^{a_N}+1\bigr)\,\omega^n\right|\\
 &\qquad\leq\int_{K_{m+1}}f_u|a_j-a|\,\omega^n
          +\int_{K_{m+1}}|e^{a_j}-e^a|\,\omega^n\\
 &\qquad\leq\int_{K_{m+1}}(f_u+e^N)|a_j-a|\,\omega^n
 \longrightarrow0.
 \end{aligned}
\]
Every $a_j$ is an admissible test function in the supremum, so
\[
 \sup_{a\in C_c(X,\RR)}
 \left\{\int_Xa\,\omega_u^n-\int_X(e^a-1)\,\omega^n\right\}
 \geq\int_{K_m}\bigl(a_Nf_u-e^{a_N}+1\bigr)\,\omega^n.
\]
At each point, $a_N$ maximizes $bf_u-e^b+1$ over $|b|\leq N$.
The allowed intervals increase to $\RR$ and contain $b=0$, so
\[
 0\leq a_Nf_u-e^{a_N}+1\uparrow h(f_u).
\]
First let $N\to\infty$, and then $m\to\infty$. Monotone
convergence implies the opposite inequality; hence,
\begin{equation}\label{eq:entropy-dual}
 \Ent_{\mathrm{rel}}(u)=\sup_{a\in C_c(X,\RR)}
 \left\{\int_Xa\,\omega_u^n-\int_X(e^a-1)\,\omega^n\right\},
\end{equation}
including the case of infinite entropy.

Now fix $a\in C_c(X,\RR)$. By local weak convergence, we have
\[
 \begin{aligned}
 \int_Xa\,\omega_u^n-\int_X(e^a-1)\,\omega^n
 &=\lim_{k\to\infty}
   \left(\int_Xa\,\omega_{u_k}^n-\int_X(e^a-1)\,\omega^n\right)\\
 &\leq\liminf_{k\to\infty}\Ent_{\mathrm{rel}}(u_k).
 \end{aligned}
\]
Taking the supremum over $a$ proves \eqref{eq:entropy-lsc}.
\end{proof}

In particular, suppose that $u_{t,l}\to u_t$ locally uniformly for
a fixed time $t$, with common local $C^{1,1}$ bounds, and that each
$u_{t,l}$ is compactly supported. Lemma~\ref{lem:local-density-weak}
ensures the required local convergence of volume measures, so
\[
 \Ent_{\mathrm{rel}}(u_t)
 \leq\liminf_{l\to\infty}\Ent_{\mathrm{rel}}(u_{t,l}).
\]
The supports may grow without bound. Since $h\geq0$, a uniform bound
for the approximating entropies implies finiteness of the limiting
entropy. Proving equality requires additional estimates and is the
subject of Section~\ref{sec:ale}.

For a smooth path $u_s\in\HH_{-\gamma}$ with
$v_s=\partial_su_s$ supported in a common compact set,
differentiation and integration by parts give
\[
 \begin{aligned}
 \frac{d}{ds}\Ent_{\mathrm{rel}}(u_s)
 &=n\int_X\log\frac{\omega_{u_s}^n}{\omega^n}\,
                      dd^cv_s\wedge\omega_{u_s}^{n-1}\\
 &=n\int_Xv_s\bigl(\vartheta-\Ric(\omega_{u_s})\bigr)
                                  \wedge\omega_{u_s}^{n-1},\\
 \frac{d}{ds}E_\vartheta(u_s)
 &=\int_Xv_s\vartheta\wedge\omega_{u_s}^{n-1}.
 \end{aligned}
\]
Subtracting yields $d\MM|_u(v)=-\int_XvS_{\omega_u}\,\omega_u^n$, for $ v\in C_c^\infty(X).$
Thus the critical metrics are scalar-flat. Write
$\mathcal D_{\omega_u}v=\bar\partial\nabla^{1,0}_{\omega_u}v$ and
$L_{\omega_u}=\mathcal D_{\omega_u}^*\mathcal D_{\omega_u}$ for the
Lichnerowicz operator. Along a smooth geodesic with variations
supported in a common compact set, we have by the second-variation formula
\cite{Mabuchi86,Mabuchi87}
$\frac{d^2}{dt^2}\MM(u_t)
=\int_X|\mathcal D_{\omega_{u_t}}\dot u_t|^2\,\omega_{u_t}^n\geq0$.
For a weak geodesic the spatial form may degenerate, so we use
the following current identities.

We call $X\times\{\tau\}$ a slice of $X\times A$. For a current
$T$ on $X\times A$, integration over the slices of $\pi$ is denoted
by $\pi_*T$ and defined by
\[
 \langle\pi_*T,\zeta\rangle=\langle T,\pi^*\zeta\rangle.
\]
Here $\zeta$ is a compactly supported test form of complementary
degree on $A$. We require $\supp T\cap\pi^{-1}(K)$ to be compact
for every compact $K\subset A$; support in $K_X\times A$ with
$K_X\Subset X$ suffices. This condition makes the pairing well-defined
although the slices are noncompact. Integration removes spatial
bidegree $(n,n)$.
For example, $\pi_*(a\,p^*\omega^n)(\tau)=\int_Xa(x,\tau)\omega^n$
when $a$ is supported in $K_X \times A$.

We write $dd^c_\tau$ for $dd^c$ acting in the annulus variable
$\tau$. Suppose that $U$ has joint $C^{1,1}$ regularity and common
compact spatial support. The support condition allows $dd^c$ to
commute with integration over the slices, so by differentiating
\eqref{eq:twisted-energy} we obtain
\begin{equation}\label{eq:twisted-energy-current}
 \begin{aligned}
 n\,dd^c_\tau E_\eta(u_\tau)
 &=\pi_*dd^c\left(
 U\,p^*\eta\wedge\sum_{j=0}^{n-1}
                      \Omega_U^{n-1-j}\wedge\Omega^j\right)\\
 &=\pi_*\bigl(p^*\eta\wedge(\Omega_U^n-\Omega^n)\bigr)
 =\pi_*\bigl(p^*\eta\wedge\Omega_U^n\bigr).
 \end{aligned}
\end{equation}
The term $p^*\eta\wedge\Omega^n$ vanishes by spatial degree.
Outside the support of $U$, we have $\Omega_U=\Omega$, so
$p^*\eta\wedge\Omega_U^n$ also vanishes there.
For $U\in C^{1,1}$, the calculation is understood in the sense of
currents.

\subsection{Entropy and the complex Hessian}\label{subsec:entropy-hessian}
We first prove elementary pointwise inequalities relating the entropy
density to the complex Hessian under $0\leq\omega_u\leq\Lambda\omega$,
allowing degeneracy. The estimate controlling the Hessian includes the
difference between the trace and determinant variations, whose integral
vanishes for compactly supported potentials. Integrating, we obtain
$L^2$ Hessian bounds without assuming pointwise Hessian decay. We use
the bounds to control the errors from spatial cutoffs in the second
derivative of the Mabuchi functional in Section~\ref{sec:ale} and to prove absolute convergence
of the mixed integrals defining $E_\vartheta$ in Section~\ref{sec:ale-finiteness}.

\begin{lemma}\label{lem:entropy-hessian}
Let $u\in C^{1,1}_{\mathrm{loc}}(X,\mathbb R)$ satisfy
$0\leq\omega_u\leq\Lambda\omega$. Write
$A_u=\omega^{-1}\omega_u$, $f_u=\det A_u$, and
$q_u=\tr(A_u-I)-(f_u-1)$, and use the Hilbert--Schmidt norm
in an $\omega$-unitary frame. Then, almost everywhere on $X$,
\begin{equation}\label{eq:entropy-hessian-pointwise}
 \begin{aligned}
 h(f_u)&\leq C_{n,\Lambda}|A_u-I|^2,\\
 |A_u-I|^2&\leq2(\Lambda+1)q_u
             +(\Lambda+1)(4n-2)h(f_u),
 \end{aligned}
\end{equation}
where $C_{n,\Lambda}=n\max\{1,\Lambda\}^{2n-2}$.
The pointwise inequalities have the following integral consequences.
If $u$ has compact support, then
\begin{equation}\label{eq:entropy-hessian-compact}
 \int_X|A_u-I|^2\,\omega^n
 \leq(\Lambda+1)(4n-2)\Ent_{\mathrm{rel}}(u).
\end{equation}
If $|u|\leq C_{\mathrm w}r^{-\gamma}$, with $\gamma$ as in
\eqref{eq:gamma}, and $\Ent_{\mathrm{rel}}(u)<\infty$, then
\begin{equation}\label{eq:hessian-l2-tail}
 \int_{\{r>2R\}}|A_u-I|^2\,\omega^n
 \leq C\int_{\{r>R\}}h(f_u)\,\omega^n
       +CR^{2n-2\gamma-4}\longrightarrow0
 \quad\text{as }R\longrightarrow\infty.
\end{equation}
Here $C$ and the lower threshold for $R$ depend only on
$n$, $(X,\omega)$, $r$, the cutoff profile, $\gamma$,
$C_{\mathrm w}$, and $\Lambda$.
\end{lemma}

\begin{proof}
We first prove the pointwise estimates for a Hermitian matrix
$0\leq A\leq\Lambda I$ with eigenvalues $\lambda_1,\ldots,\lambda_n$.
Put $f=\det A$ and $q=\tr(A-I)-(f-1)$.
Since $(\sqrt{\lambda_j}+1)^2\leq2(\lambda_j+1)\leq2(\Lambda+1)$,
factorization and the arithmetic--geometric mean inequality give
\[
 \begin{aligned}
 |A-I|^2
 &=\sum_{j=1}^n(\sqrt{\lambda_j}-1)^2(\sqrt{\lambda_j}+1)^2 \leq2(\Lambda+1)\sum_{j=1}^n(\sqrt{\lambda_j}-1)^2\\
 &\leq2(\Lambda+1)
       \bigl[q+f-1-2n(f^{1/(2n)}-1)\bigr].
 \end{aligned}
\]
For $f>0$, applying $e^x\geq1+x$ to
$x=-(2n-1)\log f/(2n)$ and multiplying by $f$ yields
\[
 f^{1/(2n)}\geq f-\frac{2n-1}{2n}f\log f.
\]
Consequently,
\[
 f-1-2n(f^{1/(2n)}-1)
 \leq(2n-1)(f\log f-f+1).
\]
By continuity the inequality also holds at $f=0$, and therefore
\begin{equation}\label{eq:matrix-cost}
 |A-I|^2\leq2(\Lambda+1)q+(\Lambda+1)(4n-2)h(f).
\end{equation}
In the other direction, $\log f\leq f-1$ implies
$h(f)\leq(f-1)^2$, including $f=0$ by continuity.
Expanding $f-1=\sum_{j=1}^n(\lambda_j-1)\prod_{a<j}\lambda_a$, we have
\begin{equation}\label{eq:ale-determinant-bound}
 |f-1|\leq\max\{1,\Lambda\}^{n-1}
                  \sum_{j=1}^n|\lambda_j-1|.
\end{equation}
Cauchy--Schwarz then yields
\begin{equation}\label{eq:bounded-matrix-cost}
 h(f)\leq(f-1)^2\leq C_{n,\Lambda}|A-I|^2,
 \qquad C_{n,\Lambda}=n\max\{1,\Lambda\}^{2n-2}.
\end{equation}

We next identify the cancellation in $q_u=\tr(A_u-I)-(f_u-1)$.
Define the closed $(n-1,n-1)$-form with locally bounded coefficients
\[
 T_u=n\omega^{n-1}
       -\sum_{j=0}^{n-1}\omega_u^j\wedge\omega^{n-1-j}.
\]
By the trace and determinant identities, we have
\[
 \begin{aligned}
 \tr(A_u-I)\,\omega^n&=n\,dd^cu\wedge\omega^{n-1},\\
 (f_u-1)\,\omega^n
 &=dd^cu\wedge\sum_{j=0}^{n-1}\omega_u^j\wedge\omega^{n-1-j},
 \end{aligned}
\]
so that $q_u\,\omega^n=dd^cu\wedge T_u=dd^c(uT_u)$.
If $u$ has compact support, integration of the last identity gives
$\int_Xq_u\,\omega^n=0$. Integrating
\eqref{eq:matrix-cost}
proves \eqref{eq:entropy-hessian-compact}.

For a decaying potential, the same cancellation improves the error
from integration by parts. More precisely,
\[
 \begin{aligned}
 T_u
 &=-\sum_{j=1}^{n-1}(\omega_u^j-\omega^j)
                      \wedge\omega^{n-1-j}\\
 &=-dd^cu\wedge
        \sum_{j=0}^{n-2}(n-1-j)\omega_u^j\wedge\omega^{n-2-j},
 \end{aligned}
\]
and $0\leq\omega_u\leq\Lambda\omega$ therefore implies
$|T_u|_\omega\leq C_{n,\Lambda}|A_u-I|$.
Let $\theta_R=1-\chi_R$, and for $L>4R$ take
$\zeta=(\theta_R\chi_L)^2$. The function $\zeta$ is compactly
supported, equals one on $\{2R<r<L\}$, and satisfies
\[
 |dd^c\zeta|_\omega
 \leq C\zeta^{1/2}
       \left(R^{-2}\mathbf1_{\{R<r<2R\}}
             +L^{-2}\mathbf1_{\{L<r<2L\}}\right).
\]
We choose the cutoff profile $\chi$ so that
$|\chi'|^2\leq C\min\{\chi,1-\chi\}$.
On the two transition annuli, $\zeta=\theta_R^2$ and
$\zeta=\chi_L^2$, respectively, so the bound follows by
differentiating the squares.

Testing $q_u\omega^n=dd^c(uT_u)$ against $\zeta$, using the
bound for $T_u$, and applying Young's inequality, we obtain for
every $\eta>0$
\[
 \begin{aligned}
 \left|\int_X\zeta q_u\,\omega^n\right|
 &=\left|\int_Xu\,dd^c\zeta\wedge T_u\right|\\
 &\leq C\sum_{a\in\{R,L\}}a^{-2}
       \int_{\{a<r<2a\}}|u|\zeta^{1/2}|A_u-I|\,\omega^n\\
 &\leq\eta\int_X\zeta|A_u-I|^2\,\omega^n
       +C_\eta\sum_{a\in\{R,L\}}a^{-4}
                      \int_{\{a<r<2a\}}|u|^2\,\omega^n\\
 &\leq\eta\int_X\zeta|A_u-I|^2\,\omega^n
       +C_\eta\bigl(R^{2n-2\gamma-4}+L^{2n-2\gamma-4}\bigr).
 \end{aligned}
\]
The final line uses $|u|\leq C_{\mathrm w}r^{-\gamma}$ and
the volume estimate for annuli. On the other hand,
\eqref{eq:matrix-cost} gives
\[
 \begin{aligned}
 \int_X\zeta|A_u-I|^2\,\omega^n
 &\leq2(\Lambda+1)\int_X\zeta q_u\,\omega^n
       +(\Lambda+1)(4n-2)\int_X\zeta h(f_u)\,\omega^n.
 \end{aligned}
\]
Choose $\eta$ so that $2(\Lambda+1)\eta\leq1/2$, absorb the
corresponding term, and use $0\leq\zeta\leq\mathbf1_{\{r>R\}}$.
It follows that
\[
 \int_X\zeta|A_u-I|^2\,\omega^n
 \leq C\int_{\{r>R\}}h(f_u)\,\omega^n
       +C\bigl(R^{2n-2\gamma-4}+L^{2n-2\gamma-4}\bigr).
\]
Since $\gamma>n-2$, the power $2n-2\gamma-4$ is negative.
Letting $L\to\infty$ and applying Fatou's lemma proves
\eqref{eq:hessian-l2-tail}. The entropy integral over $\{r>R\}$
tends to zero by its assumed finiteness. In particular,
$A_u-I\in L^2(X,\omega^n)$, since its coefficients are bounded
on each compact set.
\end{proof}

% ----------------------------------------------------------------------
% 3. Weak geodesics with compact support and convexity
% ----------------------------------------------------------------------
\section{Weak geodesics with compact support and convexity}
\label{sec:compact}

In this section, we study weak geodesics joining compactly supported
K\"ahler potentials. These are described by the homogeneous complex
Monge--Amp\`ere equation with the prescribed potentials as boundary data.
Throughout the section, we write
\[
  \HH_{\mathrm c}(X,\omega)
  =
  \left\{
    u\in C^\infty_{\mathrm c}(X,\mathbb R)
    \;\middle|\;
    \omega_u:=\omega+\ddc u>0
  \right\}.
\]
Let $u_0,u_1\in\HH_{\mathrm c}(X,\omega)$, and set
$A=\{\tau\in\CC:1<|\tau|<e\}$. Denote by
$p:X\times A\to X$ the projection and put $\Omega=p^*\omega$.
The weak geodesic joining $u_0$ to $u_1$ is the path
$u_t(x)=U(x,e^t)$, $0\leq t\leq1$, where $U$ is a bounded
$S^1$-invariant solution of the Dirichlet problem
\[
  \begin{cases}
    \Omega_U^{n+1}=0,\quad
      &\text{on }X\times A,\\
      \Omega_U\geq0, & \text{on } X \times A,\\
    U(x,\tau)=u_j(x)
      &\text{when }|\tau|=e^j,\quad j=0,1,
  \end{cases}
\]
where $\Omega_U=\Omega+\ddc U$.
The equation is understood in the sense of Bedford--Taylor.
The solution $U$ is given by the Perron envelope
\eqref{eq:perron} described in Section~\ref{sec:preliminary}.
By the existence and regularity result in
\cite[Theorem~A]{YaoGeo}, recalled there, we have continuity up to
the boundary and joint $C^{1,1}$ regularity.

We first prove that the potentials $u_t$ are supported in a
fixed compact subset of $X$, independent of $t$. We then
establish convexity of the relative Mabuchi functional along
these weak geodesics using the local Bergman kernel argument.
These compactly supported weak geodesics provide the starting
point for the approximation argument in Section~\ref{sec:ale},
where we treat weak geodesics with general decaying endpoint
potentials.

\subsection{A Liouville-type lemma}

We first prove that, for compactly supported boundary data, the weak
geodesic has slices supported in a common compact set. This is a
Liouville-type statement: the solution vanishes identically on a
sufficiently distant part of the end.

The local regularity theorem in \cite[Theorem~A]{YaoLocal}
relates the solution of the HCMA equation to the boundary data on a
spatial neighbourhood. Its geometric ingredient is a local Monge--Amp\`ere
foliation. The total K\"ahler potential is obtained
by solving a family of Dirichlet problems on each Monge--Amp\`ere leaf
\cite[Proposition~3.3]{YaoLocal}. A local convexity condition on the
reference potential, together with a sufficiently large spatial
buffer, identifies this local solution with the global Perron
envelope \cite[Theorem~6.5]{YaoLocal}. 

Theorem~A of \cite{YaoLocal} is stated over a disc. The disc
hypothesis is significant here: over an annulus, the leaves of a
local Monge--Amp\`ere foliation may remain open around the
parameter domain, so a foliation by annular graphs cannot be assumed.
Moreover, the examples in \cite[Theorem~B]{XiaoYao} exhibit weak geodesics
whose spatial forms degenerate on open sets arbitrarily far out,
despite rapidly decaying endpoint potentials in the Euclidean case.
Such degeneracy obstructs the regular Monge--Amp\`ere foliation associated
with a spatially nondegenerate solution. Thus the local regularity
argument for the disc requires additional justification in the
annular setting.

When the boundary data vanish on a spatial neighbourhood, they
admit the trivial Monge--Amp\`ere foliation by the vertical annuli
$\{x\}\times\overline A$. In this case the leaves are closed and trivial over
the parameter annulus on $(X-K)\times\partial A$. The
leafwise Dirichlet problems are therefore well-defined. We give the
local comparison argument first, followed by a proof using radial
barriers.

\begin{lemma}[Liouville-type compact-support lemma]
\label{lem:support}
Let $u_0,u_1\in\HH_{\mathrm c}(X,\omega)$, and let $U$ be the
weak geodesic defined above. There exists a compact set $K\Subset X$
such that
\[
  U(x,\tau)=0
  \qquad\text{for }x\in X\setminus K,\quad \tau\in\overline A.
\]
In particular, $\supp u_t\subset K$ for every $t\in[0,1]$.
\end{lemma}

\begin{proof}[First proof: the local foliation construction]
Choose
\[
  M>1+\max\{\|u_0\|_{L^\infty(X)},\|u_1\|_{L^\infty(X)}\}.
\]
The constant function $-M$ belongs to $\mathcal B(u_0,u_1)$,
since it is bounded, $\Omega$-psh, and lies below both boundary
values. Hence $-M\leq U$ by \eqref{eq:perron}. On the other hand,
restriction to a vertical annulus shows that $u_t(x)$ is convex
in $t$, and hence
\[
  -M\leq U(x,e^t)\leq(1-t)u_0(x)+tu_1(x)<M.
\]
Thus the global bound required in the local comparison is available.
The joint $C^{1,1}$ regularity recalled above also implies continuity up
to the temporal boundary; no higher regularity of $U$ is assumed in
this argument.

Fix $L$ with $L^2>16M$. By \cite[Lemma~2.10]{YaoAsymptotic},
every sufficiently distant point of the ALE end admits a holomorphic
coordinate ball $B_{3L}\subset\CC^n$, disjoint from
$\supp u_0\cup\supp u_1$. On each such ball, by the local potential
estimates of \cite[Propositions~8.2 and~8.4(i)]{YaoLocal}, we have
$\omega=\ddc\rho$ with
\begin{equation}\label{eq:local-support-convexity}
  D^2\rho\geq\lambda I_{2n},\qquad
  \lambda=\tfrac12,\qquad \lambda L^2>8M.
\end{equation}
Here $D^2$ denotes the real Euclidean Hessian in these holomorphic
coordinates. More generally, the local convexity threshold
$\lambda L^2$ can be made arbitrarily large on the ALE end by
increasing $L$ and choosing the centres sufficiently far out.
We now work on one such ball.

Since both endpoint potentials vanish on this ball, the boundary
graphs in the local Semmes cotangent bundle are both
$\{(z,\partial\rho(z)):z\in B_{3L}\}$. For each $y\in B_L$, the
map
\[
  G_y(\tau)=(y,\partial\rho(y)),\qquad \tau\in\overline A,
\]
is holomorphic in $\tau$ and satisfies these boundary conditions. Its
projected leaf in $X\times\overline A$ is
$\{y\}\times\overline A$. These leaves fill $B_L\times\overline A$;
each is defined over the entire closed annulus and has boundary only
on $X\times\partial A$. They are the kernel leaves of $\Omega$,
which has positive spatial restriction $\omega$.

The reconstruction from this foliation consists of the family of
Dirichlet problems
\[
  \Delta_A h_y=0\quad\text{in }A,\qquad
  h_y=\rho(y)\quad\text{on }\partial A,
  \qquad y\in B_L.
\]
By the maximum principle, we have $h_y\equiv\rho(y)$. Subtracting the
reference potential therefore gives the local HCMA solution $V=0$.
It remains to show that this local solution agrees with $U$; local
Dirichlet uniqueness alone does not establish this identification because
no values have been prescribed on $\partial B_{3L}\times A$.

For $y\in\overline{B_L}$, consider the real affine function
\[
  L_y(z)=\rho(y)+2\operatorname{Re}
  \sum_{a=1}^n\frac{\partial\rho}{\partial z_a}(y)(z_a-y_a).
\]
It is pluriharmonic and agrees with $\rho$ to first order at $y$.
By Taylor's formula and \eqref{eq:local-support-convexity}, we have
\begin{equation}\label{eq:local-support-tangent}
  L_y(z)-\rho(z)\leq-\frac{\lambda}{2}|z-y|^2,
  \qquad z\in B_{3L},\quad y\in\overline{B_L}.
\end{equation}
Set
\[
  H(z)=\sup_{y\in\overline{B_L}}\bigl(L_y(z)-\rho(z)\bigr).
\]
The family is continuous in $(y,z)$ and the parameter set is compact,
so $H$ is continuous. Moreover, $\rho+H$ is plurisubharmonic, being
a locally bounded, continuous supremum of pluriharmonic functions.
The inequality above shows that $H\leq0$, while choosing $y=z$
shows that $H=0$ on $B_L$. For $2L\leq|z|<3L$, we also have
\[
  H(z)\leq-\frac{\lambda L^2}{2}<-M.
\]
Consequently the function
\[
  F(x)=
  \begin{cases}
    \max\{-M,H(z(x))\},&x\text{ in the coordinate ball},\\
    -M,&x\text{ outside the coordinate ball},
  \end{cases}
\]
is well-defined and continuous across the boundary of the ball.
It is $\omega$-psh: inside the ball,
$\rho+F=\max\{\rho-M,\sup_y L_y\}$ is plurisubharmonic,
and $F=-M$ on a collar of the boundary. Viewed as a function on
$X\times A$, it is a bounded $\Omega$-psh function.

Since $F\leq u_j$ for $j=0,1$, we have $F\in\mathcal B(u_0,u_1)$,
and hence $U\geq F=0$ on $B_L\times A$. By the maximum principle
on each vertical annulus, we have $U\leq0$ there, so
$U=0$ on $B_L\times\overline A$. Since such balls are available
about every point outside a fixed compact set, the lemma follows.
\end{proof}

\begin{proof}[Second proof: radial barriers]
As in the first proof, choose
\[
  M>\max\{\|u_0\|_{L^\infty(X)},\|u_1\|_{L^\infty(X)}\}.
\]
Take $R$ large enough that the compact core and the supports of both endpoint potentials
lie in $\{r<2R\}$. Choose a smooth cutoff $\eta_R$ such that
\[
  0\leq\eta_R\leq1,\qquad
  \eta_R=1\text{ on }\{r\leq2R\},\qquad
  \eta_R=0\text{ on }\{r\geq4R\},\qquad
  |\ddc\eta_R|_\omega\leq CR^{-2}.
\]
These functions are obtained from a fixed radial cutoff in the ALE
coordinates, extended as one across the compact core; the constant
$C$ is independent of $R$. Put $b_R=-M\eta_R$ and increase $R$
so that $CMR^{-2}\leq\tfrac12$. Then
\[
  \omega+\ddc b_R
  \geq(1-CMR^{-2})\omega\geq\tfrac12\omega.
\]
On $\{r\leq2R\}$, $b_R=-M\leq u_j$; on its complement,
$u_j=0$ and $b_R\leq0$. Thus $b_R\leq u_j$ everywhere, for
$j=0,1$. The time-independent extension of $b_R$ therefore belongs
to $\mathcal B(u_0,u_1)$, and $b_R\leq U$ by \eqref{eq:perron}.

For every fixed $x$, restriction of $\Omega_U\geq0$ to the vertical
annulus implies that $U(x,\cdot)$ is subharmonic. By rotation invariance,
$t\mapsto u_t(x)$ is then convex, so
\[
  b_R(x)\leq u_t(x)\leq(1-t)u_0(x)+tu_1(x).
\]
Both bounds are zero when $r(x)\geq4R$. Taking
$K=\{r\leq4R\}$ completes the proof.
\end{proof}

\subsection{Local Bergman kernels}

We recall the local Bergman approximation used below. Let
\(B_0\Subset\CC^n\) be a Euclidean ball, and let \((\phi_\tau)\)
be a family of plurisubharmonic functions defined on a common
neighborhood of \(\overline{B_0}\). For \(m\geq1\), consider the
weighted Hilbert space
\[
  \mathcal A^2_{m,\tau}(B_0)
  =
  \left\{
    g\in\mathcal O(B_0)
    \;\middle|\;
    \|g\|_{m,\tau}^2
    :=\int_{B_0}|g|^2e^{-m\phi_\tau}\,dV<\infty
  \right\},
\]
where \(dV\) is Euclidean volume. Its diagonal Bergman kernel is
\[
  K_{m,\tau}^{B_0}(z,z)
  =
  \sup\left\{
    |g(z)|^2
    \;\middle|\;
    g\in\mathcal A^2_{m,\tau}(B_0),\quad
    \|g\|_{m,\tau}\leq1
  \right\}.
\]
We use the normalized Bergman density and measure
\begin{equation}\label{eq:bergman-normalization}
  c_n=n!(2\pi)^n,\qquad
  \beta_{m,\tau}^{B_0}
  =c_nm^{-n}K_{m,\tau}^{B_0}e^{-m\phi_\tau},
  \qquad
  \widetilde\beta_{m,\tau}^{B_0}
  =\beta_{m,\tau}^{B_0}\,dV.
\end{equation}
These measures approximate \((\ddc\phi_\tau)^n\)
if the coefficients of the complex Hessian are uniformly bounded in \(L^\infty\).

The estimates are due to Berman--Berndtsson \cite[Theorem~2.1 and Corollary~2.2]{BB}. 
In the following, we write \(H=(\phi_{j\bar k})\) and 
\(\beta_{m, \tau}=\beta_{m,\tau}^{B_0}\) and let \(d\sigma_r\) be the
\textit{probability measure} on \(\partial B(z,r)\). 
For holomorphic \(g\), the submean inequality for \(\log|g|^2\) and Jensen's inequality give
\begin{equation}\label{eq:bergman-submean}
  |g(z)|^2e^{-m\phi(z)}
  \exp\!\left[-m\int_{\partial B(z,r)}
    \bigl(\phi(\zeta)-\phi(z)\bigr)\,d\sigma_r(\zeta)\right]
  \leq
  \int_{\partial B(z,r)}|g|^2e^{-m\phi}\,d\sigma_r.
\end{equation}
By the complex Hessian bound, we have \(0\leq\operatorname{tr}H\leq M\),
with \(M\) independent of \(\tau\), and hence
\[
  0\leq\int_{\partial B(z,r)}
    \bigl(\phi(\zeta)-\phi(z)\bigr)\,d\sigma_r
  \leq\frac{M}{n}r^2.
\]
Now, we fix a compact set \(E\subset B_0\).
Choose \(0<R<\operatorname{dist}(E,\partial B_0)\).
For \(r\leq R/\sqrt m\), the exponential factor on the left of
\eqref{eq:bergman-submean} is bounded below by \(e^{-MR^2/n}\),
independently of \(m\). Integrating \eqref{eq:bergman-submean}
over these spheres, we therefore have
\[
  |g(z)|^2e^{-m\phi(z)}
  \leq
  \frac{e^{MR^2/n}}{\operatorname{Vol}B(z,R/\sqrt m)}
  \int_{B(z,R/\sqrt m)}|g|^2e^{-m\phi}\,dV.
\]
The ball has volume \(m^{-n}\operatorname{Vol}B(0,R)\), and its
integral is at most \(\|g\|_{m,\tau}^2\). Taking the supremum over
\(\|g\|_{m,\tau}\leq1\) and using the normalization of \(\beta_m\)
yields the uniform upper bound
\[
  \beta_{m,\tau}(z)
  \leq\frac{c_n e^{MR^2/n}}{\operatorname{Vol}B(0,R)}
  =:C_E,\qquad z\in E.
\]
This is the desired uniform upper bound for local Bergman measures, independent of \(m, \tau\).

At every common Lebesgue point of the Hessian coefficients, note that 
\( \frac1{r^2}\int_{\partial B(z,r)}
    \bigl(\phi(\zeta)-\phi(z)\bigr)\,d\sigma_r
  \longrightarrow\frac{\operatorname{tr}H(z)}n.\)
Returning to \eqref{eq:bergman-submean}, integrate up to
\(R/\sqrt m\), first let \(m\to\infty\) with \(R\) fixed, and then
let \(R\to\infty\). Together with the Gaussian integral, we have
\[
\limsup_{m\to\infty}\beta_{m,\tau}(z)
\leq n!2^n\left(\frac{\operatorname{tr}H(z)}n\right)^n.
\]
Diagonalize \(H(z)\), with eigenvalues \(\lambda_j\). The changes
\(z_j=\sqrt{\mu_j}\,w_j\), \(\prod_j\mu_j=1\), preserve Euclidean
volume and the Bergman extremal quotient, while replacing
\(\operatorname{tr}H(z)\) by \(\sum_j\mu_j\lambda_j\). The
arithmetic--geometric mean inequality gives
\begin{equation}\label{eq:bergman-trace-determinant}
  \inf_{\mu_j>0}
  \frac{\sum_j\mu_j\lambda_j}{n}
  =\left( \prod_j\lambda_j \right)^{1/n}=\left(\det H(z)\right)^{1/n},
  \qquad \prod_j\mu_j=1.
\end{equation}
For positive eigenvalues, the infimum is attained at
\(\mu_j=(\prod_k\lambda_k)^{1/n}/\lambda_j\).
Using \eqref{eq:bergman-trace-determinant} in the preceding trace
estimate, we have
\[
  \limsup_{m\to\infty}\beta_{m,\tau}(z)
  \leq n!2^n\det H(z)
  =\frac{(\ddc\phi_\tau)^n}{dV}(z).
\]

Their argument first uses
\eqref{eq:bergman-submean} on balls of radius proportional to
\(m^{-1/2}\) to bound the normalized densities on relatively
compact subsets. Optimizing the local estimate over linear changes
of coordinates using \eqref{eq:bergman-trace-determinant}, we obtain the sharp upper bound
\[
  \limsup_{m\to\infty}\beta_{m,\tau}^{B_0}
  \leq
  \frac{(\ddc\phi_\tau)^n}{dV}
  \quad\text{almost everywhere on }B_0.
\]

For fixed \(\tau\), embed \(B_0\subset\CC^n\subset\PP^n\),
where \(\PP^n\) is an auxiliary compact manifold distinct from
\(X\). Berman--Berndtsson extend the local weight to a semipositive
metric with bounded complex Hessian on \(\mathcal O_{\PP^n}(q)\).
The Hilbert--Samuel mass formula and the preceding upper estimates
imply convergence of the global Bergman measures in total variation.
Since restriction of global sections decreases their norms,
comparison with the local kernels yields
\begin{equation}\label{eq:bergman-integral-lower}
  \liminf_{m\to\infty}\int_E\beta_{m,\tau}^{B_0}\,dV
  \geq\int_E(\ddc\phi_\tau)^n.
\end{equation}
Together with the local uniform and sharp upper bounds,
\eqref{eq:bergman-integral-lower} implies convergence in \(L^1(E,dV)\).

In summary, we have the following lemma:

\begin{lemma}[Local Bergman approximation]\label{lem:bergman}
With the notation above, suppose that the coefficients of
\(\ddc\phi_\tau\) are uniformly bounded in
\(L^\infty\) on the common neighborhood of \(\overline{B_0}\).
For every relatively compact subdomain \(E\Subset B_0\), there exists
\(C_E>0\), independent of \(m\geq1\) and \(\tau\), such that
\[
  0\leq\beta_{m,\tau}^{B_0}\leq C_E
  \quad\text{on }E.
\]
For each fixed \(\tau\),
\[
  \widetilde\beta_{m,\tau}^{B_0}
  \longrightarrow(\ddc\phi_\tau)^n
  \quad\text{in total variation on }E.
\]
Equivalently, the corresponding densities converge in \(L^1(E,dV)\).
\end{lemma}
\begin{proof}
This follows directly from Berman--Berndtsson's result and the
preceding discussion.
\end{proof}

\subsection{Continuity and convexity of the Mabuchi functional}

In this subsection, we extend the convexity and continuity arguments
of Berman--Berndtsson \cite{BB} to weak geodesics on a noncompact
K\"ahler manifold whose potentials are supported in a common compact set.

Fix a weak geodesic whose potentials are supported in one compact set
\(K\Subset X\). We use the normalization
\[
  E_\eta(u)
  =
  \frac1n\sum_{j=0}^{n-1}
  \int_X u\,\eta\wedge\omega_u^{n-1-j}\wedge\omega^j,
  \qquad
  \MM(u)=\Ent_{\mathrm{rel}}(u)-nE_{\Ric(\omega)}(u),
\]
where \(\eta\) is a smooth closed real \((1,1)\)-form. If 
\(f_t=\frac{\omega_{u_t}^n}{\omega^n}\),
then \(f_t=1\) outside \(K\), and by Stokes' theorem, we have
\begin{equation}\label{eq:compact-entropy}
  \Ent_{\mathrm{rel}}(u_t) 
  =
  \int_X f_t\log f_t\,\omega^n - \int_X (\omega_{u_t}^n-\omega^n)
  =
  \int_X\log\frac{\omega_{u_t}^n}{\omega^n}\,\omega_{u_t}^n.
\end{equation}
The identity is valid for the bounded-coefficient products by
Bedford--Taylor approximation \cite{BT}. 
All these integrals are finite, since their integrands are supported in \(K\) 
and the densities are uniformly bounded.

\begin{proposition}\label{prop:compact-convexity}
Let \(u_t\), \(t\in[0,1]\), be a weak geodesic with a compact spatial
support \(K\Subset X\). Then \(t\mapsto\MM(u_t)\) is convex.
\end{proposition}

\begin{proof}
We first establish the second-derivative inequality in the sense of
distributions. Continuity of the actual values of the functional will
be proved in Proposition~\ref{prop:compact-continuity}, independently
of pointwise Mabuchi convexity.

Let \(U(x,e^{t+\sqrt{-1}s})=u_t(x)\) be the solution of the HCMA equation
\eqref{eq:hcma} corresponding to the weak geodesic \(u_t\).
For the energy $E_\eta$, set
\[
  \mathbf E_\eta(U)
  =
  \frac1n\sum_{j=0}^{n-1}
  U\,p^*\eta\wedge\Omega_U^{n-1-j}\wedge\Omega^j.
\]
Then $E_\eta(u_\tau)=\pi_*\mathbf E_\eta(U)$, and
$\mathbf E_\eta(U)$ is supported in $K\times A$. Applying
\eqref{eq:twisted-energy-current}, we obtain
\begin{equation}\label{eq:compact-energy-current}
 n\,\ddc_\tau E_\eta(u_\tau)
 =\pi_*\bigl(p^*\eta\wedge\Omega_U^n\bigr).
\end{equation}

Write $\psi_\tau=\log\omega_{u_\tau}^n$, $\psi_\omega=\log\omega^n$, and $\Psi(x,\tau)=\psi_\tau(x)$.
These are local weights with the same transition rule. 
Their difference is the global function \(\log f_\tau\), and
\(\ddc\psi_\omega=-\Ric(\omega)\).
The expression \(\ddc\Psi\wedge\Omega_U^n\) is not
a priori defined: the density may vanish, and \(\Psi\) may equal
\(-\infty\). We therefore truncate the weight before differentiating.

Choose \(V\Subset X\) containing \(K\), and take the smooth reference
weight to be \(\psi_\omega\). Since \(\overline V\) is compact, there
exists \(k_0>0\) such that $-\Ric(\omega)+k_0\omega\geq0$ on $V$.
Define $\chi=p^*\psi_\omega-k_0U$, we have
\begin{equation}\label{eq:compact-truncation-weight}
  \ddc\chi
  =
  p^*(-\Ric(\omega)+k_0\omega)-k_0\Omega_U
  \geq-k_0\Omega_U
  \quad\text{on }V\times A.
\end{equation}
For \(B>0\), we set the truncation data: 
\begin{equation}\label{eq:compact-truncation-data}
\Psi_B=\max\{\Psi,\chi-B\}, \qquad
s_B=\max\{\log f_\tau,-k_0U-B\},
\end{equation}
so \(s_B\) is the logarithmic density \(\log f_\tau\) truncated from below.
Equivalently, $s_B=\Psi_B-p^*\psi_\omega$.
The weight \(\Psi_B\) is locally bounded. Moreover, \(s_B=0\) outside
\(K\times A\), since \(U=0\) and \(f_\tau=1\) there. Define the truncated Mabuchi functional:
\[
  \MM_B(u_\tau)
  =
  \int_X s_B(\cdot,\tau)\,\omega_{u_\tau}^n
  -nE_{\Ric(\omega)}(u_\tau).
\]
The current
\[
  T_B
  =
  \ddc(\Psi_B\Omega_U^n)
  =
  \ddc\Psi_B\wedge\Omega_U^n
\]
is well-defined locally and agrees on overlaps. The transition functions of the weight $\Psi_B$ are pluriharmonic,
and \(\Omega_U^n\) is closed. This also shows that \(T_B\) is a global current. 
It vanishes outside \(K\times A\), because there
\(\Psi_B=\psi_\omega\) and
\[
  \ddc\psi_\omega\wedge\Omega^n
  =-p^*\Ric(\omega)\wedge\Omega^n=0.
\]
The compact spatial support of \(T_B\) therefore makes its fiber
integral well-defined, even though \(X\) is noncompact. To prove
\(T_B\geq0\), it suffices to work on \(V\times A\), where the
curvature bound \eqref{eq:compact-truncation-weight} holds.
Using \eqref{eq:compact-energy-current}, we obtain
\begin{equation}\label{eq:truncated-mabuchi-current}
  \ddc_\tau\MM_B(u_\tau)=\pi_*T_B.
\end{equation}

To prove $T_B\ge 0$, we follow the local positivity argument of Berman--Berndtsson
\cite[proof of Theorem~3.3]{BB}. Fix coordinate balls
\(W\Subset W'\Subset V\) in a common holomorphic chart, and let
\(\rho\) be a local potential
for \(\omega\) near \(\overline{W'}\). Set \(\Phi=\rho+U\), so that
\(\ddc\Phi=\Omega_U\geq0\). For the Bergman densities
\(\beta_m(x,\tau)=\beta_{m,\tau}^{W'}(x)\) in
\eqref{eq:bergman-normalization}, associated with
\(\phi_\tau=\Phi(\cdot,\tau)\), Berndtsson's variation theorem
\cite{Ber} gives
\[
  \ddc\log\beta_m\geq-m\Omega_U.
\]
Let \(\Psi_{B,m}=\max\{\log\beta_m,\chi-B\}\).
Together with \eqref{eq:compact-truncation-weight}, the above lower bound implies that
\(\Psi_{B,m}+m\Phi\) is locally bounded and plurisubharmonic for
\(m\geq k_0\). 
Wedging \(\ddc\Psi_{B,m}\geq-m\Omega_U\) with
\(\Omega_U^n\) and using the HCMA equation, we have
\begin{equation}\label{eq:bergman-positive-current}
  T_{B,m}
  :=
  \ddc\Psi_{B,m}\wedge\Omega_U^n
  \geq-m\Omega_U^{n+1}=0.
\end{equation}

We now pass to the limit on the total space. Fix a subannulus
\(A'\Subset A\). By Lemma~\ref{lem:bergman}, we have
\(\beta_m(\cdot,\tau)\to e^{\psi_\tau}\) in \(L^1(W,dV)\) for
each fixed \(\tau\), together with an upper bound for \(\beta_m\)
on \(W\times A'\) independent of \(m\) and \(\tau\). Applying dominated
convergence in \(\tau\) to the slice norms
\(\|\beta_m(\cdot,\tau)-e^{\psi_\tau}\|_{L^1(W,dV)}\),
we therefore have \(\beta_m\to e^\Psi\) in \(L^1(W\times A')\). For fixed \(B>0\), the positive lower bound
\(c_B:=\inf_{\overline W\times\overline{A'}}e^{\chi-B}>0\) implies
\[
  |\Psi_{B,m}-\Psi_B|
  \leq c_B^{-1}|\beta_m-e^\Psi|,
\]
and hence \(\Psi_{B,m}\to\Psi_B\) in \(L^1(W\times A')\).
Indeed, each weight is the logarithm of its density truncated from
below at \(e^{\chi-B}\). Taking the maximum with this threshold is
\(1\)-Lipschitz in the density, while \(\log\) is
\(c_B^{-1}\)-Lipschitz on \([c_B,\infty)\).
Since \(\Omega_U^n\) has bounded coefficients, we may differentiate the
convergent products \(\Psi_{B,m}\Omega_U^n\) to obtain
\(T_{B,m}\to T_B\) in the sense of currents.
Positivity is preserved under weak convergence of currents, so
\(T_B\geq0\) on \(W\times A'\). Since the coordinate balls
\(W\Subset W'\Subset V\) and the subannulus \(A'\Subset A\) are
arbitrary, these products cover \(V\times A\), and positivity of a
current can be checked locally.
Consequently \(T_B\geq0\) on \(V\times A\), and it is zero outside
\(K\times A\). Its pushforward is a positive current, so
\eqref{eq:truncated-mabuchi-current} proves the distributional
subharmonicity of \(\MM_B\).

Finally, to see the weak convergence of the Mabuchi energy, 
put \(a_\tau=e^{-k_0u_\tau-B}\). The difference of the entropy integrands is
\[
  0\leq
  f_\tau\log\max(f_\tau,a_\tau)-f_\tau\log f_\tau
  \leq\frac{a_\tau}{e}.
\]
It vanishes outside \(K\). The potentials are uniformly bounded there,
and hence
\begin{equation}\label{eq:compact-truncation-error}
  0\leq\MM_B(u_\tau)-\MM(u_\tau)\leq Ce^{-B}
  \quad\text{for all }\tau\in\overline A.
\end{equation}
Thus \(\MM_B\to\MM\) uniformly, and
\[
  \ddc_\tau\MM(u_\tau)\geq0
  \quad\text{in the sense of distributions}.
\]
By rotation invariance, we have
\[
  \frac{d^2}{dt^2}\MM(u_t)\geq0
  \quad\text{on }(0,1)
\]
in the sense of distributions. Proposition~\ref{prop:compact-continuity}
identifies these values with a continuous function on \([0,1]\).
The distributional inequality then implies the asserted pointwise
convexity, including the endpoints.
\end{proof}

\begin{proposition}\label{prop:compact-continuity}
Let \(u_t\), \(t\in[0,1]\), be a weak geodesic with one compact spatial
support \(K\Subset X\). Then \(t\mapsto\MM(u_t)\) is continuous on
\([0,1]\).
\end{proposition}

\begin{proof}
We first show lower semicontinuity. If \(t_\ell\to t\), continuity of
the total potential implies \(u_{t_\ell}\to u_t\) locally uniformly.
Bedford--Taylor continuity \cite{BT}, together with the common compact
support, implies
\[
  E_{\Ric(\omega)}(u_{t_\ell})\longrightarrow E_{\Ric(\omega)}(u_t).
\]
Indeed, the mixed Monge--Amp\`ere products converge weakly and have
uniform local mass bounds.
For the relative entropy, note that the Legendre transform gives
\[
  \Ent_{\mathrm{rel}}(u_t)
  =
  \sup_{g\in C_{\mathrm c}(X;\mathbb R)}
  \left\{
    \int_X g\,\omega_{u_t}^n
    -
    \int_X(e^g-1)\,\omega^n
  \right\}.
\]
The compact support of \(f_t-1\) makes the identity a finite-volume
duality calculation. Every expression in braces is continuous under
weak convergence of the Monge--Amp\`ere measures; hence the relative entropy is lower semi-continuous. 
Consequently,
\begin{equation}\label{eq:compact-mabuchi-lsc}
  \MM(u_t)\leq\liminf_{\ell\to\infty}\MM(u_{t_\ell}).
\end{equation}

To prove upper semicontinuity, we use the localized approximation
of \cite{BB}. Keep \(V\), \(\chi\), and \(\Psi_B\) from the preceding
proof. Choose finitely many coordinate balls
\(B_j\Subset\widehat B_j\Subset V\) covering \(K\), together with
real cutoff functions \(\theta_j\in C^\infty_{\mathrm c}(B_j)\)
satisfying \(\sum_j\theta_j^2=1\) on a neighborhood of \(K\).
Since \(f_t=1\) and \(s_B=0\) defined in \eqref{eq:compact-truncation-data}, outside \(K\). The relative and truncated entropy integrands vanish outside \(K\),
so their localized integrals over any coordinate ball disjoint
from \(K\) are zero. so it's no need to consider the contribution of entropy outside \(K\).
The first and second derivatives of $\theta_j$ are bounded on these fixed charts.
Let \(\beta_{j,m}\) be the local Bergman density on the defining ball
\(\widehat B_j\). For \(B>0\) and \(0<\varepsilon\leq1\), set
\[
  s_{j,m}
  =
  \max\{\log\beta_{j,m},\chi-B\}-\psi_\omega,
  \qquad
  \kappa_\varepsilon(s)=s+\varepsilon(e^s-1),
\]
where we introduce a convex function $\kappa_\varepsilon$ to control the errors produced by differentiating $\theta_j$. 
Define the approximating functionals:
\[
  F_{j,m,B,\varepsilon}(t)
  =
  \int_{B_j}\theta_j^2
  \kappa_\varepsilon(s_{j,m}(t))\,\omega_{u_t}^n.
\]
All weights in this formula are written in the same local frame.
By Lemma~\ref{lem:bergman} and the truncation, we have $-L_B\leq s_{j,m}\leq L$
on a neighborhood of the cutoff support, uniformly in \(m\) and \(t\).
The constants may depend on the fixed charts, the geodesic, and \(B\).

For fixed \(m\), the functionals $F_{j,m,B, \varepsilon}$ are continuous on \([0,1]\).
To see this directly, let \(\phi_t=\rho_j+u_t\) on
\(\widehat B_j\). Comparing the weighted norms, we have
\[
  \|\log\beta_{j,m,t}-\log\beta_{j,m,s}\|_{L^\infty(B_j)}
  \leq
  2m\|\phi_t-\phi_s\|_{L^\infty(\widehat B_j)}.
\]
The right side tends to zero as \(s\to t\). Thus the compactly supported
integrands vary uniformly, while \(\omega_{u_s}^n\) converges weakly
to \(\omega_{u_t}^n\).

We next establish a lower bound for the second derivatives of $F_{j,m,B, \varepsilon}$ independent of \(m\). 
To simplify notation, we write \(s=s_{j,m}\), \(\theta=\theta_j\), and \(\kappa=\kappa_\varepsilon\),
and let \(\Phi\) be a local potential of \(\Omega_U\). 
For \(m\geq k_0\), the function \(s+\psi_\omega+m\Phi\) is a maximum
of locally bounded plurisubharmonic functions, hence is plurisubharmonic
\cite[p.~2089]{BlockiKolodziej07}. 
Thus \(\ddc s\geq-m\Omega_U-C\Omega_{\mathrm E}\)
in the sense of currents, where \(\Omega_{\mathrm E}\) is the
Euclidean K\"ahler form on the product chart. 
A direct computation gives
\[
  \begin{split}
  \ddc\bigl(\theta^2\kappa(s)\bigr)
  &=
  \kappa(s)\ddc\theta^2
  +\theta^2\kappa'(s)\ddc s\\
  &\quad+\sqrt{-1}\alpha\wedge\overline\alpha
  -\frac{4(\kappa'(s))^2}{\kappa''(s)}
       \sqrt{-1}\partial\theta\wedge\bar\partial\theta.
  \end{split}
  \quad \text{ with } \alpha
  =
  \theta\sqrt{\kappa''(s)}\,\partial s
  +
  \frac{2\kappa'(s)}{\sqrt{\kappa''(s)}}\,\partial\theta.
\]
On the bounded range of \(s\), the functions
\(\kappa\), \(\kappa'\), and \((\kappa')^2/\kappa''\) are bounded for
fixed \(B,\varepsilon\), since \(\kappa''(s)=\varepsilon e^s>0\).
Hence
\begin{equation}\label{eq:localized-entropy-curvature}
  \ddc\bigl(\theta^2\kappa(s)\bigr)
  \geq
  -C_1m\Omega_U-C_2(B,\varepsilon^{-1})\Omega_{\mathrm E},
\end{equation}
Taking the Bedford--Taylor product with \(\Omega_U^n\), using the
locally bounded potentials, and applying the HCMA equation yield
\[
  \ddc
  \bigl(\theta^2\kappa(s)\Omega_U^n\bigr)
  \geq
  -C_2(B,\varepsilon)\Omega_{\mathrm E}\wedge\Omega_U^n.
\]
The bounded coefficients of \(\Omega_U\) on the total space allow
fiber integration with a fixed spatial cutoff equal to one near
\(\supp\theta\), giving
\[
  F_{j,m,B,\varepsilon}''(t)\geq-C_{j,B,\varepsilon}
  \quad\text{on }(0,1)
\]
in the sense of distributions. The constant is independent of \(m\).
By the continuity already established,
\[
  \widetilde F_{j,m,B,\varepsilon}(t)
  =F_{j,m,B,\varepsilon}(t)+\frac{C_{j,B,\varepsilon}}2t^2
\]
is continuous and convex on \([0,1]\). Hence, $\widetilde{F}_{j,m,B,\varepsilon}(t)$ satisfies the pointwise convexity property.

Using these functionals, we approximate the entropy step by step.
First keep \(B\) and \(\varepsilon\) fixed. 
Rewrite $s_B$ in \eqref{eq:compact-truncation-data}, $s_B(t)=\max\{\psi_t,\chi-B\}-\psi_\omega$
and define the intermediate functional 
\[
  F_{j,B,\varepsilon}(t)
  =
  \int_{B_j}\theta_j^2
  \kappa_\varepsilon(s_B(t))\,\omega_{u_t}^n.
\]
For every fixed \(t\), we apply Lemma~\ref{lem:bergman} in each spatial coordinate chart at time $t\in [0,1]$ to obtain
\[
  \begin{split}
  |F_{j,m,B,\varepsilon}(t)-F_{j,B,\varepsilon}(t)|
  &\leq
  C_{j,B,\varepsilon}
  \|\beta_{j,m,t}-e^{\psi_t}\|_{L^1(B_j,dV)}
  \longrightarrow0.
  \end{split}
\]
Thus we have pointwise convergence on \([0,1]\).
Set
\[
  \widetilde F_{j,B,\varepsilon}(t)
  =
  F_{j,B,\varepsilon}(t)+\frac{C_{j,B,\varepsilon}}2t^2.
\]
For \(a,b\in[0,1]\) and \(\lambda\in[0,1]\), by pointwise
convexity of $\widetilde F_{j,m,B,\varepsilon}$, we have
\[
  \widetilde F_{j,m,B,\varepsilon}((1-\lambda)a+\lambda b)
  \leq
  (1-\lambda)\widetilde F_{j,m,B,\varepsilon}(a)
  +\lambda\widetilde F_{j,m,B,\varepsilon}(b).
\]
Letting \(m\to\infty\) at these fixed times, we obtain the same inequality
for \(\widetilde F_{j,B,\varepsilon}\). This finite convex function
is continuous on \((0,1)\). Taking \(a=0\), \(b=1\), and letting
\(\lambda\) tend to \(0\) or \(1\), we obtain upper semicontinuity at
the endpoints. Hence, \(F_{j,B,\varepsilon}\) is upper semicontinuous on
\([0,1]\).

Next let \(\varepsilon\downarrow0\). On the bounded truncated range,
\(\kappa_\varepsilon(s)\to s\) uniformly. The local masses
\(\int_{B_j}\theta_j^2\omega_{u_t}^n\) are uniformly bounded, so
\[
  F_{j,B,\varepsilon}(t)\longrightarrow
  F_{j,B}(t):=\int_{B_j}\theta_j^2s_B(t)\,\omega_{u_t}^n
\]
uniformly for \(t\in[0,1]\). Hence \(F_{j,B}\) is upper semicontinuous.
Since \(s_B=0\) outside \(K\) and \(\sum_j\theta_j^2=1\) near \(K\),
\[
  \Ent_B(u_t)
  :=
  \sum_jF_{j,B}(t)
  =
  \int_X s_B(t)\,\omega_{u_t}^n
  =
  \int_X f_t\log\max\{f_t,e^{-k_0u_t-B}\}\,\omega^n
\]
is upper semicontinuous.
Finally, by the estimate \eqref{eq:compact-truncation-error}, we have
\[
  0\leq\Ent_B(u_t)-\Ent_{\mathrm{rel}}(u_t)\leq Ce^{-B}
  \quad\text{uniformly on }[0,1].
\]
It follows that \(\Ent_{\mathrm{rel}}(u_t)\) is upper semicontinuous.
$E_\vartheta(u_t)$ is continuous, so \(\MM(u_t)\) is upper semicontinuous.
Together with \eqref{eq:compact-mabuchi-lsc}, this proves continuity
on the closed interval.
\end{proof}

% ----------------------------------------------------------------------
% 4. Weighted estimates for weak geodesics
% ----------------------------------------------------------------------
\section{Weighted estimates for weak geodesics}\label{sec:weighted-geodesics}

In this section, we study the spatial decay of the weak geodesic
$u_t$ joining $u_0,u_1\in\HH_{-\gamma}(X,\omega)$. We use the
bounded Perron envelope description \eqref{eq:perron} and the joint
$C^{1,1}$ estimate of Proposition~\ref{prop:joint-regularity}.
The weighted $C^0$ bound follows from a common negative barrier.
We then prove that the natural first-derivative decay is also
preserved:
\[
 |u_t|+r|\nabla_Xu_t|_\omega\leq Cr^{-\gamma},
 \qquad 0\leq t\leq1.
\]
The corresponding assertion for second spatial derivatives is false,
even on $\CC^n$ with Schwartz endpoint potentials, as the
Euclidean example in \cite{XiaoYao} shows.

We also estimate the time derivatives in
Section~\ref{subsec:weighted-time-derivatives}. Using a lower barrier
which agrees with both endpoint potentials, we obtain
$|\dot U-(u_1-u_0)|\leq Cr^{-2\gamma-2}$, and hence
$|\dot U|\leq Cr^{-\gamma}$ uniformly on $[0,1]$.
Convexity in time and positivity of $\Omega+\ddc U$ then give
bounds of order $r^{-2\gamma-2}$ for the time integrals of
$\ddot U$, $|\ddot U|^2$, and $|\nabla_X\dot U|_\omega^2$.
The integral estimates will control the spatial cutoff errors
in Lemma~\ref{lem:quantitative-entropy-tails}.

\subsection{The weighted \texorpdfstring{$C^0$}{C0} estimate}\label{subsec:weighted-geodesics}

Let $A=\{\tau\in\CC:1<|\tau|<e\}$, and write
$U(x,\tau)=u_{\log|\tau|}(x)$. Then, the weak geodesic $u_t$ is given by $U(x, \tau)$ with $t = \log|\tau|$.
Here, $U$ is the $S^1$-invariant solution to the HCMA Dirichlet problem
\begin{equation}\label{eq:weighted-hcma}
 \begin{cases}
 (\Omega+\ddc U)^{n+1}=0&\text{on }X\times A,\\
 \Omega+\ddc U\geq0&\text{on }X\times A,\\
 U(x,\tau)=u_i(x)&\text{when }|\tau|=e^i.
 \end{cases}
\end{equation}
Here $\Omega=p^*\omega$, where $p:X\times A\longrightarrow X$
is the projection. We first prove that the weighted $C^0$ decay of the endpoint potentials is preserved uniformly along
this solution.

\begin{proposition}\label{prop:weighted-geodesic-estimates}
Let $u_t$ be the weak geodesic joining
$u_0,u_1\in\HH_{-\gamma}(X,\omega)$. Then
\begin{equation}\label{eq:weighted-geodesic-bound}
 |u_t(x)|\leq C_0r(x)^{-\gamma},\qquad 0\leq t\leq1,
\end{equation}
where $C_0$ is a uniform constant depending only on $n$, $\gamma$, $\omega$, $r$, and $M$, with
\[
 M=\max_{j=0,1}\sup_Xr^\gamma|u_j|.
\]
\end{proposition}

\begin{proof}
We prove the estimate by constructing upper and lower barriers.
Set
\[
 M=\max_{j=0,1}\sup_X r^\gamma|u_j|.
\]
The product $X\times A$
has the trivial holomorphic foliation by the vertical annuli
$\{x\}\times A$. Since $\Omega=p^*\omega$ restricts to zero on
each such annulus, $\Omega+\ddc U\geq0$ implies that
$U(x,\cdot)$ is subharmonic. Its boundary values are constant
on each boundary circle, and their harmonic extension is
\[
 h_x(\tau)=u_0(x)+\log|\tau|\bigl(u_1(x)-u_0(x)\bigr).
\]
By the maximum principle, we have $U(x,\tau)\leq h_x(\tau)$, and hence
\[
 u_t(x)\leq(1-t)u_0(x)+tu_1(x)\leq Mr(x)^{-\gamma}.
\]

For the lower bound, we construct a smooth negative
$\omega$-plurisubharmonic subsolution $w$ with
$|w|\leq C_0r^{-\gamma}$ and $w\leq-Mr^{-\gamma}$. Fix, depending
only on $\gamma$, a smooth positive function $b$ on $[0,\infty)$
such that $b=1$ on $[0,1]$, $b(s)=s^{-\gamma}$ for $s\geq2$, and
\begin{equation}\label{eq:weighted-barrier-profile}
 c_\gamma\min\{1,s^{-\gamma}\}\leq b(s)
 \leq B_\gamma\min\{1,s^{-\gamma}\},\qquad s>0,
\end{equation}
where $0<c_\gamma\leq1$ and $B_\gamma$ depend only on $\gamma$.
Consider
\[
 w=-A_0b(r/L),
\]
where the parameters $A_0>0$ and $L\geq1$ will be chosen below.
Since $r,L\geq1$, the lower bound on $b$ implies
$b(r/L)\geq c_\gamma r^{-\gamma}$. Thus choosing
$A_0=\max\{1,M/c_\gamma\}$ ensures that
\[
 w\leq-c_\gamma A_0r^{-\gamma}\leq-Mr^{-\gamma}\leq u_j,
 \qquad j=0,1.
\]
Fix $L_*\geq1$ so that the compact core is contained in
$\{r<L_*\}$, and constants $A_r,B_r$ such that
$|\nabla r|\leq A_r$ and $|\nabla^2r|\leq B_rr^{-1}$ on the end.
These constants depend only on $(X,\omega)$ and the
radius function $r$. For $L\geq L_*$, the function $w$ is constant
on $\{r\leq L\}$, and by the chain rule, we have
\[
 |\ddc w|_\omega\leq C_bA_0L^{-2},
\]
where $C_b=C_b(n,\gamma,A_r,B_r)$ is independent of $A_0$ and $L$.
We therefore take
\[
 L=\max\{L_*,(2C_bA_0)^{1/2}\}.
\]
By the upper bound on $b$, we have $|w|\leq B_\gamma A_0L^\gamma r^{-\gamma}$.
We may therefore choose
\[
 C_0=\max\{M,B_\gamma A_0L^\gamma\},
\]
whose geometric dependence is through $L_*,A_r,B_r$ alone. With
these parameters fixed,
\begin{equation}\label{eq:negative-floor}
 -C_0r^{-\gamma}\leq w\leq-Mr^{-\gamma},\qquad
 \omega+\ddc w\geq\tfrac12\omega.
\end{equation}

The time-independent extension of $w$ is bounded,
$\Omega$-plurisubharmonic, and lies below both endpoint potentials.
It therefore belongs to $\mathcal B(u_0,u_1)$, so the Perron
description \eqref{eq:perron} implies $U\geq w$. Together with
the upper barrier, this yields
\[
 -C_0r(x)^{-\gamma}\leq w(x)\leq u_t(x)\leq Mr(x)^{-\gamma}.
\]
For fixed $\gamma$ and $(X,\omega)$, the same choices
of $A_0$, $L$, and $w$ work for every pair of endpoint potentials with a common
bound $M$, proving the asserted uniformity.
\end{proof}

We will also use the freedom to enlarge $M$ in this construction.
For every prescribed $M_*>0$, changing the sign of the resulting
subsolution gives a positive function $w$ such that
\begin{equation}\label{eq:adjustable-negative-floor}
 M_*r^{-\gamma}\leq w\leq C_{M_*}r^{-\gamma},\qquad
 \omega-\ddc w\geq\tfrac12\omega.
\end{equation}
Here $C_{M_*}=C_0(n,\gamma,\omega,r,M_*)$ is given by the same
construction. The amplitude is chosen first and the transition
radius $L$ afterwards. Thus the lower bound on $w$ does not prevent
its complex Hessian from being small.

\subsection{The weighted gradient estimate}\label{subsec:weighted-gradient}

We estimate spatial derivatives by comparing $U$ with its composition
with the flows of compactly supported real vector fields.
Using the barrier function in \eqref{eq:adjustable-negative-floor},
we modify the composed potential to obtain an $\Omega$-plurisubharmonic
subsolution, from which we obtain a weighted bound for the difference.

Recall from Proposition~\ref{prop:joint-regularity} that $U$ has a
joint $C^{1,1}$ bound. Writing $\tau=e^{t+\sqrt{-1}s}$, with
$t\in[0,1]$ and $s\in\mathbb R/2\pi\mathbb Z$, we use the fixed
product metric
\[
 \widehat g=g+dt^2+ds^2,\qquad
 \|U\|_{C^{1,1}(X\times\overline A,\widehat g)}\leq M_{1,1}.
\]
Here $g$ is the background metric on $X$. This norm controls all
first and second joint derivatives. Since $U$ is independent of $s$, we also write
$U(x,t)=u_t(x)$ in these logarithmic coordinates.

Fix a point $p$ sufficiently far out on the ALE end and put
$R=r(p)$. Choose real ALE coordinates $y=(y^1,\ldots,y^{2n})$
on a neighbourhood $B_p$ of $p$ such that
$y(B_p)=B(y(p),cR)$. The constant $c>0$ can be chosen independently
of $p$, sufficiently small that this ball lifts to the Euclidean
cover of the end without identifications.
Choose $\chi\in C^\infty_{\mathrm c}(B(0,1))$ with
$0\leq\chi\leq1$ and $\chi=1$ on $B(0,1/2)$. For
$\alpha=1,\ldots,2n$, define
\begin{equation}\label{eq:cutoff-coordinate-fields}
 Z_\alpha(x)=
 \begin{cases}
 R\chi\!\left(\dfrac{y(x)-y(p)}{cR}\right)
     \left.\dfrac{\partial}{\partial y^\alpha}\right|_x,
       &x\in B_p,\\[4pt]
 0,&x\notin B_p.
 \end{cases}
\end{equation}
Because the cutoff vanishes near $\partial B_p$, each $Z_\alpha$
is a smooth compactly supported real vector field on $X$. At $p$,
\[
 \nabla_{Z_\alpha}U(p,t)
 :=du_t(p)\bigl(Z_\alpha(p)\bigr)
 =R\frac{\partial u_t}{\partial y^\alpha}(p).
\]
The ALE metric is equivalent to the Euclidean metric on the end,
so estimates for these directional derivatives control
$|\nabla_Xu_t(p)|_\omega$.

Fix $\alpha$ and write $Z=Z_\alpha$. Consider the flow $F_h:X\to X$
defined by $\frac{d}{dh}F_h(x)=Z(F_h(x))$, with $F_0(x)=x$.
Compact support of $Z$ ensures that the flow is defined for all
real $h$. Each $F_h$ is a diffeomorphism, equals the identity
outside $B_p$. We call it a spatial flow in the sequel.
Let $\widetilde F_h(x,\tau)=(F_h(x),\tau)$.
Set $W_h=U\circ\widetilde F_h$, and
\[
 W_h(x,t)=U(F_h(x),t),\qquad
 \left.\frac{d}{dh}W_h(p,t)\right|_{h=0}
 =\nabla_ZU(p,t).
\]

For a spatial flow, let
$\mathcal P_h:T_xX\longrightarrow T_{F_h(x)}X$ denote parallel
transport along $\ell\mapsto F_\ell(x)$ from $\ell=0$ to $\ell=h$, using
the Levi-Civita connection of $g$. Both $dF_h$ and $\mathcal P_h$
map $T_xX$ to $T_{F_h(x)}X$, so their difference is well-defined.
We denote by $\nabla dF_h$ the covariant derivative of $dF_h$.
The following lemma applies to any real spatial flow satisfying
the stated bounds.

\begin{lemma}\label{lem:spatial-flow-positivity}
Let $U$ be as above, and let $F_h$ be the flow of a smooth real
vector field $Z$, defined on all of $X$ for $|h|<h_0$, where
$0<h_0\leq1$. Use the
notation $W_h=U\circ\widetilde F_h$ and $\mathcal P_h$ introduced
above. Suppose that,
\begin{equation}\label{eq:spatial-flow-bounds}
 |dF_h-\mathcal P_h|+|\nabla dF_h|\leq C_F|h|,
\end{equation}
where $C_F$ is uniform on $X$ for $|h|<h_0$ with $C_F=C_F(n,g,c,\|\chi\|_{C^2}).$
Then
\begin{equation}\label{eq:spatial-flow-positivity}
 \Omega+\ddc W_h\geq-C|h|\Omega
\end{equation}
in the sense of currents, where $C=C(n,M_{1,1},C_F)$ is independent
of $h$ and $h_0$.
\end{lemma}

\begin{proof}
Fix $h$ with $|h|<h_0$ and put $y=F_h(x)$. We compute at
almost every $(x,t)$ where the second derivatives exist. Choose
holomorphic normal coordinates $z=(z^i)$ at $x$ and
$w=(w^\alpha)$ at $y$, with unitary frames identified by
$\mathcal P_h$. This is possible because parallel transport for
the reference K\"ahler metric $g$ is isometric and complex linear.
For this pointwise calculation, write $F=F_h$ and $W=W_h$, and
set
\[
 P_i^\alpha=\frac{\partial F^\alpha}{\partial z^i},\qquad
 Q_i^\alpha=\frac{\partial\overline{F^\alpha}}{\partial z^i}.
\]
Note that $\mathcal P_h$ has the identity matrix
in the chosen frames. By inequality \eqref{eq:spatial-flow-bounds}, we have
\[
 \left\|
 \begin{pmatrix}P-I&Q\\\overline Q&\overline P-I\end{pmatrix}
 \right\|_{\mathrm{op}}
 =\|dF_h-\mathcal P_h\|_{\mathrm{op}}\leq C_F|h|.
\]
Each block has operator norm at most that of the full matrix;
hence
\[
 \|P-I\|_{\mathrm{op}}\leq C_F|h|,\qquad
 \|Q\|_{\mathrm{op}}\leq C_F|h|,\qquad
 \|P\|_{\mathrm{op}}\leq1+C_F|h|.
\]
We use operator norms for matrices below. At the centers of
these normal coordinates, the connection coefficients vanish,
so the second coordinate derivatives of $F$ are components of
$\nabla dF_h$. In the following, we will use $|(F^\alpha)_{i\bar j}|+
 |(\overline{F^\alpha})_{i\bar j}|\leq C(n)C_F|h|$.

We next compute the spatial Hessian of $W$. All derivatives of
$U$ in the following formulas are evaluated at $(y,t)$, and those
of $W$ at $(x,t)$; repeated indices are summed. By the chain rule,
we first have $W_i=P_i^\alpha U_\alpha+Q_i^\alpha U_{\bar\alpha}$,
and then
\[
 \begin{aligned}
 W_{i\bar j}
 ={}&P_i^\alpha\overline{P_j^\beta}\,U_{\alpha\bar\beta}
     +Q_i^\alpha\overline{Q_j^\beta}\,U_{\bar\alpha\beta}\\
 &+P_i^\alpha\overline{Q_j^\beta}\,U_{\alpha\beta}
     +Q_i^\alpha\overline{P_j^\beta}\,U_{\bar\alpha\bar\beta}\\
 &+U_\alpha(F^\alpha)_{i\bar j}
     +U_{\bar\alpha}(\overline{F^\alpha})_{i\bar j}.
 \end{aligned}
\]
Let $G(y,t)$ be the matrix of $\omega+\ddc_X U(\cdot,t)$ at $y$,
and let $G_h(x,t)$ be that of $\omega+\ddc_X W(\cdot,t)$ at $x$.
Thus
\[
 G(y,t)=\bigl(\delta_{\alpha\beta}+U_{\alpha\bar\beta}(y,t)\bigr),
 \qquad
 G_h(x,t)=\bigl(\delta_{ij}+W_{i\bar j}(x,t)\bigr).
\]
The first term in the chain-rule formula contributes
$P(G-I)P^*$, so we can write
\begin{equation}\label{eq:flow-transformed-blocks}
 G_h(x,t)=P\,G(y,t)\,P^*+E_h,\qquad |E_h|\leq C|h|.
\end{equation}
Here a star denotes conjugate transpose, and $E_h$ consists of
$I-PP^*$ and the last five terms in the chain-rule formula.
To verify the error bound, first note that
\[
 |I-PP^*|\leq2C_F|h|+C_F^2|h|^2.
\]
At the normal-coordinate centers, the joint $C^{1,1}$ bound
controls the first derivatives of $U$ and all its second
spatial derivatives, including $U_{\alpha\beta}$. Consequently,
\[
 |E_h|\leq |I-PP^*|
 +C(n)M_{1,1}\bigl(|P|\,|Q|+|Q|^2+|\nabla dF_h|\bigr)
 \leq C|h|,
\]
where $C=C(n,M_{1,1},C_F)$, since $|h|<h_0\leq1$.

It remains to control the full space-time form. Use a local
logarithmic coordinate $\tau=t+\sqrt{-1}s$ on the annulus and suppress
the evaluation points in the matrix notation. At $(y,t)$,
positivity of $\Omega_U$ reads
\[
 \begin{pmatrix}G&b\\b^*&a\end{pmatrix}\geq0,\qquad
 b_\alpha =\tfrac12(U_{\alpha t}),\qquad a=\tfrac14U_{tt}.
\]
The corresponding matrix of $\Omega+\ddc W_h$ at $(x,t)$ has
blocks $G_h,b_h,a_h$. Since $F_h$ is independent of $\tau$
and $U_{\bar\tau}=\frac12\partial_tU$ is real, we have $b_h=Pb+Q\overline b$, $a_h=a$ by the chain rule.
If $a>0$, positivity of the original block matrix implies
\begin{equation}\label{eq:flow-schur-bound}
 \frac{bb^*}{a}\leq G,\qquad
 \frac{|b|^2}{a}\leq\tr G\leq C(n)(1+M_{1,1}).
\end{equation}
Therefore
\[
 \left|\frac{b_hb_h^*}{a}
       -P\frac{bb^*}{a}P^*\right|
 \leq\bigl(2|P|\,|Q|+|Q|^2\bigr)\frac{|b|^2}{a}
 \leq C|h|.
\]
Combining this with \eqref{eq:flow-transformed-blocks} yields
\[
 G_h-\frac{b_hb_h^*}{a}
 \geq P\left(G-\frac{bb^*}{a}\right)P^*-C|h|I
 \geq-C|h|I.
\]
A direct computation gives
\[
\begin{aligned}
\begin{pmatrix}
G_h&b_h\\
b_h^*&a
\end{pmatrix}
&=
\begin{pmatrix}
G_h-\dfrac{b_hb_h^*}{a}&0\\
0&0
\end{pmatrix}
+\frac1a
\begin{pmatrix}b_h\\a\end{pmatrix}
\begin{pmatrix}b_h^*&a\end{pmatrix}\\
&\ge
-C|h|
\begin{pmatrix}
I&0\\
0&0
\end{pmatrix}.
\end{aligned}
\]
If $a=0$, positivity of the original matrix forces
$b=0$, hence $b_h=0$, and
\eqref{eq:flow-transformed-blocks} implies $G_h\geq-C|h|I$.
In both cases, we have $\Omega+\ddc W_h\geq-C|h|\Omega$ almost everywhere,
with $C=C(n,M_{1,1},C_F)$.
\end{proof}

\begin{proposition}\label{prop:weighted-gradient}
Let $u_t$ be the weak geodesic joining
$u_0,u_1\in\HH_{-\gamma}(X,\omega)$. Then
\begin{equation}\label{eq:weighted-gradient-bound}
 |\nabla_Xu_t(x)|_\omega\leq C_1r(x)^{-\gamma-1},
 \qquad 0\leq t\leq1.
\end{equation}
For fixed $(X,\omega)$ and $\gamma$, the constant $C_1$ depends
only on a bound for
\[
 \max_{j=0,1}\sup_X
 r^\gamma\bigl(|u_j|+r|\nabla_Xu_j|_\omega\bigr)
\]
and on the joint $C^{1,1}$ bound for $U$.
\end{proposition}

\begin{proof}
We first show that $r(x)^\gamma|U(F_h(x),t)-U(x,t)|$ converges uniformly to zero in $(x,t)\in X\times[0,1]$ as $h\to0$, with
a bound proportional to $|h|$. Fix $x$ sufficiently far out on
the ALE end and consider the flow $F_h$ generated by a compactly supported vector field satisfying
\eqref{eq:spatial-flow-bounds} and $r(F_h(x))\sim r(x)$
uniformly on their supports for $|h|<h_0$. Precisely, we write it as
$c r(x)\leq r(F_h(x))\leq c^{-1}r(x)$, where
$0<c\leq1$ is independent of $x$ and $h$.
Note that for endpoint potentials, we have
\begin{equation}\label{eq:flow-endpoint-bounds}
 |u_j(F_h(x))-u_j(x)|\leq D|h|r(x)^{-\gamma},
 \quad j=0,1,
\end{equation}
with $D$ independent of $x$ and $h$. Let $Z$, as defined in \eqref{eq:cutoff-coordinate-fields}, be a vector field generating $F_h$.
\eqref{eq:flow-endpoint-bounds} follows by integrating $Zu_j$ along $F_h$: by the weighted endpoint gradient bounds and $|Z|=O(R)$, we have
$|Zu_j|=O(R^{-\gamma})$ on its support, where $r\sim R$ by the preceding discussion.
Since we assume the flow $F_h$ satisfies \eqref{eq:spatial-flow-bounds}, Lemma~\ref{lem:spatial-flow-positivity} bounds the loss of
positivity of $W_h=U\circ\widetilde F_h$ such that
\[
 \Omega+\ddc W_h\geq-C_L|h|\Omega,
\]
where $C_L=C_L(n,M_{1,1},C_F)$. Thus $W_h$ need not be
$\Omega$-plurisubharmonic, and its boundary values need not lie
below $u_j$. We use the barrier in
\eqref{eq:adjustable-negative-floor} to correct both properties.
The resulting function belongs to the Perron class, so we will obtain
the desired estimate by comparison with $U$.
Recall $M=\max_{j=0,1}\sup_Xr^\gamma|u_j|$. The radius comparison
implies
\[
 u_j(F_h(x))\geq-Mc^{-\gamma}r(x)^{-\gamma}.
\]
Choose $M_*=1+Mc^{-\gamma}$ in
\eqref{eq:adjustable-negative-floor}, and let $w$ be the resulting
positive function. Then
\[
 w(x)+u_j(F_h(x))\geq r(x)^{-\gamma},\qquad j=0,1.
\]
Fix $K=\max\{2C_L,D,1\}$. For
$|h|<\min\{h_0,(2K)^{-1}\}$, put $\eta=K|h|<1/2$ and define
\[
 V_h=(1-\eta)W_h-\eta w,
\]
where $w$ is regarded as a function on $X\times A$ independent
of the annulus variable. The barrier supplies a positive spatial term that absorbs the
positivity defect of $W_h$.
Indeed, by Lemma~\ref{lem:spatial-flow-positivity} and
\eqref{eq:adjustable-negative-floor}, we have
\[
 \begin{split}
 \Omega+\ddc V_h
 &=(1-\eta)(\Omega+\ddc W_h)
      +\eta(\Omega-\ddc w)\\
 &\geq\bigl(-C_L|h|(1-\eta)+\eta/2\bigr)\Omega\geq0,
 \end{split}
\]
since $K\geq2C_L$. Thus $V_h$ is $\Omega$-plurisubharmonic.
The same correction also ensures the required boundary inequalities.
At $t=j$, we have
\[
 \begin{split}
 V_h(x,j)-u_j(x)
 &=u_j(F_h(x))-u_j(x)
      -\eta\bigl(u_j(F_h(x))+w(x)\bigr)\\
 &\leq(D-K)|h|r(x)^{-\gamma}\leq0.
 \end{split}
\]
Moreover, $V_h$ is bounded and rotation invariant. Hence
$V_h\in\mathcal B(u_0,u_1)$, and the Perron description
\eqref{eq:perron} implies $V_h\leq U$.
Rearranging this comparison yields
\[
 U(F_h(x),t)-U(x,t)
 \leq\eta\bigl(U(F_h(x),t)+w(x)\bigr).
\]
Proposition~\ref{prop:weighted-geodesic-estimates} and the
radius comparison imply
\[
 |U(F_h(x),t)|\leq C_0c^{-\gamma}r(x)^{-\gamma},
 \qquad w(x)\leq C_{M_*}r(x)^{-\gamma}.
\]
Thus the expression in parentheses is bounded above by a constant
times $r(x)^{-\gamma}$, while the factor $\eta=K|h|$ supplies
the factor $|h|$. In particular,
\[
 U(F_h(x),t)-U(x,t)\leq C_*|h|r(x)^{-\gamma},
 \qquad C_*=K\bigl(C_0c^{-\gamma}+C_{M_*}\bigr).
\]
To obtain the opposite inequality, apply this estimate with
parameter $-h$ at the point $F_h(x)$. Since
$F_{-h}(F_h(x))=x$, we have
\[
 U(x,t)-U(F_h(x),t)
 \leq C_*|h|r(F_h(x))^{-\gamma}
 \leq C_*c^{-\gamma}|h|r(x)^{-\gamma}.
\]
The constant $C=C_*c^{-\gamma}$ depends only on
$n,\gamma,\omega,r,M,M_{1,1},C_F,D,c$ and is uniform over flows
with common constants in the flow bounds, the radius comparison,
and the endpoint bounds. Combining the two inequalities gives
\begin{equation}\label{eq:weighted-flow-comparison}
 |U(F_h(x),t)-U(x,t)|\leq C|h|r(x)^{-\gamma},
 \qquad 0\leq t\leq1.
\end{equation}

We need to estimate the spatial derivative in every direction
at a point $p$ sufficiently far out on the end. Set $R=r(p)$
and, for each $g$-unit vector
$v=\sum_\alpha v^\alpha\partial_{y^\alpha}|_p$, define
$Z=\sum_\alpha v^\alpha Z_\alpha$, where the fields $Z_\alpha$
are given by \eqref{eq:cutoff-coordinate-fields}.
For every such $v$, the vector field $Z$ satisfies
\begin{equation}\label{eq:scaled-spatial-vector-field}
 Z(p)=Rv,\qquad |\nabla^kZ|\leq CR^{1-k},\qquad k=0,1,2,
\end{equation}
where $C=C(n,g,c,\|\chi\|_{C^2})$ is independent of $p$ and $v$.
Let $F_h$ be its flow. We verify the two bounds in
\eqref{eq:spatial-flow-bounds}. Since $Z$ vanishes near
$\partial B_p$ and outside $B_p$, the flow is defined for all
$h\in\mathbb R$ and is the identity outside $B_p$. Its inverse
is $F_{-h}$, so $F_h(B_p)=B_p$. Fix $x\in X$ and write
$J_h=dF_h|_x$. For $v\in T_xX$, differentiating the flow equation
in the spatial direction, we have
\[
 \nabla_h(J_h v)=\nabla_{J_h v}Z,\qquad J_0 v=v,
\]
where $\nabla_h$ denotes covariant differentiation along
$h\mapsto F_h(x)$. By the bound $|\nabla Z|\leq C$ in
\eqref{eq:scaled-spatial-vector-field} and Gr\"onwall's inequality,
we have $|J_h|\leq Ce^{C|h|}$. To compare $J_h$ with parallel
transport, use
\[
 \frac{d}{dh}\bigl(\mathcal P_h^{-1}J_h\bigr)
 =\mathcal P_h^{-1}(\nabla Z)_{F_h(x)}J_h.
\]
Since $\mathcal P_h$ is an isometry, integrating from $0$ to $h$,
we have
\[
 |dF_h-\mathcal P_h|\leq C|h|,\qquad |h|\leq1.
\]

For the second bound, set $B_h=(\nabla dF_h)_x$. Taking a spatial
covariant derivative of the preceding flow equation and commuting
it with $\nabla_h$, we obtain, for $a,b\in T_xX$,
\[
 \begin{split}
 \nabla_h B_h(a,b)
 &=\nabla_{B_h(a,b)}Z
   +(\nabla^2Z)(J_ha,J_hb)\\
 &\quad+\mathrm{Rm}(Z,J_ha)J_hb,
 \qquad B_0(a,b)=0.
 \end{split}
\]
Here the tensors on the right are evaluated at $F_h(x)$, and
we use the convention
$\mathrm{Rm}(A,B)=[\nabla_A,\nabla_B]-\nabla_{[A,B]}$.
The curvature term arises from commuting the covariant
derivatives. On $B_p$, by the bounds
$|\mathrm{Rm}|\leq CR^{-2}$ and
\eqref{eq:scaled-spatial-vector-field}, we have
\[
 |\nabla^2Z|+|\mathrm{Rm}|\,|Z|\leq CR^{-1}.
\]
Together with the bound for $J_h$, this implies
\[
 |\nabla_hB_h|\leq C|B_h|+CR^{-1},\qquad |h|\leq1.
\]
Since $B_0=0$, a second application of Gr\"onwall's inequality yields
\[
 |\nabla dF_h|\leq C|h|R^{-1}.
\]
Outside $B_p$, the map $F_h$ is the identity, so both
$|dF_h-\mathcal P_h|$ and $|\nabla dF_h|$ vanish. Thus these
estimates hold on all of $X$, uniformly for $|h|\leq1$.
As $R\geq1$, they imply \eqref{eq:spatial-flow-bounds} with
$C_F=C_F(n,g,c,\|\chi\|_{C^2})$, independently of $p$, $v$,
and $h$.

Differentiate \eqref{eq:weighted-flow-comparison} at $h=0$.
By the $C^{1,1}$ regularity of $U$, we have
\[
 R\,|du_t(p)(v)|=|Z u_t(p)|\leq CR^{-\gamma}.
\]
Taking the supremum over unit vectors proves
\eqref{eq:weighted-gradient-bound} on the end. The uniform
$C^1$ bound on the compact core completes the proof.
\end{proof}

\subsection{Failure of weighted Hessian estimates}\label{subsec:weighted-hessian-failure}

The construction in \cite[Theorem~B]{XiaoYao} gives
smooth Schwartz potentials $u_0,u_1$ on
$(\CC^n,\omega_0)$, where
$\omega_0=\ddc|z|^2$, such that
\[
 \kappa\omega_0\leq\omega_0+\ddc u_j\leq\Lambda\omega_0,
 \qquad j=0,1,
\]
for fixed $0<\kappa\leq\Lambda$. Nevertheless, the weak
geodesic has spatial rank $n-1$ on open sets arbitrarily far
out, including at the fixed time $t=1/2$. Here Schwartz decay
means that every spatial derivative of the endpoint potentials decays
faster than every inverse power of $|z|$.

The local form of the example explains how this can happen.
Write $z_1=x_1+\sqrt{-1}y_1$ and $z'=(z_2,\ldots,z_n)$. The planar
construction uses a domain whose end is the narrow region
\[
 N=\{x_1>R_0,\ |y_1|<a(x_1)\} \subset \CC^1,\qquad a(x_1)=e^{-x_1}.
\]
The width $a(x_1)$ may also be chosen with sufficiently fast
polynomial decay; the exponential choice yields the Schwartz
example considered here. The boundary contact sets and their
harmonic measures determine the interior time intervals on which
the geodesic agrees with a common obstacle. In $\CC^n$, $n\geq2$,
a cutoff in $z'$ makes the endpoint potentials compactly supported in
these transverse directions, while preserving the construction
where $|z'|<L$ for a sufficiently large fixed $L$.
There are open spatial subsets of $N\times\{|z'|<L\}$
arbitrarily far out and corresponding open time intervals
containing $1/2$ on which the geodesic is smooth, stationary,
and given by
\begin{equation}\label{eq:xiao-yao-hessian-example}
 u_t(z)=-2y_1^2.
\end{equation}
This formula is local to these narrow regions and their
corresponding time intervals. On each such product,
\begin{equation}\label{eq:xiao-yao-nondecaying-hessian}
 \ddc u_t=-\sqrt{-1}\,dz_1\wedge d\bar z_1,\qquad
 \omega_0+\ddc u_t
   =\sqrt{-1}\sum_{k=2}^n dz_k\wedge d\bar z_k.
\end{equation}
Thus $|\ddc u_{1/2}|_{\omega_0}$ stays bounded below by a fixed
positive constant along a sequence tending to infinity:
even decay of the spatial complex Hessian to zero fails.

There is no conflict with the weighted $C^0$ and $C^1$
estimates. On these narrow regions, by the explicit formula, we have
$|u_t|\leq2a(x_1)^2$, and the Euclidean coordinate gradient has
norm $4|y_1|\leq4a(x_1)$. Since $r\sim x_1$ there, the
exponentially shrinking width makes both quantities decay
faster than every inverse power of $r$, despite the nonzero
Hessian. Globally, by Propositions~\ref{prop:weighted-geodesic-estimates}
and~\ref{prop:weighted-gradient}, we have the weighted $C^0$ and
$C^1$ estimates for every $\gamma>0$ for these Schwartz
endpoint potentials.

\subsection{Weighted estimates for the time derivatives}
\label{subsec:weighted-time-derivatives}

We construct a lower barrier which agrees with $u_0$ at $t=0$
and with $u_1$ at $t=1$. Comparing its derivatives at the two
boundary times with those of $U$, we obtain a bound for $\dot U$
on the entire interval. We write $\dot U=\partial_tU$ and
$\ddot U=\partial_t^2U$.

\begin{proposition}\label{prop:weighted-time-derivative}
Let $U$ be the weak geodesic joining
$u_0,u_1\in\HH_{-\gamma}(X,\omega)$. Choose $\kappa>0$ such that
$\omega_{u_j}\geq\kappa\omega$, $j=0,1$, and set
\[
 M_1=\max_{j=0,1}\sup_X r^\gamma
       \bigl(|u_j|+r|\partial u_j|_\omega\bigr).
\]
There is a constant $C$, depending only on $n$, $\gamma$,
$(X,\omega)$, the radius function $r$, $\kappa$, and $M_1$, such that
\begin{equation}\label{eq:weighted-time-first-bound}
 |\dot U(x,t)|\leq Cr(x)^{-\gamma},\qquad
 x\in X,\quad t\in[0,1],
\end{equation}
with one-sided time derivatives at $0$ and $1$.
The constant is uniform for families with common bounds for
$M_1$ and a common positive lower bound $\kappa$.
\end{proposition}

\begin{proof}
Set $v=u_1-u_0$. We seek a lower barrier of the form
\[
 V(x,t)=(1-t)u_0(x)+tu_1(x)-t(1-t)q(x).
\]
We will construct a smooth positive function $q$ with
$q\leq C_qr^{-2\gamma-2}$ satisfying
\begin{equation}\label{eq:time-barrier-conditions}
 |\partial v|_\omega^2\leq\frac\kappa8q,\qquad
 \frac{|\partial q|_\omega^2}{q}\leq\frac\kappa8,\qquad
 |\ddc q|_{\omega,\mathrm{op}}\leq\frac\kappa4,
\end{equation}
where we recall that $\omega_{u_j}\geq\kappa\omega$, $j=0,1$.

Assume first that $q$ satisfies \eqref{eq:time-barrier-conditions}.
The spatial form of $V$ is positive, since
\[
 \omega_{V_t}
 =(1-t)\omega_{u_0}+t\omega_{u_1}-t(1-t)\ddc q
 \geq\left(\kappa-\frac\kappa4t(1-t)\right)\omega
 \geq\frac\kappa2\omega.
\]
Also,
\[
 \dot V=v-(1-2t)q,\qquad \ddot V=2q,
\]
so by the first two conditions, we have
\[
 \begin{aligned}
 |\partial_X\dot V|_{\omega_{V_t}}^2
 &\leq\frac2\kappa
       |\partial v-(1-2t)\partial q|_\omega^2\\
 &\leq\frac4\kappa
       \bigl(|\partial v|_\omega^2+|\partial q|_\omega^2\bigr)
 \leq q\leq\ddot V.
 \end{aligned}
\]
Since $\omega_{V_t}>0$ and
$\ddot V\geq|\partial_X\dot V|_{\omega_{V_t}}^2$, $V$ is a
subgeodesic, so $\Omega+\ddc V\geq0$.
Since $q$ is bounded, $V$ is bounded
and $S^1$-invariant, with $V(x,0)=u_0(x)$ and $V(x,1)=u_1(x)$.
It belongs to $\mathcal B(u_0,u_1)$ in \eqref{eq:perron} and
hence is a lower barrier for $U$.

We now construct $q$. Choose the smooth positive function $b$
as in \eqref{eq:weighted-barrier-profile},
with exponent $2\gamma+2$ in place of $\gamma$. Thus $b=1$ on
$[0,1]$, $b(s)=s^{-2\gamma-2}$ for $s\geq2$, and
\begin{equation}\label{eq:time-barrier-profile}
 c_\gamma\min\{1,s^{-2\gamma-2}\}\leq b(s)
 \leq B_\gamma\min\{1,s^{-2\gamma-2}\},\qquad s>0,
\end{equation}
where $0<c_\gamma\leq1$ and $B_\gamma$ depend only on $\gamma$.
By smoothness and positivity on $[1,2]$ and the explicit
power for $s\geq2$, we may choose a constant $C_\gamma\geq1$, depending
only on $\gamma$, such that
\[
 \frac{|b'(s)|^2}{b(s)}\leq C_\gamma,\qquad
 |b''(s)|\leq C_\gamma,\qquad
 \frac{|b'(s)|}{s}\leq C_\gamma,\qquad s\geq1.
\]
We then pick
\[
 q(x)=A_0b(r(x)/L),
\]
where we first choose $A_0$ and then $L\geq1$.
Because $r,L\geq1$, \eqref{eq:time-barrier-profile} implies
$b(r/L)\geq c_\gamma r^{-2\gamma-2}$: for $r\leq L$ the left
side is one, while for $r>L$ it is at least
$c_\gamma(r/L)^{-2\gamma-2}$. We therefore take
\[
 A_0=\max\left\{1,\frac{32M_1^2}{\kappa c_\gamma}\right\}.
\]
Then, for every $L\geq1$,
\[
 |\partial v|_\omega^2\leq4M_1^2r^{-2\gamma-2}
 \leq\frac\kappa8A_0c_\gamma r^{-2\gamma-2}
 \leq\frac\kappa8q.
\]
The first condition in \eqref{eq:time-barrier-conditions} is
therefore preserved when we enlarge $L$.

To obtain the remaining two conditions, choose $L_*\geq1$ so
that the compact core lies in $\{r<L_*\}$ and
\[
 |\partial r|_\omega\leq C_0,\qquad
 |\ddc r|_{\omega,\mathrm{op}}\leq C_0/r
 \quad\hbox{on }\{r\geq L_*\},
\]
where $C_0\geq1$ and $L_*$ depend only on $(X,\omega)$ and $r$.
For $L\geq L_*$, the function $q$ is constant on $\{r\leq L\}$.
On $\{r>L\}$, write $s=r/L$. Direct differentiation gives
\[
 \partial q=\frac{A_0}{L}b'(s)\partial r,\qquad
 \ddc q=\frac{A_0}{L^2}b''(s)\sqrt{-1}\partial r\wedge\bar\partial r
             +\frac{A_0}{L}b'(s)\ddc r.
\]
Consequently,
\[
 \begin{aligned}
 \frac{|\partial q|_\omega^2}{q}
 &\leq\frac{A_0C_0^2}{L^2}\frac{|b'(s)|^2}{b(s)}
 \leq\frac{A_0C_0^2C_\gamma}{L^2},\\
 |\ddc q|_{\omega,\mathrm{op}}
 &\leq\frac{A_0}{L^2}
       \left(C_0^2|b''(s)|+C_0\frac{|b'(s)|}{s}\right)
 \leq\frac{2A_0C_0^2C_\gamma}{L^2}.
 \end{aligned}
\]
The bounds also hold on $\{r\leq L\}$, where both derivatives
vanish. We now fix
\[
 L=\max\left\{1,L_*,
       \left(\frac{8A_0C_0^2C_\gamma}{\kappa}\right)^{1/2}\right\}.
\]
With this choice, we have $|\partial q|_\omega^2/q\leq\kappa/8$ and
$|\ddc q|_{\omega,\mathrm{op}}\leq\kappa/4$, completing the
verification of \eqref{eq:time-barrier-conditions}. Finally, by
\eqref{eq:time-barrier-profile}, we have
\begin{equation}\label{eq:time-barrier-q-bounds}
 0<q\leq C_qr^{-2\gamma-2},\qquad
 C_q=B_\gamma A_0L^{2\gamma+2}.
\end{equation}
All parameters have the dependencies stated in the proposition.

The constructed $q$ therefore makes $V$ a lower barrier.
Comparing with the Perron envelope and using convexity of
$U(x,\cdot)$, we have
\begin{equation}\label{eq:time-barrier-comparison}
 (1-t)u_0+tu_1-t(1-t)q
       \leq U(x,t)\leq(1-t)u_0+tu_1.
\end{equation}
To estimate the time derivative, fix $x$ and take $0<h<1$.
The comparison \eqref{eq:time-barrier-comparison} gives
\[
 \begin{aligned}
 v-(1-h)q&\leq\frac{U(x,h)-u_0(x)}{h}\leq v,\\
 v&\leq\frac{u_1(x)-U(x,1-h)}{h}\leq v+(1-h)q.
 \end{aligned}
\]
Letting $h\downarrow0$ bounds the two one-sided derivatives.
The derivative of a convex function is nondecreasing, so
\[
 v-q\leq\dot U(x,0+)\leq\dot U(x,t)
          \leq\dot U(x,1-)\leq v+q.
\]
Thus $|\dot U-v|\leq q\leq C_qr^{-2\gamma-2}$.
Since $|v|\leq2M_1r^{-\gamma}$ and $r\geq1$,
\eqref{eq:weighted-time-first-bound} follows with $C=2M_1+C_q$.
\end{proof}

\begin{corollary}\label{cor:integrated-time-derivatives}
Under the assumptions of Proposition~\ref{prop:weighted-time-derivative},
let $M_{1,1}$ be a common joint $C^{1,1}$ bound as in
\eqref{eq:joint-c11-bound}, and suppose that
$0\leq\omega_{u_t}\leq\Lambda\omega$.
With $C_q$ as in \eqref{eq:time-barrier-q-bounds}, for almost every
$x\in X$,
\begin{equation}\label{eq:integrated-time-derivatives}
 \begin{aligned}
 0\leq\int_0^1\ddot U(x,t)\,dt
       &\leq2C_qr(x)^{-2\gamma-2},\\
 \int_0^1|\ddot U(x,t)|^2\,dt
       &\leq2M_{1,1}C_qr(x)^{-2\gamma-2},\\
 \int_0^1|\nabla_X\dot U(x,t)|_g^2\,dt
       &\leq4\Lambda C_qr(x)^{-2\gamma-2}.
 \end{aligned}
\end{equation}
Here the gradient norm is taken with respect to the associated
real metric $g$. The constants depend only on the data defining
$C_q$, $M_{1,1}$, and $\Lambda$, and are uniform for families
with common bounds for those quantities.
\end{corollary}

\begin{proof}
By convexity and joint $C^{1,1}$ regularity, we have
$0\leq\ddot U\leq M_{1,1}$ almost everywhere. The endpoint
bounds proved above imply
\[
 \int_0^1\ddot U(x,t)\,dt
   =\dot U(x,1-)-\dot U(x,0+)
   \leq2q(x)\leq2C_qr(x)^{-2\gamma-2}.
\]
Since $|\ddot U|^2\leq M_{1,1}\ddot U$, the first estimate
implies the second. Finally, by positivity of $\Omega+\ddc U$ and
$\omega_{u_t}\leq\Lambda\omega$, we have
\[
 |\nabla_X\dot U|_g^2
 =2|\partial_X\dot U|_\omega^2\leq2\Lambda\ddot U
 \quad\hbox{almost everywhere on }X\times(0,1).
\]
Integration in time and the first estimate prove the last bound
for almost every $x$.
\end{proof}

% ----------------------------------------------------------------------
% 5. Finiteness of the Mabuchi energy on ALE Kahler manifolds
% ----------------------------------------------------------------------
\section{Finiteness of the Mabuchi energy on ALE K\"ahler manifolds}\label{sec:ale-finiteness}

In this section, we consider the weak geodesic $u_t$ joining
$u_0,u_1\in\HH_{-\gamma}(X,\omega)$. We prove that the relative
Mabuchi functional $\MM(u_t)$ is finite for every $t\in[0,1]$,
as asserted in Theorem~\ref{thm:main-a}.

The main difficulty is that the spatial complex Hessian of a
weak geodesic need not inherit the weighted decay of the endpoint
Hessians. Such decay would directly imply finiteness of the
relative entropy. In the absence of such decay, finiteness on the
infinite-volume manifold requires a separate argument.
Our approach is to approximate the given weak geodesic by weak
geodesics with compact spatial support, for which finiteness and
convexity are already established. We first construct these
approximations and prove convergence of the endpoint Mabuchi
energies, then establish finiteness along the limiting geodesic.

\subsection{Cutoffs and weighted bounds}
We construct compactly supported approximations of the weak geodesic
and prove convergence of the Mabuchi energies at the endpoints.

Fix a smooth function $\chi:[0,\infty)\to[0,1]$ which is one on
$[0,1]$ and zero on $[2,\infty)$. For $R>0$, set
\[
 \chi_R(x)=\chi\bigl(r(x)/R\bigr),\qquad
 u_{j,R}=\chi_Ru_j,\quad j=0,1.
\]
Choose radii $R_1<R_2<\cdots$ tending to infinity and define,
for each endpoint potential, the sequence
\[
 u_{j,l}:=u_{j,R_l}=\chi_{R_l}u_j,
 \qquad j=0,1,\quad l=1,2,\ldots.
\]
Thus $u_{j,l}=u_j$ on $\{r\leq R_l\}$ and
$u_{j,l}=0$ on $\{r\geq2R_l\}$.
Taking $R_1$ sufficiently large ensures that
$\omega+\ddc u_{j,l}>0$ for both endpoint potentials and every $l$.
We therefore obtain a sequence of compactly supported boundary
pairs $(u_{0,l},u_{1,l})\in\HH_{\mathrm c}(X,\omega)^2$.

For each $l$, let $U_l$, also denoted by $U_{R_l}$,
be the weak geodesic joining $u_{0,l}$ to $u_{1,l}$.
More precisely, $U_l$ is the solution of \eqref{eq:hcma}
with boundary values
\[
 U_l(x,\tau)=u_{j,l}(x)
 \qquad\text{when }|\tau|=e^j,\quad j=0,1.
\]
Its existence and joint $C^{1,1}$ regularity follow from
Proposition~\ref{prop:joint-regularity}.
By Lemma~\ref{lem:support}, each $U_l$ vanishes outside a
compact spatial set, independent of $\tau\in\overline A$.
Choose these compact sets successively, enlarging them so that
$K_l$ contains $\{r\leq2R_l\}$, $\supp U_l \subseteq K_l$, and form a compact
exhaustion $K_1,K_2,\ldots$ of $X$:
\[
 K_l\subset K^\circ_{l+1},\qquad
 \{r\leq2R_l\}\subset K_l,\qquad
 \bigcup_{l=1}^\infty K_l=X.
\]
For every $l$, the support property then reads
\[
 U_l=0\quad\text{on }(X\setminus K_l)\times\overline A.
\]
Thus every time slice of the $l$-th approximating geodesic is
supported in the same compact set $K_l$.

\begin{lemma}[Convergence of the endpoint Mabuchi energies]
\label{lem:endpoint-energy-convergence}
Let $u\in\HH_{-\gamma}(X,\omega)$ and set $u_l=\chi_{R_l} u$.
Write $\vartheta=\Ric(\omega)$ for the fixed Ricci form.
Then $u_l\in\HH_{\mathrm c}(X,\omega)$ for all sufficiently large $R_l$,
and, as $R_l\to\infty$,
\[
 \Ent_{\mathrm{rel}}(u_l)\longrightarrow\Ent_{\mathrm{rel}}(u),
 \qquad E_{\vartheta}(u_l)\longrightarrow E_{\vartheta}(u).
\]
Hence the Mabuchi energies also converge:
\begin{equation}\label{eq:smooth-cutoff-mabuchi} 
  \MM(u_l)\longrightarrow\MM(u).
\end{equation}
\end{lemma}

\begin{proof}
We first prove positivity. Set
$A_l=\{R_l<r<2R_l\}$. By the cutoff bounds, we have
\begin{equation}\label{eq:smooth-cutoff-hessian}
 |\nabla^ju_l|_\omega\leq Cr^{-\gamma-j},\qquad 0\leq j\leq k,
\end{equation}
where $C$ depends only on $(X, \omega)$ and radius
function, $\gamma,k,\chi$, and the weighted $C^k$ bound of $u$,
and is independent of $l$. Thus $|\ddc u_l|_\omega\leq
CR_l^{-\gamma-2}$ on $A_l$, so $\omega_{u_l}>0$ there for
large $l$. On $\{r\leq R_l\}$ and $\{r\geq2R_l\}$ the forms
are respectively $\omega_u$ and $\omega$. Hence
$u_l\in\HH_{\mathrm c}(X,\omega)$.

To prove the convergence, $u_l=u$ on $\{r\leq R_l\}$, while the relative
entropy integrand and the integrands defining $E_\vartheta(u_l)$ vanish on
$\{r\geq2R_l\}$. Since the original endpoint integrals are
absolutely convergent, it remains only to estimate the cutoff
contributions on $A_l$. By \eqref{eq:smooth-entropy-tail},
\begin{equation}\label{eq:smooth-cutoff-entropy}
 \int_{A_l}h(f_{u_l})\,\omega^n
 \leq CR_l^{2n-2\gamma-4}\longrightarrow0.
\end{equation}
For the Ricci term, expanding \eqref{eq:twisted-energy}, we obtain
\begin{equation}\label{eq:ricci-linear-split}
 nE_{\vartheta}(u)=\int_Xu S_\omega\,\omega^n+\mathcal E_{\mathrm{rem}}(u),
\end{equation}
where every term in the finite remainder $\mathcal E_{\mathrm{rem}}(u)$ contains a factor
$\ddc u$. For $u_l$, the corresponding remainder integrands on
$A_l$ are bounded by
$C|u_l|\,|\vartheta|_\omega|\ddc u_l|_\omega\omega^n
\leq Cr^{-2\gamma-\delta-4}\omega^n$.
Since $|u_l|\leq|u|$, the cutoff contribution to
$nE_{\vartheta}(u_l)$ is therefore bounded in absolute value by
\[
 \int_{A_l}|u|\,|S_\omega|\,\omega^n
   +CR_l^{2n-2\gamma-\delta-4}\longrightarrow0,
\]
using $S_\omega\in L^1$ and $2\gamma+\delta+4>2n$.
Together with the vanishing tails of the original integrals,
these estimates prove both convergence statements and hence
$\MM(u_l)\to\MM(u)$.
\end{proof}

Let $U_l$ be the lift of the weak geodesic joining $u_{0,l}$ to $u_{1,l}$,
and let $U$ be the lift of the weak geodesic joining $u_0$ to $u_1$.
Since $|u_{j,l}|\leq|u_j|$, by Proposition~\ref{prop:weighted-geodesic-estimates} we have
\[
 |U_l(x,\tau)|\leq Cr(x)^{-\gamma},
\]
where $C=C(n,\gamma,\omega,r,M)$, with
$M=\max_{j=0,1}\sup_Xr^\gamma|u_j|$, is independent of $l$ and
$\tau\in\overline A$. As discussed before, for each $l$, $U_l$ is supported in $K_l$.

By the uniform stability estimate \eqref{eq:perron-stability}
with respect to the boundary data, we have
\begin{equation}\label{eq:cutoff-geodesic-convergence}
 \|U_l-U\|_{L^\infty(X\times\overline A)}
 \leq\max_{j=0,1}\|u_{j,l}-u_j\|_\infty
 \leq CR_l^{-\gamma}.
\end{equation}
The time derivatives also converge locally uniformly:
\begin{equation}\label{eq:cutoff-time-derivative-convergence}
 \partial_tU_l\longrightarrow\partial_tU
 \quad\text{uniformly on }K\times I,\qquad
 K\Subset X,\ I\Subset(0,1).
\end{equation}
Indeed, the cutoff endpoint potentials have uniformly bounded $C^3$
norms, and their K\"ahler forms have a common positive lower bound.
By Proposition~\ref{prop:joint-regularity}, we therefore have a common joint
$C^{1,1}$ bound $M_{1,1}$ for $U_l$ and $U$. Fix a compact set
$K\Subset X$ and a compact interval $I\Subset(0,1)$. By Taylor's formula
in the time variable, for
$0<h<\operatorname{dist}(I,\{0,1\})$, we have
\[
 \|\partial_tU_l-\partial_tU\|_{L^\infty(K\times I)}
 \leq\frac{2}{h}\|U_l-U\|_{L^\infty(X\times\overline A)}
       +M_{1,1}h,
\]
since $|\partial_t^2(U_l-U)|\leq2M_{1,1}$ almost everywhere.
By \eqref{eq:cutoff-geodesic-convergence}, taking
$h=R_l^{-\gamma/2}$ for all sufficiently large $l$ makes the
right-hand side tend to zero, proving
\eqref{eq:cutoff-time-derivative-convergence}.

The uniform convergence in \eqref{eq:cutoff-geodesic-convergence}
will be used in the approximation argument below.

\subsection{Energy estimates and finiteness}\label{sec:energy-finiteness}
In this subsection, we use approximation by compactly
supported weak geodesics to prove finiteness of the Mabuchi energy
along the limiting weak geodesic. We first use the first-derivative
expression for $E_\vartheta$ to establish its finiteness and prove
$E_\vartheta(u_{t,l})\to E_\vartheta(u_t)$ uniformly in $t$ by
interpolation. The argument applies to every $n\geq2$.
The entropy term is treated in Lemma~\ref{lem:entropy-finiteness},
after which we prove Theorem~\ref{thm:main-a}.

\begin{lemma}[Convergence of $E_\vartheta$]\label{lem:twisted-energy-convergence}
Let $u_t$ be the weak geodesic joining
$u_0,u_1\in\HH_{-\gamma}(X,\omega)$, and let $u_{t,l}$ be the
compactly supported approximating weak geodesics constructed above.
Then
\[
 \lim_{l\to\infty}\sup_{t\in[0,1]}
 \left|E_\vartheta(u_{t,l})-E_\vartheta(u_t)\right|=0.
\]
\end{lemma}

\begin{proof}
Let $\vartheta=\Ric(\omega)$. By Propositions~\ref{prop:weighted-geodesic-estimates},
\ref{prop:weighted-gradient}, and \ref{prop:joint-regularity},
for $v=u_t$ or $v=u_{t,l}$, we have
\[
 |v|+r|\nabla_Xv|_\omega\leq C_{\mathrm w}r^{-\gamma},\qquad
 \|v\|_{C^{1,1}(X,g)}\leq M_{1,1},\qquad
 0\leq\omega_v\leq\Lambda\omega.
\]
These constants are independent of $l$ and $t$. The cutoff endpoint potentials
have common weighted $C^1$ bounds, and the common joint $C^{1,1}$ bound.
By \eqref{eq:twisted-energy} and integration by parts, we have
\begin{equation}\label{eq:energy-first-derivative}
 \begin{aligned}
 nE_\vartheta(v)
 &=\int_XvS_\omega\,\omega^n\\
 &\quad-\sum_{j=0}^{n-2}(n-1-j)
 \int_X\sqrt{-1}\partial v\wedge\bar\partial v\wedge\vartheta
                  \wedge\omega_v^j\wedge\omega^{n-2-j}.
 \end{aligned}
\end{equation}
Every integral on the right is absolutely convergent. Indeed,
$S_\omega\in L^1$, and each gradient integrand is bounded by
$Cr^{-2\gamma-\delta-4}\omega^n$, with
$2\gamma+\delta+4>2n$. Integration by parts on $\{r<L\}$,
justified for almost every sufficiently large $L$ by local smooth
approximation, gives boundary integrands on $\{r=L\}$ bounded by
$C|v||\nabla_Xv|_\omega|\vartheta|_\omega
\leq CL^{-2\gamma-\delta-3}$. Since these spheres have area
$O(L^{2n-1})$, the boundary integrals are bounded by
\[
 CL^{2n-1}L^{-2\gamma-\delta-3}
 =CL^{2n-2\gamma-\delta-4}\longrightarrow0
 \qquad\text{as }L\longrightarrow\infty.
\]
Thus, the limit of the second integral over $\{r<L\}$ exists and
equals \eqref{eq:energy-first-derivative}.
In particular, by \eqref{eq:energy-first-derivative}, we have a bound for
$|E_\vartheta(v)|$ independent of $l$ and $t$.

Set $\varepsilon_l=\|U_l-U\|_{L^\infty(X\times\overline A)}$, which
tends to zero by \eqref{eq:cutoff-geodesic-convergence}, and write
$w_{t,l}=u_{t,l}-u_t$. For $x\in X$ and a unit vector $\xi\in T_xX$,
Taylor's formula along the $g$-geodesic with initial velocity $\xi$
gives, for $0<h\leq1$,
\[
 |dw_{t,l}(x)[\xi]|
 \leq\frac{2\varepsilon_l}{h}
       +\frac h2\|\nabla_X^2w_{t,l}\|_{L^\infty(X)}
 \leq\frac{2\varepsilon_l}{h}+M_{1,1}h,
\]
since the $C^{1,1}$ bound on each potential implies
$\|\nabla_X^2w_{t,l}\|_{L^\infty(X)}\leq2M_{1,1}$.
Taking $h=\sqrt{\varepsilon_l}$ and the supremum over $x$, $\xi$,
and $t$, we obtain, for all sufficiently large $l$,
\[
 \sup_{t\in[0,1]}
 \|\nabla_X(u_{t,l}-u_t)\|_{L^\infty(X)}
 \leq(2+M_{1,1})\sqrt{\varepsilon_l}.
\]
The constant is independent of $l$ and $t$.

For fixed $l,t$, set $v=u_t+s w_{t,l}$, $0\leq s\leq1$.
Differentiate the mixed expression over $\{r<L\}$ along this segment
and integrate by parts, using $d\vartheta=0$. The segment satisfies the same metric
and weighted bounds, so the boundary terms have the same vanishing
bound as above. Passing to the limit as $L\to\infty$, we obtain
\[
 \begin{aligned}
 n\bigl(E_\vartheta(u_{t,l})-E_\vartheta(u_t)\bigr)
 &=\int_Xw_{t,l}S_\omega\,\omega^n\\
 &\quad-n\sum_{j=0}^{n-2}\int_0^1\!\int_X
 \sqrt{-1}\partial w_{t,l}\wedge\bar\partial v\wedge\vartheta
       \wedge\omega_v^j\wedge\omega^{n-2-j}\,ds.
 \end{aligned}
\]
These calculations are valid for $C^{1,1}$ potentials by local
smooth approximation and integration by parts on compact sets.
The interpolation estimate,
$|\nabla_Xv|\leq C_{\mathrm w}r^{-\gamma-1}$,
and the common metric upper bound yield
\[
 \begin{aligned}
 n|E_\vartheta(u_{t,l})-E_\vartheta(u_t)|
 &\leq\varepsilon_l\|S_\omega\|_{L^1(X,\omega^n)}
  +C\sqrt{\varepsilon_l}
       \int_Xr^{-\gamma-1}|\vartheta|_\omega\,\omega^n\\
 &\leq C(\varepsilon_l+\sqrt{\varepsilon_l})\longrightarrow0.
 \end{aligned}
\]
All constants in this energy estimate depend only on $n$, $(X,\omega)$,
$r$, $\gamma$, $\delta$, $C_{\mathrm w}$, $M_{1,1}$, and $\Lambda$;
they are independent of $l$ and $t$. Hence
$E_\vartheta(u_{t,l})\to E_\vartheta(u_t)$ uniformly in $t$.
\end{proof}

\begin{lemma}[Finiteness of the relative entropy]\label{lem:entropy-finiteness}
Let $u_t$ be the weak geodesic joining
$u_0,u_1\in\HH_{-\gamma}(X,\omega)$. Then the relative entropy
is uniformly bounded along the geodesic:
\[
 \sup_{t\in[0,1]}\Ent_{\mathrm{rel}}(u_t)<\infty.
\]
\end{lemma}

\begin{proof}
We first bound the relative entropies of the compactly supported
approximating geodesics $u_{t,l}$ uniformly in $l$ and $t$. 
By Lemma~\ref{lem:endpoint-energy-convergence}, we have
\[
 \MM(u_{j,l})\longrightarrow\MM(u_j)\in\mathbb R,
 \qquad j=0,1.
\]
Both endpoint-energy sequences are therefore uniformly bounded.
Propositions~\ref{prop:compact-convexity} and \ref{prop:compact-continuity} imply that there is a uniform constant $C_0$, independent of $l$ and $t$, such that
\[
 \MM(u_{t,l})\leq(1-t)\MM(u_{0,l})+t\MM(u_{1,l})\leq C_0.
\]

The preceding estimates bound $E_\vartheta(u_{t,l})$
independently of $l$ and $t$. Consequently,
\[
 \Ent_{\mathrm{rel}}(u_{t,l})
 =\MM(u_{t,l})+nE_\vartheta(u_{t,l})\leq C.
\]

For each fixed $t$, \eqref{eq:cutoff-geodesic-convergence} and local
Bedford--Taylor continuity imply
$\omega_{u_{t,l}}^n\rightharpoonup\omega_{u_t}^n$ against compactly
supported tests. By lower semicontinuity of the relative entropy
\eqref{eq:entropy-lsc}, we have
\[
 0\leq\Ent_{\mathrm{rel}}(u_t)
 \leq\liminf_{l\to\infty}\Ent_{\mathrm{rel}}(u_{t,l})\leq C,
\]
where $C$ is independent of $t$. This proves the lemma.
\end{proof}

\begin{proof}[Proof of Theorem~\ref{thm:main-a}]
By Lemma~\ref{lem:entropy-finiteness}, the relative
entropy is finite along the weak geodesic. Fix $t\in[0,1]$ and put $u=u_t$.
The weighted $C^0$ and metric bounds allow us to apply
Lemma~\ref{lem:entropy-hessian}. By \eqref{eq:hessian-l2-tail},
$|A_u-I|\in L^2(X,\omega^n)$. We verify absolute convergence of
the original mixed integrals defining $E_\vartheta(u)$.
For each summand, we have
\[
 \begin{aligned}
 u\vartheta\wedge\omega_u^j\wedge\omega^{n-1-j}
 &=\frac{uS_\omega}{n}\,\omega^n\\
 &\quad+u\vartheta\wedge\ddc u\wedge
              \sum_{a=0}^{j-1}\omega_u^a\wedge\omega^{n-2-a},
 \end{aligned}
\]
with an empty sum when $j=0$. The first term is integrable because $u$ is bounded and $S_\omega\in L^1$.
The remaining term is bounded by
$Cr^{-\gamma-\delta-2}|A_u-I|\,\omega^n$, which is integrable
by Cauchy--Schwarz, since $|A_u-I|\in L^2$ and
$\gamma+\delta+2>n$.
The integration by parts justified in the proof of
Lemma~\ref{lem:twisted-energy-convergence} identifies the mixed
integrals with \eqref{eq:energy-first-derivative}.
Hence both terms in
$\MM(u_t)=\Ent_{\mathrm{rel}}(u_t)-nE_{\Ric(\omega)}(u_t)$ are finite
for every $t\in[0,1]$, proving Theorem~\ref{thm:main-a}.
\end{proof}

% ----------------------------------------------------------------------
% 6. Mabuchi convexity on ALE Kahler manifolds
% ----------------------------------------------------------------------
\section{Mabuchi convexity on ALE K\"ahler manifolds}\label{sec:ale}

We now prove that $t\mapsto\MM(u_t)$ is continuous and convex
on $[0,1]$, as asserted in Theorem~\ref{thm:main-b}, for the weak
geodesic joining $u_0,u_1\in\HH_{-\gamma}(X,\omega)$.

Passing compact-support convexity to the limit requires further
control of the relative entropy, since uniform convergence of
the potentials implies only lower semicontinuity. The main
difficulty is to prove convergence of the relative entropy at
every interior slice. We establish convergence on compact subsets
and a uniform estimate for the entropy outside a large compact
set, and then combine these two statements.

\subsection{Framework of the proof}\label{subsec:convexity-framework}

We retain the compactly supported weak geodesics $U_l$ introduced
in Section~\ref{sec:ale-finiteness}, with slices $u_{t,l}$ and
endpoint potentials $u_{j,l}=\chi_{R_l}u_j$, $j=0,1$.
In particular, the cutoff endpoint potentials agree exactly with $u_j$ on
$\{r\leq R_l\}$, and $U_l\to U$ uniformly on
$X\times\overline A$ by \eqref{eq:cutoff-geodesic-convergence}.
Throughout this section, we use the following notations:
\[
 d\mu=\omega^n,\qquad
 f_{t,l}=\frac{\omega_{u_{t,l}}^n}{\omega^n},\qquad
 f_t=\frac{\omega_{u_t}^n}{\omega^n},\qquad
 h(s)=s\log s-s+1.
\]
Thus
$\Ent_{\mathrm{rel}}(u_t)=\int_Xh(f_t)\,d\mu$.
By the estimates already proved, there are constants $C_{\mathrm w}$,
$M_{1,1}$, $\Lambda\geq1$, and $H_*$, independent of $l,t$, such that
\begin{equation}\label{eq:ale-common-bounds}
 \begin{split}
 |U_l|+r|\nabla_XU_l|_\omega &\leq C_{\mathrm w}r^{-\gamma},
 \\ 
 \|U_l\|_{C^{1,1}(X\times\overline A,\widehat g)} &\leq M_{1,1},
 \\
 0\leq\omega_{u_{t,l}} &\leq\Lambda\omega,
 \\
 \Ent_{\mathrm{rel}}(u_{t,l}) &\leq H_*.
 \end{split}
\end{equation}
The same bounds hold for $U$ and its slices. We also fix a common
$\kappa>0$ with $\omega_{u_{j,l}}\geq\kappa\omega$, $j=0,1$,
as established for the cutoff endpoint potentials in
Section~\ref{sec:ale-finiteness}. Together with
\eqref{eq:ale-common-bounds} at $t=0,1$, the lower bound makes
the constants in Corollary~\ref{cor:integrated-time-derivatives}
uniform in $l$. We use
$\vartheta=\Ric(\omega)$, $\gamma=2n-4+\varepsilon'$ with
$\varepsilon'>0$, and $\delta>n-1$ as in
\eqref{eq:ale-setting}--\eqref{eq:gamma}.
Lemma~\ref{lem:twisted-energy-convergence} already establishes uniform
convergence of $E_\vartheta(u_{t,l})$ on $[0,1]$.

The first lemma concerns the density on a fixed compact set.
Uniform convergence of potentials alone does not imply its
$L^1$ convergence. The exact agreement of the cutoff endpoint potentials
on expanding compact sets is also used.

\begin{lemma}[Convergence on compact subsets]
\label{lem:compact-density-convergence}
For every $t\in(0,1)$ and every compact set $K\Subset X$,
\[
 \|f_{t,l}-f_t\|_{L^1(K,d\mu)}\longrightarrow0,\qquad
 \|h(f_{t,l})-h(f_t)\|_{L^1(K,d\mu)}\longrightarrow0.
\]
\end{lemma}

For the second lemma, fix a sufficiently large radius $R$,
independently of the approximation index $l$, and set
$\theta_R$ such that $\theta_R=0$ on $\{r\leq R\}$, $\theta_R=1$ on
$\{r\geq2R\}$. Note by definition, $\theta_R = 1- \chi_R$ , and its first and second derivatives are bounded
by $C R^{-1}$ and $C R^{-2}$. Define the entropy outside the
cutoff region by
\[
 H_{l,R}(t)=\int_X\theta_R h(f_{t,l})\,d\mu.
\]
Then, we have the following lemma:

\begin{lemma}[Uniform estimates for the entropy tails]
\label{lem:uniform-entropy-tails}
There exist constants $C>0$ and $R_0\geq1$ such that, for every
$l$, $R\geq R_0$, and $0<a<b<1$,
\begin{equation}\label{eq:ale-tail-chord}
 H_{l,R}(t)
 \leq\frac{b-t}{b-a}H_{l,R}(a)
       +\frac{t-a}{b-a}H_{l,R}(b)+\varepsilon_R,
 \qquad t\in[a,b],
\end{equation}
where
\begin{equation}\label{eq:ale-tail-error-rate}
 \varepsilon_R=C(1+\log R)R^{-(n-2+\varepsilon')}
       \longrightarrow0\qquad\text{as }R\longrightarrow\infty.
\end{equation}
The constants $C,R_0$ depend only on $n$, $(X,\omega)$, $r$,
the chosen cutoff profiles, $\gamma$, $\delta$, $\kappa$, and the
bounds $C_{\mathrm w}$, $M_{1,1}$, $\Lambda$, and $H_*$ in
\eqref{eq:ale-common-bounds}.
\end{lemma}

The estimate concerns the convexity inequality for the entropy
tails. It does not assert that those tails are already uniformly
small. The latter conclusion at a fixed interior time follows
only after the next proposition.

\begin{proposition}[Convergence of the entropy at interior slices]
\label{prop:interior-entropy-convergence}
For every $t\in(0,1)$,
\begin{equation}\label{eq:ale-entropy-convergence}
 \Ent_{\mathrm{rel}}(u_{t,l})
       \longrightarrow\Ent_{\mathrm{rel}}(u_t).
\end{equation}
Consequently $\MM(u_{t,l})\to\MM(u_t)$ at every such $t$.
\end{proposition}

Here is the logic connecting the two lemmas to the proposition.
The convex functions $\MM(u_{t,l})$ have locally uniformly
convergent subsequences on $(0,1)$. If $G$ is such a limit,
by lower semicontinuity we have $D(t):=G(t)-\MM(u_t)\geq0$.
Convergence on compact subsets identifies $D(t)$ with the
remaining entropy at infinity. Passing
\eqref{eq:ale-tail-chord} to the limit shows that $D$ is convex.
Convergence of the endpoint energies and lower semicontinuity
at the endpoints imply $\lim_{t\to0,1}D(t)=0$, so $D=0$.
This proves the proposition for the entire sequence.
Theorem~\ref{thm:main-b} then follows by passing the
compact-support convexity inequality to the limit and checking
continuity at the endpoints.

\subsection{Convergence on compact subsets}
\label{subsec:compact-density-convergence}

In this subsection, we prove Lemma~\ref{lem:compact-density-convergence}. The argument uses
the exact agreement of the cutoff endpoint potentials with the original endpoint potentials
on a fixed compact set. To prove this, we consider the two cases $\dot u_t\ne0$ and
$\dot u_t=0$. In the first case, we show that $u_{t,l}$ agrees
with $u_t$ in a neighborhood of each point for all sufficiently
large $l$. On the set $\{\dot u_t=0\}$, we prove convergence
of the Monge--Amp\`ere densities using the contact-set
properties of the rooftop envelopes.

\begin{proof}[Proof of Lemma~\ref{lem:compact-density-convergence}]
For $\lambda\in\RR$, set $P_\lambda=P_\omega(u_0,u_1-\lambda)$, $P_{l,\lambda}=P_\omega(u_{0,l},u_{1,l}-\lambda)$,
using the rooftop envelope defined in \eqref{eq:rooftop-definition}. 
We write $U(x,t)=u_t(x)$ in the logarithmic annulus coordinate.
By Proposition~\ref{prop:rooftop-representation}, applied to
$U$ and $U_l$, we have
\begin{equation}\label{eq:ale-rooftop-duality}
 P_\lambda(x)=\inf_{0\leq s\leq1}\{U(x,s)-\lambda s\},
 \qquad
 U(x,t)=\sup_{\lambda\in\RR}\{P_\lambda(x)+\lambda t\},
\end{equation}
and the same identities hold for $P_{l,\lambda}$ and $U_l$.
We use $P_{\cdot,\lambda}$ to denote either $P_{l,\lambda}$,
for a fixed $l$, or $P_\lambda$, with the same choice throughout
each formula. The infimum formulas imply
\begin{equation}\label{eq:ale-rooftop-stability}
 \|P_{l,\lambda}-P_\lambda\|_{L^\infty(X)}
 \leq\|U_l-U\|_{L^\infty(X\times\overline A)},
 \qquad
 \|P_{\cdot,\lambda}-P_{\cdot,\lambda'}\|_{L^\infty(X)}
 \leq|\lambda-\lambda'|.
\end{equation}
By \eqref{eq:rooftop-c11} and the common joint $C^{1,1}$ bound,
for every compact set $K\Subset X$ we have
\begin{equation}\label{eq:ale-rooftop-derivative-bound}
 \|\nabla_X P_{\cdot,\lambda}\|_{L^\infty(K)}
 +\|\nabla_X^2P_{\cdot,\lambda}\|_{L^\infty(K)}
 \leq C_K.
\end{equation}
Here derivatives and norms use the reference metric $g$, with
the Hessian understood almost everywhere. The constant $C_K$
depends only on $n$, $K$, $(X,\omega)$, and $M_{1,1}$ in
\eqref{eq:ale-common-bounds}, independently of $l,\lambda$.

Define their densities by
\[
 g_{\cdot,\lambda}
   =\frac{\omega_{P_{\cdot,\lambda}}^n}{\omega^n}.
\]
Here $g_{\cdot,\lambda}$ denotes $g_{l,\lambda}$ or $g_\lambda$,
according to the choice of $P_{\cdot,\lambda}$.
By \eqref{eq:ale-rooftop-derivative-bound}, these densities satisfy
common bounds on each compact set.
By Lemma~\ref{lem:local-density-weak} and
\eqref{eq:ale-rooftop-stability}, for every fixed relatively
compact measurable set $E\subset X$, we have
\begin{equation}\label{eq:ale-rooftop-density-weak}
 \begin{aligned}
 \int_E g_{l,\lambda}\,\omega^n
 &\longrightarrow\int_E g_\lambda\,\omega^n
       &&(l\to\infty,\ \lambda\text{ fixed}),\\
 \int_E g_\lambda\,\omega^n
 &\longrightarrow\int_E g_0\,\omega^n
       &&(\lambda\to0).
 \end{aligned}
\end{equation}
Applying the same lemma to \eqref{eq:cutoff-geodesic-convergence},
we also have $\int_E f_{t,l}\,\omega^n\to\int_E f_t\,\omega^n$
for each fixed $t$ and $E$.
Finally, \eqref{eq:rooftop-contact-measures} yields
\begin{equation}\label{eq:ale-rooftop-measures}
 \begin{cases}
 g_\lambda=0
 \quad&\hbox{a.e.} \quad \text{ on }
 \{P_\lambda<\min(u_0,u_1-\lambda)\},\\
 g_\lambda=\frac{\omega_{u_j}^n}{\omega^n}
 \quad&\hbox{a.e.} \quad \text{ on }
 \{P_\lambda=u_j-j\lambda\},\quad j=0,1.
 \end{cases}
\end{equation}
The same assertions hold for $g_{l,\lambda}$ with the cutoff
endpoint potentials.

\subsubsection*{The case $\dot u_t(x)\ne0$}
The K\"ahler forms $\omega_{u_j}$ and $\omega_{u_{j,l}}$ have a common
positive lower bound
\[
 \omega_{u_j},\ \omega_{u_{j,l}}\geq\kappa\omega,\qquad j=0,1.
\]
Here, $R_1$ is chosen sufficiently large such that on $\{r\leq R_l\}$
$\omega_{u_{j,l}}=\omega_{u_j}$, while on $\{r>R_l\}$ it is uniformly
close to $\omega$. Fix $a>0$ and a compact set $K\Subset X$.
Choose \(C_0>0\), independent of \(l\), such that
$\|u_1-u_0\|_{L^\infty}$,
$\|u_{1,l}-u_{0,l}\|_{L^\infty}<C_0$
for every \(l\). 
The endpoint decay and the derivative bounds for $\chi_R$
allow us to choose $R$ such that $K\subset\{r<R\}$ and
\[
 |u_{1,l} (x)-u_{0,l}(x)|\leq a/2\quad(r(x)>R),\qquad
 C_0|\ddc\chi_R|_\omega\leq\kappa/2,
\]
uniformly in $l$, with the same bounds for the original endpoint potentials.

For $\lambda\geq a$, the function
\[
 v=u_{1,l}-\lambda-C_0\chi_R
\]
is bounded and $\omega$-psh. It can be checked that $v$ lies below $u_{1,l}-\lambda$
everywhere and below $u_{0,l}$. In $\{r\leq R\}$ this follows from
the choice of $C_0$, and in $\{r>R\}$ from
$u_{1,l}-\lambda\leq u_{0,l}-a/2$. Since $v$ equals the lower
obstacle outside $\{r<2R\}$, $v = u_{1,l}-\lambda \le u_{0,l}$, the envelope satisfies
\[
 P_{l,\lambda}=u_{1,l}-\lambda\quad(r>2R).
\]
For $\lambda\leq-a$, the same argument with
$v=u_{0,l}-C_0\chi_R$ gives $P_{l,\lambda}=u_{0,l}$ there.
The identical constructions apply to $P_\lambda$.

For sufficiently large $l>l_0$, we can assume the cutoff endpoint potentials agree with
the original endpoint potentials on $\{r<4R\}$. By the preceding identities,
we therefore have
\[
 P_{l,\lambda}=P_\lambda
 \qquad\text{on }\{2R<r<3R\},\qquad |\lambda|\geq a:
\]
both envelopes equal $u_1-\lambda$ when $\lambda\geq a$, and $u_0$
when $\lambda\leq-a$.

We can consequently glue $P_{l,\lambda}$ on $\{r<3R\}$ to $P_\lambda$
on $\{r>2R\}$, as illustrated in Figure~\ref{fig:rooftop-gluing}. Since the two functions agree on the overlap, they
define a globally bounded $\omega$-psh function. This function lies
below both original obstacles $u_0$ and $u_1-\lambda$. On $\{r<3R\}$,
this follows from the agreement of the endpoint potentials, and on
$\{r>2R\}$, from the definition of $P_\lambda$. The defining property
of $P_\lambda$ then implies $P_{l,\lambda}\leq P_\lambda$ on $K$.
Conversely, gluing $P_\lambda$ inside to $P_{l,\lambda}$ outside, we obtain
a bounded $\omega$-psh function below the cutoff obstacles, and hence
the reverse inequality. Thus
\begin{equation}\label{eq:ale-rooftop-exact}
 P_{l,\lambda}=P_\lambda\quad\hbox{on }K,\qquad |\lambda|\geq a,
\end{equation}
for all sufficiently large $l$, with the lower bound for $l$
independent of $\lambda$.

\begin{figure}[htbp]
\centering
\begin{tikzpicture}[x=1cm,y=1cm,
  every node/.style={font=\small,inner sep=2pt}]
\path[use as bounding box] (0,-0.15) rectangle (14.4,4.1);
\foreach \x/\inner/\outer in
 {3.0/{P_{l,\lambda}}/{P_\lambda},
  10.2/{P_\lambda}/{P_{l,\lambda}}}{
 \begin{scope}[shift={(\x,2.15)}]
  \draw[densely dashed,line width=0.45pt] (0,0) circle (1.45);
  \draw[densely dashed,line width=0.45pt] (0,0) circle (1.02);
  \node[font=\normalsize] at (0,0) {$\inner$};
  \node[font=\normalsize] at (2.05,-0.55) {$\outer$};
  \draw[line width=0.4pt] (0.72,0.72)--(1.30,1.25)--(1.90,1.25);
  \node[anchor=west] at (1.90,1.25) {$r=2R$};
  \draw[line width=0.4pt] (1.23,0.77)--(1.75,0.85)--(2.22,0.85);
  \node[anchor=west] at (2.22,0.85) {$r=3R$};
 \end{scope}
 \node at (\x,0.17) {Glue $\inner$ inside $\outer$.};
}
\end{tikzpicture}
\caption{The gluing step.}
\label{fig:rooftop-gluing}
\end{figure}

Consider a relatively compact neighborhood in $X\times(0,1)$ on
which $|\partial_tU|\geq2a>0$. By
\eqref{eq:cutoff-time-derivative-convergence}, we have
$|\partial_tU_l|\geq a$ there for all sufficiently large $l$.
Thus both supporting slopes
\[
 \lambda=\partial_tU(x,t),\qquad
 \lambda_l=\partial_tU_l(x,t)
\]
have absolute value at least $a$. By convexity in time, each is a
maximizing slope for the corresponding identity
\eqref{eq:ale-rooftop-duality}. Using \eqref{eq:ale-rooftop-exact}
at these two slopes, we have
\[
 \begin{aligned}
 U_l(x,t)&\geq P_{l,\lambda}(x)+t\lambda
                  =P_\lambda(x)+t\lambda=U(x,t),\\
 U(x,t)&\geq P_{\lambda_l}(x)+t\lambda_l
                  =P_{l,\lambda_l}(x)+t\lambda_l=U_l(x,t).
 \end{aligned}
\]
Consequently,
\begin{equation}\label{eq:ale-geodesic-exact}
 U_l=U\quad\hbox{locally on }\{\partial_tU\ne0\}
 \quad\hbox{for all sufficiently large }l.
\end{equation}

\subsubsection*{The case $\dot{u}_t=0$}
Fixing $t\in(0,1)$, throughout the discussion of the case, we consider the sets: $E_t=\{x:\partial_tU(x,t)=0\}$ and $Z=\{P_0=u_0=u_1\}$.
For $x\in E_t$, by convexity of $s\mapsto u_s(x)$, we have
\[
 u_s(x)\geq u_t(x)+(s-t)\dot u_t(x)=u_t(x),
 \qquad s\in[0,1].
\]
Thus $u_t(x)$ is the minimum along this interval, and by Legendre duality
\eqref{eq:ale-rooftop-duality} we have
\[P_0(x)=\inf_{s\in[0,1]}u_s(x)=u_t(x).\]
Moreover, if $x \in Z$, $P_0\leq u_s\leq(1-s)u_0+su_1$ shows that $u_s=u_0=u_1$
on $Z$ for every $s\in[0,1]$. In particular, $Z\subset E_t$.
Since $u_t=P_0$ on $E_t$, the second part of the proof of
Lemma~\ref{lem:rooftop-contact-measures} gives
$\ddc u_t=\ddc P_0$ and hence $f_t=g_0$ almost everywhere on $E_t$.

In the following proof, we separate $Z$ and $E_t\setminus Z$, as the limiting density
behaves differently on these two sets. 
On $E_t\setminus Z$, we will
prove that $f_t=0$ almost everywhere. Lemma~\ref{lem:local-density-weak}
and nonnegativity of $f_{t,l}$ then imply local $L^1$ convergence. On $Z$,
the limiting density equals $\omega_{u_0}^n/\omega^n>0$ almost
everywhere. We will show that, for every compact set $K\Subset X$,
the $\omega^n$-volume of
$\{x\in Z\cap K:u_{t,l}(x)\ne u_t(x)\}$ tends to zero as $l\to\infty$.
The Hessian agreement proved in Lemma~\ref{lem:rooftop-contact-measures},
together with the common density bound, then implies
$f_{t,l}\to f_t$ in $L^1(Z\cap K,\omega^n)$.

We claim that
\begin{equation}\label{eq:ale-zero-slope-density}
 f_t=0\quad\hbox{almost everywhere on }E_t\setminus Z.
\end{equation}
On a compact set $K$, define, for $m\geq1$,
\[
 F_{0,m}=E_t\cap K\cap\{P_0=u_0,\ u_1-u_0\geq1/m\}.
\]
On $F_{0,m}\subset E_t$, we have $u_t=P_0$ by the preceding argument.
By the definition of $F_{0,m}$, we also have $P_0=u_0$ there.
Thus $u_t=u_0$ on $F_{0,m}$.
For $0<\lambda<1/(2m)$, by \eqref{eq:ale-rooftop-duality} we have
\[
 P_\lambda\leq u_t-t\lambda
             =u_0-t\lambda<u_0<u_1-\lambda
 \quad\hbox{on }F_{0,m}.
\]
Thus $g_\lambda=0$ almost everywhere on $F_{0,m}$ by
\eqref{eq:ale-rooftop-measures}. Applying
\eqref{eq:ale-rooftop-density-weak} to the fixed set $F_{0,m}$, we obtain
\[
 \int_{F_{0,m}}g_0\,\omega^n
 =\lim_{\lambda\downarrow0}\int_{F_{0,m}}g_\lambda\,\omega^n=0.
\]
Hence $g_0=0$ almost everywhere on $F_{0,m}$.

For
\[
 F_{1,m}=E_t\cap K\cap\{P_0=u_1,\ u_0-u_1\geq1/m\},
\]
take $-1/(2m)<\lambda<0$. Then
\[
 P_\lambda\leq u_1-t\lambda<u_1-\lambda,
 \qquad
 P_\lambda\leq u_1+|\lambda|<u_0.
\]
Again $g_\lambda=0$ almost everywhere on $F_{1,m}$, and
\eqref{eq:ale-rooftop-density-weak} gives
\[
 \int_{F_{1,m}}g_0\,\omega^n
 =\lim_{\lambda\uparrow0}\int_{F_{1,m}}g_\lambda\,\omega^n=0.
\]
Thus $g_0=0$ almost everywhere on $F_{1,m}$.
The countable unions of these sets cover the one-sided contacts
in $(E_t\setminus Z)\cap K$. At all remaining points,
$P_0<\min(u_0,u_1)$, where \eqref{eq:ale-rooftop-measures}
already implies $g_0=0$. Since $f_t=g_0$ almost everywhere on $E_t$, this proves
\eqref{eq:ale-zero-slope-density}.

It remains to control the common contact set $Z$. For large $l$,
both cutoff endpoint potentials agree with the original endpoint potentials on a
neighborhood of $K$. Set
\[
 Z_l=Z\cap K\cap\{P_{l,0}=u_0\},\qquad
 f_0=\frac{\omega_{u_0}^n}{\omega^n}.
\]
On $Z\cap K$, we have $u_{0,l}=u_{1,l}=u_0=u_1$.
By the definition of $Z_l$, at every point of $(Z\cap K)\setminus Z_l$,
\[
 P_{l,0}<u_{0,l}=u_{1,l}.
\]
Thus, by \eqref{eq:ale-rooftop-measures}, we have $g_{l,0}=0$ almost everywhere
on $(Z\cap K)\setminus Z_l$. On $Z_l$, we have $P_{l,0}=u_0=u_{0,l}$, so the contact-set
identity in the same formula implies $g_{l,0}=f_0$ almost everywhere.
Since $P_0=u_0$ on $Z\cap K$, we also have $g_0=f_0$ there by that identity.
Consequently,
\[
 g_{l,0}=f_0\mathbf1_{Z_l},
 \qquad g_0=f_0
 \quad\hbox{almost everywhere on }Z\cap K.
\]
Applying \eqref{eq:ale-rooftop-density-weak} to the fixed measurable
set $Z\cap K$, we obtain
\[
 \int_{Z_l}f_0\,\omega^n
 =\int_{Z\cap K}g_{l,0}\,\omega^n
 \longrightarrow\int_{Z\cap K}g_0\,\omega^n
 =\int_{Z\cap K}f_0\,\omega^n.
\]
Since $Z_l\subset Z\cap K$, subtracting from the last integral yields
\[
 \int_{(Z\cap K)\setminus Z_l}f_0\,\omega^n
 =\int_{Z\cap K}f_0\,\omega^n-\int_{Z_l}f_0\,\omega^n
 \longrightarrow0.
\]
Since $f_0\geq\kappa^n$, the volume of $(Z\cap K)\setminus Z_l$ tends
to zero. On $Z_l$,
\[
 u_0=P_{l,0}\leq u_{t,l}
 \leq(1-t)u_{0,l}+tu_{1,l}=u_0.
\]
Thus $u_{t,l}=u_t=u_0$ on $Z_l$. The second part of the proof of
Lemma~\ref{lem:rooftop-contact-measures} gives $f_{t,l}=f_t$
almost everywhere on $Z_l$. The common upper bound for
the slice densities therefore proves
\[
 \int_{Z\cap K}|f_{t,l}-f_t|\,\omega^n\longrightarrow0.
\]

On $(E_t\setminus Z)\cap K$, the limiting density is zero by
\eqref{eq:ale-zero-slope-density}. Applying
Lemma~\ref{lem:local-density-weak} to the fixed measurable set
$(E_t\setminus Z)\cap K$ and using $f_{t,l}\geq0$, we have
\[
 \begin{aligned}
 \int_{(E_t\setminus Z)\cap K}|f_{t,l}-f_t|\,\omega^n
 &=\int_{(E_t\setminus Z)\cap K}f_{t,l}\,\omega^n\\
 &\longrightarrow\int_{(E_t\setminus Z)\cap K}f_t\,\omega^n=0.
 \end{aligned}
\]
Finally, on each compact set
$K\cap\{|\partial_tU(\cdot,t)|\geq1/m\}$, a finite covering and
\eqref{eq:ale-geodesic-exact} imply equality of the densities for
all sufficiently large $l$. The remaining subsets of
$K\setminus E_t$ have volume tending to zero as $m\to\infty$.
The common density bound proves convergence there as well.
Combining the three regions yields
\[
 \|f_{t,l}-f_t\|_{L^1(K,\omega^n)}\longrightarrow0.
\]
Uniform continuity of $h$ on $[0,\Lambda^n]$ then implies
$h(f_{t,l})\to h(f_t)$ in $L^1(K,\omega^n)$.
\end{proof}

\subsection{Uniform estimates for the entropy tails}
\label{subsec:entropy-tail-proof}

In this part, we prove Lemma~\ref{lem:uniform-entropy-tails}. 
The positivity of the Bergman kernel argument from Section~\ref{sec:compact} remains valid for each $U_l$.
We apply the convexity argument to $H_{l,R}$, but derivatives of the
spatial cutoff $\theta_R$ produce additional error terms. The following
lemma establishes a convexity inequality for $H_{l,R}$ with explicit bounds
on the cutoff and truncation errors, uniform in $l$.
Lemma~\ref{lem:quantitative-entropy-tails} asserts the same convexity inequality
as Lemma~\ref{lem:uniform-entropy-tails}, with a more detailed estimate
of the error $\varepsilon_R$ in \eqref{eq:ale-tail-chord}.

\begin{lemma}\label{lem:quantitative-entropy-tails}
Fix $k_0>0$ such that $k_0\omega-\vartheta\geq0$ on $X$.
There exist constants $C>0$ and $R_0\geq1$ such that, for every
$l$, $R\geq R_0$, $0<a<b<1$, and $B>k_0C_{\mathrm w}+2$,
\begin{equation}\label{eq:ale-tail-chord-quantitative}
 \begin{aligned}
 H_{l,R}(t)
 &\leq\frac{b-t}{b-a}H_{l,R}(a)
       +\frac{t-a}{b-a}H_{l,R}(b)\\
 &\quad+Ce^{-B}+2\zeta_R
       +C(1+B)R^{n-\gamma-2},
 \qquad t\in[a,b],
 \end{aligned}
\end{equation}
where $C_{\mathrm w}$ is the weighted bound in
\eqref{eq:ale-common-bounds}, and
\begin{equation}\label{eq:ale-energy-tail-rate}
 \zeta_R=C\left(
 R^{-\gamma}\int_{\{r>R\}}|S_\omega|\,d\mu
       +R^{n-\gamma-\delta-2}\right)\longrightarrow0
       \qquad\text{as }R\longrightarrow\infty.
\end{equation}
The constants $C,R_0$ depend only on the data listed in
Lemma~\ref{lem:uniform-entropy-tails} and $k_0$; in particular,
they are independent of $a,b,l,B$ and the spatial supports.
\end{lemma}

\begin{proof}[Proof of Lemma~\ref{lem:quantitative-entropy-tails}]
We apply the Bergman kernel argument of Section~\ref{sec:compact}
to the full relative Mabuchi energy outside the cutoff region.
The spatial cutoff produces error terms in its second time
derivative. We estimate those terms, obtain a convexity inequality
with a small error, and then recover the entropy by showing that
the corresponding $E_\vartheta$ term tends to zero.

Fix $I=[a,b]\Subset(0,1)$ and a sufficiently large radius $R$,
independently of $l$. Recall that $\theta_R=0$ on $\{r\leq R\}$
and $\theta_R=1$ on $\{r\geq2R\}$. With the convention
$\mathcal E=nE_\vartheta$ following \eqref{eq:mabuchi}, define
\[
 \begin{aligned}
 \E_{l,R}(t)
   &=\sum_{j=0}^{n-1}\int_X\theta_R u_{t,l}\,
       \vartheta\wedge\omega_{u_{t,l}}^j\wedge\omega^{n-1-j},\\
 \M_{l,R}(t)&=H_{l,R}(t)-\E_{l,R}(t).
 \end{aligned}
\]
Thus, $\E_{l,R}$ is the expression for $nE_\vartheta$ with the
factor $\theta_R$, and $\M_{l,R}$ is the full relative Mabuchi
energy with the same factor. As in
\eqref{eq:compact-truncation-data}, put the truncated entropy integrand as follows:
\[
 s_{l,B}=\max\{\log f_{t,l},-k_0U_l-B\},
\]
and define
\[
 \begin{aligned}
 H^B_{l,R}(t)
   &=\int_X\theta_R(f_{t,l}s_{l,B}-f_{t,l}+1)\,d\mu,\\
 \M^B_{l,R}(t)&=H^B_{l,R}(t)-\E_{l,R}(t).
 \end{aligned}
\]
Only the logarithmic term is truncated. In particular,
$\M^B_{l,R}-\M_{l,R}=H^B_{l,R}-H_{l,R}\geq0$.
For each fixed $l$, all the integrands have compact spatial
support, so the displayed integrals are finite.

\subsubsection*{The second derivative of $\M^B_{l,R}$.}
Use the logarithmic coordinate $\tau=t+\sqrt{-1}s$, where $s$ has period
$2\pi$, and write $\Omega=p^*\omega$ and
$\Omega_l=\Omega+\ddc U_l$. We retain the notation
$\mathbf E_\vartheta(U_l)$ from the proof of
Proposition~\ref{prop:compact-convexity}:
\[
 n\mathbf E_\vartheta(U_l)
 =U_l\,p^*\vartheta\wedge
       \sum_{j=0}^{n-1}\Omega_l^j\wedge\Omega^{n-1-j}.
\]
In particular, $\E_{l,R}=n\pi_*(\theta_R\mathbf E_\vartheta(U_l))$.
For the truncated Mabuchi functional, define the $(n,n)$-form
\[
 \mathbf M_B(U_l)
   =s_{l,B}\Omega_l^n-n\mathbf E_\vartheta(U_l).
\]
Both forms have bounded measurable coefficients and compact
spatial support. In a local holomorphic frame, write
$\psi_\omega=\log\omega^n$. The positive current from approximation Bergman kernel argument of
Section~\ref{sec:compact} is
\[
 T_{l,B}=\ddc(\psi_\omega+s_{l,B})\wedge\Omega_l^n\geq0.
\]
Its definition is independent of the frame, and its positivity
follows from \eqref{eq:bergman-positive-current} and the subsequent
limiting argument. The computation in
\eqref{eq:twisted-energy-current} gives
\begin{equation}\label{eq:ale-current-primitive}
 n\ddc\mathbf E_\vartheta(U_l)
       =p^*\vartheta\wedge\Omega_l^n,\qquad
 \ddc\mathbf M_B(U_l)=T_{l,B}.
\end{equation}
Here $p^*\vartheta\wedge\Omega^n=0$ and
$\ddc\psi_\omega=-p^*\vartheta$ have been used.
The identities hold in the sense of currents, without
estimating derivatives of $s_{l,B}$.

The scalar functional is
\begin{equation}\label{eq:ale-localized-functional}
 \M^B_{l,R}
   =\pi_*(\theta_R\mathbf M_B(U_l))
       -\pi_*\bigl(\theta_R(\Omega_l^n-\Omega^n)\bigr).
\end{equation}
The last term is the mass correction in the relative entropy.
Let $\eta\in C_c^\infty((a,b))$, regarded as an $S^1$-invariant
test function on $A$, and set
\[
 L_R=\frac{\sqrt{-1}}{2}\bigl(\partial\theta_R\wedge d\bar\tau
                       +d\tau\wedge\bar\partial\theta_R\bigr).
\]
We calculate the second derivative against $\eta$:
\begin{equation}\label{eq:ale-test-product-rule}
 \begin{aligned}
 \int_A\eta\,\ddc\M^B_{l,R}
 &=\int_{X\times A}\theta_R\ddc\eta\wedge
       \bigl(\mathbf M_B(U_l)-(\Omega_l^n-\Omega^n)\bigr)\\
 &=\int_{X\times A}
       \bigl(\ddc(\theta_R\eta)-\eta\,\ddc\theta_R-\eta'L_R\bigr)
             \wedge\mathbf M_B(U_l)\\
 &\quad-\int_{X\times A}\eta\,\ddc\theta_R
                       \wedge(\Omega_l^n-\Omega^n)\\
 &=\int_{X\times A}\theta_R\eta\,T_{l,B}
       -\int_{X\times A}\eta\,\ddc\theta_R
                        \wedge\bigl(\mathbf M_B(U_l)+\Omega_l^n\bigr)\\
 &\quad-\int_{X\times A}\eta'L_R\wedge\mathbf M_B(U_l).
 \end{aligned}
\end{equation}
Here we used closedness of $\Omega_l^n-\Omega^n$,
$\ddc\theta_R\wedge\Omega^n=0$, and
$\ddc\mathbf M_B(U_l)=T_{l,B}$ from
\eqref{eq:ale-current-primitive}.
Since $\ddc\M^B_{l,R}=\tfrac14(\M^B_{l,R})''
\sqrt{-1}\,d\tau\wedge d\bar\tau$, we obtain
\begin{equation}\label{eq:ale-distribution-errors}
 (\M^B_{l,R})''=\nu_{l,R,B}+e_0+\partial_t e_1,
 \qquad \nu_{l,R,B}\geq0,
\end{equation}
where
\begin{equation}\label{eq:ale-errors-defined}
 \begin{aligned}
 \nu_{l,R,B}\,\sqrt{-1}\,d\tau\wedge d\bar\tau
   &=4\pi_*(\theta_RT_{l,B}),\\
 e_0\,\sqrt{-1}\,d\tau\wedge d\bar\tau
   &=-4\pi_*\bigl(\ddc\theta_R
                 \wedge(\mathbf M_B(U_l)+\Omega_l^n)\bigr),\\
 e_1\,\sqrt{-1}\,d\tau\wedge d\bar\tau
   &=4\pi_*(L_R\wedge\mathbf M_B(U_l)).
 \end{aligned}
\end{equation}
All coefficients are $S^1$-invariant. The main term
$\nu_{l,R,B}$ is a positive distribution. The integrands defining
$e_0,e_1$ contain, respectively, second and first spatial derivatives
of $\theta_R$ and vanish outside $A_R=\{R<r<2R\}$. More explicitly,
\eqref{eq:ale-distribution-errors} means
\[
 \int_I\M^B_{l,R}\eta''\,dt
  =\langle\nu_{l,R,B},\eta\rangle
       +\int_Ie_0\eta\,dt-\int_Ie_1\eta'\,dt.
\]
All integrations by parts are justified by compact spatial
support for each $l$; no bound for the size of that support
is used.

We will prove the following estimates for the error terms:
\begin{equation}\label{eq:ale-uniform-error-bounds}
 \begin{aligned}
 \sup_l\|e_0\|_{L^1(I)}
      &\leq C(1+B)R^{2n-2\gamma-4},\\
 \sup_l\|e_1\|_{L^1(I)}
      &\leq C(1+B)R^{n-\gamma-2}.
 \end{aligned}
\end{equation}
We will also show that
\[
 \sup_{l,t}|\E_{l,R}(t)|\leq\zeta_R\longrightarrow0,
 \qquad
 0\leq\M^B_{l,R}(t)-\M_{l,R}(t)\leq Ce^{-B},
\]
where $\zeta_R$ is defined in \eqref{eq:ale-energy-tail-rate}.
Because $\gamma+2-n>0$, both powers in
\eqref{eq:ale-uniform-error-bounds} are negative.
We will use the error estimates to prove a convexity inequality for
$\M^B_{l,R}$; the bounds for $\E_{l,R}$ and the truncation error
will then imply the required inequality for $H_{l,R}$.
\subsubsection*{Estimates for the cutoff errors.}
In the estimates below, $C$ denotes a positive constant that may
change from line to line and depends only on $n$, $(X,\omega)$,
$r$, the chosen cutoff profiles, $\gamma$, $\delta$, $k_0$,
and the common constants $C_{\mathrm w}$, $M_{1,1}$, $\Lambda$,
$H_*$, and $\kappa$. It is independent of $x,t,l,R,B$,
the interval $I=[a,b]$, and the supports of the approximating
geodesics. We use three ingredients in proving the error estimates. The integral control of the complex
Hessian by the entropy in Lemma~\ref{lem:entropy-hessian},
positivity of the total complex Hessian matrix, and the time
integral estimates of Corollary~\ref{cor:integrated-time-derivatives}.
We recall the bounds before applying them to $e_0,e_1$.

Write $A_{t,l}=\omega^{-1}\omega_{u_{t,l}}$; then
$0\leq A_{t,l}\leq\Lambda I$ and $f_{t,l}=\det A_{t,l}$.
All matrix norms below are Hilbert--Schmidt norms in an
$\omega$-unitary spatial frame. By the entropy bound in
\eqref{eq:ale-common-bounds} and
Lemma~\ref{lem:entropy-hessian}, we have
\begin{equation}\label{eq:ale-global-hessian}
 \sup_{t,l}\int_X|A_{t,l}-I|^2\,d\mu
       \leq C_{n,\Lambda}H_*.
\end{equation}
In the same frame, the coefficient matrix of $\Omega_l$ is
\begin{equation}\label{eq:ale-joint-block}
 \begin{pmatrix}A_{t,l}&b_l\\b_l^*&c_l\end{pmatrix},
 \qquad
 (b_l)_\alpha=\tfrac12\partial_\alpha\partial_tU_l,
 \qquad c_l=\tfrac14\partial_t^2U_l.
\end{equation}
The matrix is nonnegative and its coefficients are uniformly
bounded by the joint $C^{1,1}$ estimate. In particular,
$c_l\geq0$ and $|b_l|^2\leq\Lambda c_l$.
Integrating \eqref{eq:integrated-time-derivatives} over $A_R$
and using $\operatorname{Vol}_\omega(A_R)\leq CR^{2n}$, we obtain
\begin{equation}\label{eq:ale-annular-time-bounds}
 \int_{A_R\times[0,1]}(c_l+c_l^2+|b_l|^2)\,d\mu\,dt
       \leq CR^{2n-2\gamma-2}.
\end{equation}
Consequently,
\begin{equation}\label{eq:ale-mixed-l2-bounds}
 \|c_l\|_{L^2(A_R\times I)}
       +\|b_l\|_{L^2(A_R\times I)}
       \leq CR^{n-\gamma-1}.
\end{equation}
The constants are independent of $I$, since the time integral
bounds hold on $[0,1]$.

We first express the cutoff terms in terms of $b_l,c_l$.
Let $\mathcal C_l$ be the spatial block of the adjugate of
\eqref{eq:ale-joint-block}. The adjugate is nonnegative:
for a positive definite matrix $H$,
$\operatorname{adj}(H)=\det(H)H^{-1}\geq0$, and adding a positive
multiple of the identity and passing to the limit proves the
claim for a nonnegative matrix. The diagonal cofactor
$(\mathcal C_l)_{\alpha\bar\alpha}$ is the determinant of the
principal submatrix obtained by deleting the $\alpha$-th
spatial row and column. Hadamard's inequality therefore gives
\[
 0\leq(\mathcal C_l)_{\alpha\bar\alpha}
 \leq c_l\prod_{\beta\ne\alpha}(A_{t,l})_{\beta\bar\beta}
 \leq\Lambda^{n-1}c_l.
\]
By positivity, we also have
\[
 |(\mathcal C_l)_{\alpha\bar\beta}|^2
 \leq(\mathcal C_l)_{\alpha\bar\alpha}
       (\mathcal C_l)_{\beta\bar\beta}.
\]
The mixed entries of the full adjugate are
$-\operatorname{adj}(A_{t,l})b_l$ and their conjugates. Hence
\begin{equation}\label{eq:ale-cofactor-bounds}
 |\mathcal C_l|\leq C_{n,\Lambda}c_l,
 \qquad
 |\operatorname{adj}(A_{t,l})b_l|\leq C_{n,\Lambda}|b_l|.
\end{equation}
The identities remain valid at singular matrices by continuity.

If $\ddc\theta_R=\sqrt{-1}(D_R)_{\alpha\bar\beta}
 dz^\alpha\wedge d\bar z^\beta$, expanding the wedge product,
we obtain
\[
 \ddc\theta_R\wedge\Omega_l^n
   =\tr(D_R\mathcal C_l)\,
        \omega^n\wedge \sqrt{-1}\,d\tau\wedge d\bar\tau.
\]
The wedge product $L_R\wedge\Omega_l^n$ involves only the mixed
cofactors. Since $|\ddc\theta_R|\leq CR^{-2}$ and
$|d\theta_R|\leq CR^{-1}$, by \eqref{eq:ale-cofactor-bounds} we have
\[
 \begin{aligned}
 |\ddc\theta_R\wedge\Omega_l^n|
     &\leq CR^{-2}c_l\,d\mu\wedge \sqrt{-1}\,d\tau\wedge d\bar\tau,\\
 |L_R\wedge\Omega_l^n|
     &\leq CR^{-1}|b_l|\,d\mu\wedge \sqrt{-1}\,d\tau\wedge d\bar\tau.
 \end{aligned}
\]
Absolute values here refer to the top-degree coefficients.

Substituting the definition of $\mathbf M_B(U_l)$ into
\eqref{eq:ale-errors-defined}, we separate the entropy, mass,
and $E_\vartheta$ contributions:
\begin{equation}\label{eq:ale-error-components}
 \begin{aligned}
 e_0^{\rm ent}\,\sqrt{-1}\,d\tau\wedge d\bar\tau
   &=-4\pi_*(s_{l,B}\ddc\theta_R\wedge\Omega_l^n),\\
 e_0^{\rm mass}\,\sqrt{-1}\,d\tau\wedge d\bar\tau
   &=-4\pi_*(\ddc\theta_R\wedge\Omega_l^n),\\
 e_0^E\,\sqrt{-1}\,d\tau\wedge d\bar\tau
   &=4n\pi_*(\ddc\theta_R\wedge\mathbf E_\vartheta(U_l)),\\
 e_1^{\rm ent}\,\sqrt{-1}\,d\tau\wedge d\bar\tau
   &=4\pi_*(s_{l,B}L_R\wedge\Omega_l^n),\\
 e_1^E\,\sqrt{-1}\,d\tau\wedge d\bar\tau
   &=-4n\pi_*(L_R\wedge\mathbf E_\vartheta(U_l)).
 \end{aligned}
\end{equation}
Thus $e_0=e_0^{\rm ent}+e_0^{\rm mass}+e_0^E$ and
$e_1=e_1^{\rm ent}+e_1^E$.

To estimate the entropy contributions, recall that
$B>k_0C_{\mathrm w}+2$. By the definition of $s_{l,B}$,
an elementary calculation using \eqref{eq:ale-determinant-bound} shows that
\begin{equation}\label{eq:ale-truncated-log-bound}
 |s_{l,B}|\leq C(1+B)|f_{t,l}-1|
        \leq C(1+B)|A_{t,l}-I|.
\end{equation}
The matrix bound $A_{t,l}\leq\Lambda I$ also implies
$|s_{l,B}|\leq C(1+B)$.

In the following calculations, the spatial integrals are over
$A_R$ and the time integrals are over $I$. The estimates for
the second derivatives of the cutoff give
\[
 \begin{aligned}
 \|e_0^{\rm ent}\|_{L^1(I)}
 &\leq C(1+B)R^{-2}\int_{A_R\times I}c_l\,d\mu\,dt
       \leq C(1+B)R^{2n-2\gamma-4},\\
 \|e_0^{\rm mass}\|_{L^1(I)}
 &\leq CR^{-2}\int_{A_R\times I}c_l\,d\mu\,dt
       \leq CR^{2n-2\gamma-4},
 \end{aligned}
\]
by \eqref{eq:ale-annular-time-bounds}.
For $e_1^{\rm ent}$ we keep the factor $|A_{t,l}-I|$ in
\eqref{eq:ale-truncated-log-bound}. By Cauchy--Schwarz, we have
\[
 \begin{aligned}
 \|e_1^{\rm ent}\|_{L^1(I)}
 &\leq C(1+B)R^{-1}
          \int_{A_R\times I}|A_{t,l}-I||b_l|\,d\mu\,dt\\
 &\leq C(1+B)R^{-1}
       \|A_{t,l}-I\|_{L^2(A_R\times I)}
       \|b_l\|_{L^2(A_R\times I)}\\
 &\leq C(1+B)R^{n-\gamma-2}.
 \end{aligned}
\]
Here \eqref{eq:ale-global-hessian} and $|I|\leq1$ bound the
first norm uniformly; \eqref{eq:ale-mixed-l2-bounds} bounds the
second. No uniform smallness of the first norm on the annulus
is required.

For $e_0^E$ and $e_1^E$, recall that $p^*\vartheta$ and
$\Omega$ have only spatial components. A nonzero top-degree
term in $\ddc\theta_R\wedge\mathbf E_\vartheta(U_l)$ must
therefore contain both $d\tau$ and $d\bar\tau$ from the powers
of $\Omega_l$. Their coefficients contain either $c_l$ or two
mixed entries, bounded by $|b_l|^2\leq\Lambda c_l$.
In $L_R\wedge\mathbf E_\vartheta(U_l)$, the form $L_R$
already supplies one time differential, so the coefficient
contains one mixed entry from $b_l$.
The remaining spatial coefficients are bounded in terms of
$n,\Lambda$. Since
$|U_l||\vartheta|_\omega\leq CR^{-\gamma-\delta-2}$ on $A_R$,
\[
 \begin{aligned}
 \|e_0^E\|_{L^1(I)}
 &\leq CR^{-\gamma-\delta-4}
                \int_{A_R\times I}c_l\,d\mu\,dt
       \leq CR^{2n-3\gamma-\delta-6},\\
 \|e_1^E\|_{L^1(I)}
 &\leq CR^{-\gamma-\delta-3}
                \int_{A_R\times I}|b_l|\,d\mu\,dt\\
 &\leq CR^{n-\gamma-\delta-3}
             \|b_l\|_{L^2(A_R\times I)}
       \leq CR^{2n-2\gamma-\delta-4}.
 \end{aligned}
\]
We used \eqref{eq:ale-annular-time-bounds} in the first estimate;
the second follows from Cauchy--Schwarz,
$\operatorname{Vol}(A_R)\leq CR^{2n}$, and
\eqref{eq:ale-mixed-l2-bounds}.
Because $\gamma>0$, $\delta>n-1$, and $\gamma+2-n>0$,
the two $E_\vartheta$ contributions decay faster than the
respective entropy contributions. Adding the five estimates
proves \eqref{eq:ale-uniform-error-bounds}.

\subsubsection*{The energy term and the truncation error.}
We next prove that $\E_{l,R}$ tends to zero as $R \to \infty$ uniformly in $l,t$.
We separate the scalar-curvature term by subtracting the reference
powers of $\omega$:
\[
 \begin{aligned}
 &\sum_{j=0}^{n-1}u_{t,l}\vartheta
       \wedge\omega_{u_{t,l}}^j\wedge\omega^{n-1-j}\\
 &\quad=n u_{t,l}\vartheta\wedge\omega^{n-1}
       +u_{t,l}\vartheta\wedge\sum_{j=1}^{n-1}
          \bigl(\omega_{u_{t,l}}^j-\omega^j\bigr)
                   \wedge\omega^{n-1-j}\\
 &\quad=u_{t,l}S_\omega\,\omega^n
       +O\bigl(r^{-\gamma-\delta-2}|A_{t,l}-I|\bigr)\omega^n.
 \end{aligned}
\]
For the last equality, each difference
$\omega_{u_{t,l}}^j-\omega^j$ contains a factor
$\omega_{u_{t,l}}-\omega$, whose norm is $|A_{t,l}-I|$.
The remaining form factors are uniformly bounded by
$0\leq\omega_{u_{t,l}}\leq\Lambda\omega$, and
$|u_{t,l}|\,|\vartheta|_\omega\leq Cr^{-\gamma-\delta-2}$.
Thus Cauchy--Schwarz and \eqref{eq:ale-global-hessian} yield
\begin{equation}\label{eq:ale-local-energy-tail}
 \begin{aligned}
 |\E_{l,R}(t)|
 &\leq CR^{-\gamma}\int_{\{r>R\}}|S_\omega|\,d\mu\\
 &\quad+C\left(\int_{\{r>R\}}
             r^{-2\gamma-2\delta-4}\,d\mu\right)^{1/2}
                \|A_{t,l}-I\|_{L^2(X)}\\
 &\leq C\left(R^{-\gamma}\int_{\{r>R\}}|S_\omega|\,d\mu
                    +R^{n-\gamma-\delta-2}\right)
 =\zeta_R.
 \end{aligned}
\end{equation}
The integral of the weight is finite because
$2\gamma+2\delta+4>2n$. Also $S_\omega\in L^1$ and
$n-\gamma-\delta-2<0$, so $\zeta_R\to0$ as $R\to\infty$.

Lastly, we remove the truncation by estimating the difference
\[
 \begin{aligned}
 \M^B_{l,R}(t)-\M_{l,R}(t)
   =\int_X\theta_R f_{t,l}
                  (s_{l,B}-\log f_{t,l})\,d\mu.
 \end{aligned}
\]
To bound the difference, put
$a_B=e^{-B+k_0C_{\mathrm w}}<1/2$.
Where $s_{l,B}>\log f_{t,l}$, we have $f_{t,l}<a_B$ and
\[
 0\leq f_{t,l}(s_{l,B}-\log f_{t,l})
      \leq f_{t,l}\log(a_B/f_{t,l})\leq a_B/e.
\]
The value at $f_{t,l}=0$ is by continuous extension.
The region $\{s_{l,B}>\log f_{t,l}\}$ is precisely where
truncation changes the logarithm, and $f_{t,l}<1/2$ there.
Since $h$ is decreasing on $[0,1]$ and $h(1/2)>0$, by the uniform
entropy bound in \eqref{eq:ale-common-bounds} we have
\[
 h(1/2)\operatorname{Vol}_\omega\{s_{l,B}>\log f_{t,l}\}
 \leq\int_{\{s_{l,B}>\log f_{t,l}\}}h(f_{t,l})\,d\mu
 \leq H_*.
\]
Thus the volume of the region where truncation changes the
logarithm is at most $H_*/h(1/2)$, uniformly in $l,t,B$.
The truncation error vanishes on the complement. Integrating
the bound $f_{t,l}(s_{l,B}-\log f_{t,l})\leq a_B/e$ over
that region and using $0\leq\theta_R\leq1$, we obtain
\begin{equation}\label{eq:ale-truncation-uniform}
 0\leq\M^B_{l,R}(t)-\M_{l,R}(t)
   =H^B_{l,R}(t)-H_{l,R}(t)\leq Ce^{-B}.
\end{equation}
In \eqref{eq:ale-truncation-uniform}, one may take
$C=e^{k_0C_{\mathrm w}}H_*/(e\,h(1/2))$.
Hence $C$ depends only on $k_0$, $C_{\mathrm w}$, and $H_*$,
and is independent of $l,R,t,B$.
In particular, the error does not depend on the size of the
spatial support of $U_l$.
Together with \eqref{eq:ale-local-energy-tail}, the estimate gives
\begin{equation}\label{eq:ale-q-versus-tail}
 H_{l,R}-\zeta_R\leq\M^B_{l,R}
             \leq H_{l,R}+Ce^{-B}+\zeta_R
\end{equation}
at every time.

\subsubsection*{Completion of the proof.}
To prove the ``approaching'' pointwise convex inequality in the lemma,
we first check continuity of $\M^B_{l,R}$ in $t$ for fixed $l,R,B$. 
All slices of $U_l$ are supported in $K_l$.
Let $t_j\to t$ and write
 $f_j=f_{t_j,l}$, $f=f_{t,l}$,
 $v_j=u_{t_j,l}$, $v=u_{t,l}$.
To prove continuity of $H^B_{l,R}$ in $t$, define
\[
 h_B(f,v)=f\max\{\log f,-k_0v-B\}-f+1,
 \qquad h_B(0,v)=1.
\]
Since $f_j=f=1$ and $v_j=v=0$ outside $K_l$, the integrand
of the difference is supported in $K_l$, and
\[
 |H^B_{l,R}(t_j)-H^B_{l,R}(t)|
 \leq\int_{K_l}\theta_R
       |h_B(f_j,v_j)-h_B(f,v)|\,d\mu.
\]
We therefore estimate the difference of the integrands.

For $0\leq f\leq\Lambda^n$ and $|v|\leq C_{\mathrm w}$,
\[
 h_B(f,v)=
 \begin{cases}
  (-k_0v-B-1)f+1,&0\leq f\leq e^{-k_0v-B},\\
  f\log f-f+1,&f\geq e^{-k_0v-B}.
 \end{cases}
\]
Their derivatives
with respect to $f$ are $-k_0v-B-1$ and $\log f$, respectively.
On the second branch,
$-B-k_0C_{\mathrm w}\leq\log f\leq n\log\Lambda$.
Both derivatives are therefore bounded in absolute value by
 $C_B=1+B+k_0C_{\mathrm w}+n\log\Lambda$.
Thus $h_B$ is Lipschitz in $f$ with constant $C_B$,
uniformly for $|v|\leq C_{\mathrm w}$.
By the triangle inequality,
\[
 \begin{aligned}
 |h_B(f_j,v_j)-h_B(f,v)|
 &\leq|h_B(f_j,v_j)-h_B(f,v_j)|
       +|h_B(f,v_j)-h_B(f,v)|\\
 &\leq C_B|f_j-f|+k_0\Lambda^n|v_j-v|.
 \end{aligned}
\]
Integrating the preceding estimate and using
$0\leq\theta_R\leq1$, we obtain
\[
 \begin{aligned}
 |H^B_{l,R}(t_j)-H^B_{l,R}(t)|
 &\leq C_B\|f_j-f\|_{L^1(K_l)}
       +k_0\Lambda^n\operatorname{Vol}_\omega(K_l)
              \|v_j-v\|_{L^\infty(K_l)}
 \end{aligned}
\]
Here $\operatorname{Vol}_\omega(K_l)<\infty$ and $C_B$ are
fixed because $l,B$ are fixed. Joint $C^{1,1}$ regularity implies
$v_j\to v$ uniformly as $t_j\to t$.

It remains to prove $\|f_j-f\|_{L^1(K_l)}\to0$.
Since $v_j=u_{t_j,l}\to v=u_{t,l}$ locally uniformly on $X$,
Bedford--Taylor continuity gives
$f_j\,d\mu\rightharpoonup f\,d\mu$.
By Proposition~\ref{prop:compact-continuity}, we have
$\M(v_j)\to\M(v)$, while uniform convergence of the potentials
and Bedford--Taylor continuity on the fixed compact support imply
$E_\vartheta(v_j)\to E_\vartheta(v)$.
As $f_j=f=1$ outside $K_l$ and $h(1)=0$, it follows that
\[
 \begin{aligned}
 \int_{K_l}h(f_j)\,d\mu
 &=\M(v_j)+nE_\vartheta(v_j)\\
 &\longrightarrow\M(v)+nE_\vartheta(v)
   =\int_{K_l}h(f)\,d\mu.
 \end{aligned}
\]

To obtain the $L^1$ estimate, note that
$h''(s)=1/s\geq\Lambda^{-n}$ for $0<s\leq\Lambda^n$.
Applying Taylor's inequality at $(x+y)/2$ to $x$ and $y$ and
adding the resulting inequalities, we obtain
\[
 |x-y|^2\leq8\Lambda^n
 \left[\frac{h(x)+h(y)}2-h\!\left(\frac{x+y}2\right)\right],
 \qquad 0\leq x,y\leq\Lambda^n,
\]
Applying the inequality to $f_j,f$ and using Cauchy--Schwarz yields
\[
 \begin{aligned}
 \|f_j-f\|_{L^1(K_l)}^2
 &\leq\operatorname{Vol}_\omega(K_l)
             \int_{K_l}|f_j-f|^2\,d\mu\\
 &\leq8\Lambda^n\operatorname{Vol}_\omega(K_l)
       \int_{K_l}\left[
       \frac{h(f_j)+h(f)}2
       -h\!\left(\frac{f_j+f}2\right)\right]d\mu.
 \end{aligned}
\]
The last integral is nonnegative by convexity of $h$.
Moreover, $((f_j+f)/2)\,d\mu\rightharpoonup f\,d\mu$.
The variational formula \eqref{eq:entropy-dual} gives lower
semicontinuity of entropy for the midpoint measures; since their
densities equal one outside $K_l$, we obtain
\[
 \int_{K_l}h(f)\,d\mu
 \leq\liminf_{j\to\infty}
       \int_{K_l}h\!\left(\frac{f_j+f}2\right)d\mu.
\]
Together with the convergence of $\int_{K_l}h(f_j)\,d\mu$,
the lower bound shows that the nonnegative integral in the
$L^1$ estimate tends to zero. Hence
$\|f_j-f\|_{L^1(K_l)}\to0$.
The previous bound for
$|H^B_{l,R}(t_j)-H^B_{l,R}(t)|$, together with uniform convergence
of $v_j$ to $v$, now proves continuity of $H^B_{l,R}$ in $t$.
The mixed measures in $\E_{l,R}$ converge by local
Bedford--Taylor continuity. Their common coefficient bounds
also allow multiplication by the uniformly convergent factors
$u_{t_j,l}\theta_R\vartheta$, which proves continuity of
$\E_{l,R}$. Thus $\M^B_{l,R}=H^B_{l,R}-\E_{l,R}$ is continuous.

To remove the error terms, set
\[
 \begin{aligned}
 w(t)&=\int_a^t\bigl((t-s)e_0(s)+e_1(s)\bigr)\,ds -\frac{t-a}{b-a}
       \int_a^b\bigl((b-s)e_0(s)+e_1(s)\bigr)\,ds.
 \end{aligned}
\]
Then $w$ is continuous on $[a,b]$, satisfies $w(a)=w(b)=0$,
and solves $w''=e_0+\partial_te_1$ in distributions.
By \eqref{eq:ale-uniform-error-bounds},
\[
 \|w\|_{L^\infty(a,b)}
 \leq\frac{b-a}{4}\|e_0\|_{L^1(I)}+\|e_1\|_{L^1(I)}
 \leq C(1+B)R^{n-\gamma-2}.
\]
Here $2n-2\gamma-4<n-\gamma-2<0$.
By \eqref{eq:ale-distribution-errors},
$(\M^B_{l,R}-w)''=\nu_{l,R,B}\geq0$.
The continuous function $\M^B_{l,R}-w$ is therefore convex.
Using $w(a)=w(b)=0$, we obtain
\[
 \M^B_{l,R}(t)
 \leq\frac{b-t}{b-a}\M^B_{l,R}(a)
       +\frac{t-a}{b-a}\M^B_{l,R}(b)
       +C(1+B)R^{n-\gamma-2}.
\]
Finally, by \eqref{eq:ale-q-versus-tail} at $a,t,b$, we have
\[
 \begin{aligned}
 H_{l,R}(t)
 &\leq\M^B_{l,R}(t)+\zeta_R\\
 &\leq\frac{b-t}{b-a}H_{l,R}(a)
       +\frac{t-a}{b-a}H_{l,R}(b)
       +Ce^{-B}+2\zeta_R+C(1+B)R^{n-\gamma-2}.
 \end{aligned}
\]
This is \eqref{eq:ale-tail-chord-quantitative}.
\end{proof}

\begin{proof}[Proof of Lemma~\ref{lem:uniform-entropy-tails}]
Fix $k_0=1+\|\vartheta\|_{L^\infty(X,\omega)}$.
For $R\geq R_0\geq1$, choose
\[
 B=k_0C_{\mathrm w}+3+(\gamma+2-n)\log R.
\]
Since $\gamma+2-n=n-2+\varepsilon'>0$, this choice satisfies
$B>k_0C_{\mathrm w}+2$, and we have
\[
 e^{-B}=e^{-k_0C_{\mathrm w}-3}R^{n-\gamma-2},
 \qquad 1+B\leq C(1+\log R).
\]
The first and last error terms in
\eqref{eq:ale-tail-chord-quantitative} are therefore bounded by
$C(1+\log R)R^{n-\gamma-2}$.
For the remaining term, since $S_\omega\in L^1$, we have
\[
 R^{-\gamma}\int_{\{r>R\}}|S_\omega|\,d\mu
       \leq\|S_\omega\|_{L^1}R^{-\gamma}.
\]
Since $n\geq2$, $\delta>0$, and $R\geq1$,
by \eqref{eq:ale-energy-tail-rate} we have
\[
 \zeta_R\leq C\bigl(R^{-\gamma}+R^{n-\gamma-\delta-2}\bigr)
       \leq CR^{n-\gamma-2}.
\]
Since $n-\gamma-2=-(n-2+\varepsilon')$, absorbing all three
terms into $C(1+\log R)R^{-(n-2+\varepsilon')}$ proves
\eqref{eq:ale-tail-chord} with
\eqref{eq:ale-tail-error-rate}.
The parameter $k_0$ is determined by $(X,\omega)$, so the
dependencies are exactly those stated in the lemma.
Since $\varepsilon'>0$, this error tends to zero as $R\to\infty$,
uniformly in $l$ and $t\in[a,b]$.
\end{proof}

\subsection{Convergence of entropy and proof of convexity}
\label{subsec:entropy-convergence-convexity}

We now combine convergence on compact subsets, Lemma \ref{lem:compact-density-convergence}, with uniform tail estimates,
Lemma~\ref{lem:uniform-entropy-tails} to prove convergence of the relative entropy at every interior time \(t\in(0,1)\).

\begin{proof}[Proof of Proposition~\ref{prop:interior-entropy-convergence}]
In the proof, we write 
$H_l(t)=\Ent_{\mathrm{rel}}(u_{t,l})$, $H(t)=\Ent_{\mathrm{rel}}(u_t)$, $M_l(t)=\MM(u_{t,l})$, $M(t)=\MM(u_t)$.
The functions $M_l$ are convex and continuous on $[0,1]$ by
Propositions~\ref{prop:compact-convexity} and~\ref{prop:compact-continuity}.
Lemma~\ref{lem:endpoint-energy-convergence} bounds their endpoint
values uniformly, so by convexity we have a common upper bound.
The estimates in Lemma~\ref{lem:twisted-energy-convergence} bound
$E_\vartheta(u_{t,l})$ uniformly. Since $H_l\geq0$, these estimates also imply
a common lower bound for $M_l$. Thus $|M_l(t)|\leq C$ for a
constant $C$ independent of $l,t$.
To extract a convergent subsequence, fix $[a,b]\Subset(0,1)$.
Convexity bounds the one-sided slopes by the secant slopes to
$0$ and $1$, and hence
\[
 |M_l(t)-M_l(s)|\leq\frac{2C}{\min\{a,1-b\}}|t-s|,
 \qquad s,t\in[a,b].
\]
The Arzel\`a--Ascoli theorem on compact subintervals, followed
by a diagonal argument, therefore shows that every subsequence
has a further subsequence converging locally uniformly to a
finite convex function $G$ on $(0,1)$.
It remains to prove that every such limit $G$ equals $M$.
The uniform convergence of $E_\vartheta(u_{t,l})$ will then imply
$H_l(t)\to H(t)$ for the entire sequence at every $t\in(0,1)$.
Fix such a subsequence, and assume without loss of generality that it converges to $G$.

By Lemma~\ref{lem:twisted-energy-convergence},
\begin{equation}\label{eq:ale-entropy-subsequence}
 H_l(t)=M_l(t)+nE_\vartheta(u_{t,l})
       \longrightarrow G(t)+nE_\vartheta(u_t)
       \quad(0<t<1).
\end{equation}
By lower semicontinuity of entropy, we have $M(t)\leq G(t)$.
Define the nonnegative discrepancy
\[
 D(t)=G(t)-M(t)\geq0.
\]
It measures a possible difference between the limiting entropy
values in \eqref{eq:ale-entropy-subsequence} and the entropy
of the limiting slice.

For each fixed $R$, by Lemma~\ref{lem:compact-density-convergence}
we have
\[
 \int_X(1-\theta_R)h(f_{t,l})\,d\mu
       \longrightarrow\int_X(1-\theta_R)h(f_t)\,d\mu.
\]
Subtract this convergence from
\eqref{eq:ale-entropy-subsequence}. Recalling
$H_{l,R}(t)=\int_X\theta_Rh(f_{t,l})\,d\mu$, we obtain
\begin{equation}\label{eq:ale-tail-limit}
 \lim_{l\to\infty}H_{l,R}(t)
       =D(t)+\int_X\theta_Rh(f_t)\,d\mu.
\end{equation}
The last integral tends to zero as $R\to\infty$ because
$H(t)<\infty$. Thus
\[
 \lim_{R\to\infty}\lim_{l\to\infty}H_{l,R}(t)=D(t).
\]

Fix $0<a<t<b<1$ and apply \eqref{eq:ale-tail-chord}.
First let $l\to\infty$, using \eqref{eq:ale-tail-limit} at
the three fixed times. Then let $R\to\infty$. The entropy
tails of the limiting slices vanish by their finiteness,
and $\varepsilon_R\to0$ by \eqref{eq:ale-tail-error-rate}.
We obtain
\[
 D(t)\leq\frac{b-t}{b-a}D(a)
             +\frac{t-a}{b-a}D(b).
\]
Thus $D$ is convex on $(0,1)$.

Convergence of the endpoint Mabuchi energies, proved in
Lemma~\ref{lem:endpoint-energy-convergence}, and convexity of
$M_l$ imply
\begin{equation}\label{eq:ale-limit-endpoint-chord}
 G(t)\leq(1-t)M(0)+tM(1).
\end{equation}
The joint $C^1$ bound implies $u_{t_j}\to u_t$ uniformly on $X$
when $t_j\to t$. Continuity of $t\mapsto E_\vartheta(u_t)$,
obtained as in Lemma~\ref{lem:twisted-energy-convergence}, and
lower semicontinuity of the entropy by \eqref{eq:entropy-lsc}
imply that $M(t)$ is lower semicontinuous on $[0,1]$.
Letting $t\downarrow0$ in \eqref{eq:ale-limit-endpoint-chord} yields
\[
 0\leq\liminf_{t\downarrow0}D(t)
 \leq\limsup_{t\downarrow0}D(t)
 \leq M(0)-\liminf_{t\downarrow0}M(t)\leq0.
\]
The same argument applies as $t\uparrow1$.
A nonnegative convex function with both boundary limits zero
vanishes identically.
Thus $D=0$, and \eqref{eq:ale-entropy-subsequence} implies
$H_l(t)\to H(t)$ along the subsequence. This proves that
\[
\lim_{l} H_l (t) = H (t).
\]
\end{proof}

\begin{proof}[Proof of Theorem~\ref{thm:main-b}]
By Proposition~\ref{prop:interior-entropy-convergence} and
Lemma~\ref{lem:twisted-energy-convergence}, we have
$\MM(u_{t,l})\to\MM(u_t)$ for every interior time.
Lemma~\ref{lem:endpoint-energy-convergence} gives the same
convergence at $t=0,1$. Convexity passes to the limit; in particular,
\begin{equation}\label{eq:original-chord}
 \MM(u_t)\leq(1-t)\MM(u_0)+t\MM(u_1).
\end{equation}

For $0<a<b<1$, convexity bounds every secant slope of $\MM(u_t)$
on $[a,b]$: 
\[
 \frac{\MM(u_a)-\MM(u_0)}{a}\leq\frac{\MM(u_{t})-\MM(u_s)}{t-s} \le \frac{\MM(u_1)-\MM(u_b)}{1-b},
 \qquad a\leq s<t\leq b.
\]
Hence we have local Lipschitz continuity on $(0,1)$. For continuity at the endpoints, let $t_j\to t$ with
$t\in\{0,1\}$. By the uniform bound for $\dot U$, we have
$\varepsilon_j:=\|u_{t_j}-u_t\|_\infty\leq C|t_j-t|\to0$.
The interpolation and integration-by-parts argument in the proof of
Lemma~\ref{lem:twisted-energy-convergence} applies to these time slices
and proves $E_\vartheta(u_{t_j})\to E_\vartheta(u_t)$. Local uniform convergence of the
potentials also implies local weak convergence of their Monge--Amp\`ere
measures by Bedford--Taylor continuity. By Lemma~\ref{lem:entropy-lsc},
we therefore have lower semicontinuity of the entropy, and continuity
of $E_\vartheta$ implies lower semicontinuity of $\MM$ at both
endpoints. Taking the upper limit in \eqref{eq:original-chord}
as $t\downarrow0$ gives the opposite inequality, so
\[
 \MM(u_0)\leq\liminf_{t\downarrow0}\MM(u_t)
       \leq\limsup_{t\downarrow0}\MM(u_t)\leq\MM(u_0).
\]
The argument at $t=1$ is identical. Hence the Mabuchi
functional is continuous on the closed interval.
\end{proof}

\section{Lichnerowicz equation and uniqueness of sfK metrics}
\label{sec:lichnerowicz}
\label{sec:uniqueness}

The linear ingredient in passing from convexity to uniqueness is the
construction of compactly supported approximate solutions to the
linearized scalar curvature equation, the Lichnerowicz equation. In the compact case, by Fredholm theory,
solvability is obstructed by the kernel of the Lichnerowicz operator, which consists of
holomorphy potentials. We need to establish a corresponding theory of the
Lichnerowicz equation on ALE K\"ahler manifolds.

The relation between holomorphy potentials and the linearized
scalar-curvature equation has been studied by
Lichnerowicz, Matsushima, and Calabi
\cite{Lichnerowicz57,Matsushima57,Calabi85}.
Weighted elliptic theory on noncompact manifolds was developed by
Lockhart--McOwen \cite{LockhartMcOwen}.
Arezzo--Pacard \cite[Section 5.2]{ArezzoPacard} use this theory for the
scalar operator on ALE spaces in their gluing construction. Their
holomorphic coordinates at infinity and potential expansion are
stronger than the general ALE chart used here. The argument below for
decaying Killing fields does not require a holomorphic ALE chart.

We follow the perturbative organization of
Berman--Berndtsson \cite[Section 4]{BB}. There are two additional
points at infinity. The endpoint slope inequality must be proved for
the relative functional, and the perturbation term must respect the
range obstruction in \eqref{eq:zero-flux}. We first prove
Lemma~\ref{prop:weighted-density}, which supplies compactly
supported perturbation terms. We then establish the slope
inequality by changing the reference metric and complete the
uniqueness argument. Exact solvability in weighted spaces is discussed
in Appendix~\ref{app:weighted-lichnerowicz}.

\subsection{The Lichnerowicz equation in ALE K\"ahler manifolds}
In this subsection only, we use $\omega$ to denote a smooth complete
scalar-flat K\"ahler form on the connected manifold $X$. Write
$g_\omega$ for its associated Riemannian metric. We work in the ALE
chart $\Upsilon$ fixed in the introduction, with the same radius
function $r$.
For $v\in C_c^\infty(X,\RR)$, set $\omega_s=\omega+sdd^cv$ for
sufficiently small $s$, and write
$\Delta_\omega v=\omega^{a\bar b}v_{a\bar b}$.
Differentiating the inverse metric and
$\Ric(\omega_s)=-\ddc\log\det((\omega_s)_{a\bar b})$,
we have
\[
 \left.\frac{d}{ds}\right|_{s=0}\omega_s^{a\bar b}
 =-\omega^{a\bar d}v_{c\bar d}\omega^{c\bar b},
 \qquad
 \left.\frac{d}{ds}\right|_{s=0}\Ric(\omega_s)
 =-dd^c(\Delta_\omega v).
\]
Taking the variation of the trace of the Ricci form, we obtain
\begin{equation}\label{eq:lich-formula}
 \left.\frac{d}{ds}\right|_{s=0}S_{\omega_s}
 =-\Delta_\omega^2v-
       \langle\Ric(\omega),dd^cv\rangle_\omega
 =:-L_\omega v.
\end{equation}

To identify $L_\omega$, define
$\mathcal D_\omega v=\bar\partial\nabla^{1,0}_\omega v$,
with formal adjoint taken using the Hermitian metric induced by
$\omega$ and the volume form $\omega^n$.
Commuting covariant derivatives gives
\[
 \mathcal D_\omega^*\mathcal D_\omega v
 =\Delta_\omega^2v+
   \langle\Ric(\omega),dd^cv\rangle_\omega
   +\omega^{a\bar b}(\partial_aS_\omega)\partial_{\bar b}v.
\]
Since $S_\omega=0$, the last term vanishes. Thus
$L_\omega=\mathcal D_\omega^*\mathcal D_\omega$ is the Lichnerowicz
operator, and the scalar-curvature linearization is $-L_\omega$.
In the following, we also use the linearization of the scalar-curvature measure.
Since $S_\omega=0$, variation of the volume form contributes no term,
and by \eqref{eq:lich-formula}, we have
\begin{equation}\label{eq:scalar-density-linearization}
 \left.\frac{d}{ds}\right|_{s=0}
       (-S_{\omega_s})\omega_s^n
       =(L_\omega v)\,\omega^n .
\end{equation}

We first prove that a Killing field which decays at infinity vanishes,
and then apply this result to the kernel of $L_\omega$.
For these two vanishing results, it suffices to assume the
differentiated ALE estimates
\begin{equation}\label{eq:finite-order-linear-ale}
 |\nabla_{\mathrm{euc}}^{\,j}(g_\omega-g_{\mathrm{euc}})|
 \leq C_jr^{-\delta_\omega-j},\qquad
 0\leq j\leq4,\qquad \delta_\omega>0.
\end{equation}
The estimates are understood in the fixed ALE chart $\Upsilon$,
lifted to its Euclidean cover. In particular,
$|\mathrm{Rm}(g_\omega)|\leq Cr^{-2-\delta_\omega}$.

\begin{lemma}\label{lem:killing}
Let $Z$ be a Killing field on a connected complete ALE manifold
satisfying \eqref{eq:finite-order-linear-ale}. If
$|Z|+|\nabla Z|\to0$ at infinity, then $Z=0$.
\end{lemma}
\begin{proof}
Fix $x\in X$. We first construct a minimizing ray from $x$ which
leaves every compact set. Choose a sequence $\{x_j\}$ with
$d_\omega(x,x_j)\to\infty$, and let $\sigma_j$ be a minimizing
unit-speed geodesic from $x$ to $x_j$. The initial velocities lie in
the compact unit sphere of $T_xX$. After passing to a subsequence,
$\dot\sigma_j(0)\to v$ with $|v|=1$.
Completeness defines $\sigma(t)=\exp_x(tv)$ for every $t\geq0$.
For each fixed $t$, the segment $\sigma_j|_{[0,t]}$ is defined and
minimizing for large $j$. By continuous dependence on the initial
velocity, we have
\[
 d_\omega(x,\sigma(t))
 =\lim_{j\to\infty}d_\omega(x,\sigma_j(t))=t.
\]
Thus $\sigma$ is a minimizing ray. Since distance from $x$ is bounded
on any compact set, $\sigma$ eventually leaves every compact set.

We next compare $r(\sigma(t))$ with $t$. Enlarge $R_0$ so that the
compact region bounded by $\{r=R_0\}$ contains $x$, and
\[
 \tfrac12g_{\mathrm{euc}}\leq g_\omega\leq2g_{\mathrm{euc}}
 \qquad\text{on }\{r\geq R_0\}.
\]
For a point $y$ on this end, every path from $x$ to $y$ crosses
$\{r=R_0\}$. After the last crossing, the Euclidean length is at
least $r(y)-R_0$. Conversely, let $m_y$ be the point at $\{r=R_0\}$ such that $y$ can be connected to $m_y$ by a radial line in Euclidean coordinates of the end. 
Join $x$ to $m_y$, and then follow the radial line to $y$. By compactness of the coordinate sphere and metric
comparison, there is a constant $C_x$ such that
\[
 \frac{r(y)-R_0}{\sqrt2}
 \leq d_\omega(x,y)
 \leq C_x+\sqrt2\bigl(r(y)-R_0\bigr).
\]
Setting $y=\sigma(t)$ yields
\[
 R_0+\frac{t-C_x}{\sqrt2}
 \leq r(\sigma(t))\leq R_0+\sqrt2t.
\]
Consequently, after increasing $C_x$, we have
$C_x^{-1}t\leq r(\sigma(t))\leq C_xt$ for sufficiently large $t$.
Below, constants denoted by $C_x$ may increase and depend on $x$,
$(X,\omega)$ and the bounds in \eqref{eq:finite-order-linear-ale}.

Write $Z(t)=Z(\sigma(t))$. To obtain an equation for $Z(t)$, let
$\Phi_s$ be the local flow generated by the Killing field $Z$.
By the Killing equation, we have
\[
 \frac{d}{ds}\Phi_s^*g_\omega
 =\Phi_s^*(\mathcal L_Zg_\omega)=0,
\]
so $\Phi_s$ is a local isometry. On each finite segment of $\sigma$,
$F(s,t)=\Phi_s(\sigma(t))$ is therefore a geodesic in $t$ for each
sufficiently small $s$. The variation field is
$\partial_sF(0,t)=Z(t)$; hence, $Z(t)$ satisfies the Jacobi equation:
\[
 D_t^2Z(t)+R_\omega(Z(t),\dot\sigma(t))\dot\sigma(t)=0,
\]
where $D_t$ denotes covariant differentiation along $\sigma$ and
$R_\omega$ is the curvature tensor. Choose a parallel orthonormal
frame along $\sigma$, and continue to write $Z(t)$ for its component
vector in that frame. Covariant differentiation then becomes
ordinary differentiation. By the curvature bound and the radius
comparison, we have
\[
 |Z''(t)|\leq |\mathrm{Rm}(g_\omega)|(\sigma(t))\,|Z(t)|
 \leq C_xt^{-2-\delta_\omega}|Z(t)|
\]
for large $t$. Moreover, by the assumed decay, we have
\[
 Z(t)\longrightarrow0,\qquad Z'(t)\longrightarrow0.
\]

It suffices to show that this solution of the Jacobi equation
vanishes. For sufficiently large $T$, set
$M_T=\sup_{s\geq T}|Z(s)|<\infty$.
Integrating from infinity and using $Z'(t)\to0$ and $Z(t)\to0$, we have, for $t\geq T$,
\[
 Z(t)=-\int_t^\infty Z'(s)\,ds
     =\int_t^\infty(s-t)Z''(s)\,ds.
\]
It follows that
\[
 |Z(t)|
 \leq C_xM_T\int_t^\infty(s-t)s^{-2-\delta_\omega}\,ds =\frac{C_xM_T}{\delta_\omega(1+\delta_\omega)}
                         t^{-\delta_\omega}.
\]
Taking the supremum over $t\geq T$ and choosing $T$ sufficiently
large, we obtain
\[
 M_T\leq\frac{C_x}{\delta_\omega(1+\delta_\omega)}
                 T^{-\delta_\omega}M_T\leq\tfrac12M_T.
\]
Hence $M_T=0$, so $Z(t)=0$ for all $t\geq T$.
In particular, $Z(T)=Z'(T)=0$.
Uniqueness for the linear Jacobi equation then implies $Z(t)=0$ also
for $0\leq t\leq T$. Thus $Z(x)=0$.
Since $x$ was arbitrary, the Killing field vanishes on $X$.
\end{proof}

The identity $L_\omega=\mathcal D_\omega^*\mathcal D_\omega$ now
gives a vanishing result for potentials. 

\begin{lemma}\label{lem:decaying-kernel}
Let $(X,J,\omega)$ be connected, complete, scalar-flat, and ALE of
complex dimension $n\geq2$, satisfying
\eqref{eq:finite-order-linear-ale}. If a real distributional solution
of $L_\omega h=0$ satisfies
\[
 |h|\leq C_hr^{-\beta},\qquad \beta>n-2,
\]
then $h=0$.
\end{lemma}
\begin{proof}
By local elliptic regularity, $h$ is smooth. From the bound
$|h|\leq C_hr^{-\beta}$, we obtain the weighted derivative estimates
\begin{equation}\label{eq:kernel-derivatives}
 |\nabla h|\leq CC_hr^{-\beta-1},\qquad
 |\nabla^2h|+|\mathcal D_\omega h|
       \leq CC_hr^{-\beta-2}
\end{equation}
on the ALE end. To see this, work on the ALE cover and rescale
$B(x_0,R/2)$, where $|x_0|=R$, by setting
$h_R(y)=h(x_0+Ry)$ for $y\in B_{1/2}$.
After multiplying the equation by $R^4$,
\eqref{eq:finite-order-linear-ale} implies uniform ellipticity and
uniform $C^1$ bounds for the rescaled coefficients.
For a fixed $p>2n$, the interior $L^p$ estimates
\cite[Section~6.2]{Morrey66} give
\[
 \|h_R\|_{W^{4,p}(B_{1/4})}
 \leq C\|h_R\|_{L^p(B_{1/2})}
 \leq CC_hR^{-\beta}.
\]
By Sobolev embedding and rescaling back, we obtain
\eqref{eq:kernel-derivatives}.
Here and below, $C$ depends on $n,p,\beta,\delta_\omega$, the bounds
in \eqref{eq:finite-order-linear-ale}, the uniform comparison bounds
between $g_\omega$ and $g_{\mathrm{euc}}$ on the end, and the cutoff
profile, but not on $x_0,R$, or $h$.

Choose $0\leq\chi_R\leq1$, equal to one on $\{r<R\}$ and zero on
$\{r>2R\}$, with $|\nabla\chi_R|\leq CR^{-1}$ and
$|\nabla^2\chi_R|\leq CR^{-2}$.
Testing the equation against $\chi_Rh$ and taking the real part
of the Hermitian pairing, we have
\[
 \int_X\chi_R|\mathcal D_\omega h|^2\,\omega^n
 =-\operatorname{Re}\int_X
 \left\langle\mathcal D_\omega h,
 \mathcal D_\omega(\chi_Rh)-\chi_R\mathcal D_\omega h
 \right\rangle\,\omega^n.
\]
By the product rule, the second factor vanishes outside
$\{R<r<2R\}$ and has norm at most
$C(R^{-2}|h|+R^{-1}|\nabla h|)$ there.
Using \eqref{eq:kernel-derivatives} and the annular volume bound,
we obtain
\begin{equation}\label{eq:kernel-cutoff}
 0\leq\int_X\chi_R|\mathcal D_\omega h|^2\,\omega^n
 \leq CC_h^2R^{2n-2\beta-4}\longrightarrow0.
\end{equation}
Since $\chi_R=1$ on every fixed compact subset for large $R$,
we conclude that $\mathcal D_\omega h=0$ everywhere.

Thus, $\nabla^{1,0}_\omega h$ is holomorphic. Since $h$ is real,
$Z=J\nabla h$ is a real holomorphic vector field and
$\iota_Z\omega=-dh$. Hence, $\mathcal L_Z\omega=0$ and
$\mathcal L_ZJ=0$, which implies that $Z$ is Killing.
By \eqref{eq:kernel-derivatives}, $|Z|+|\nabla Z|\to0$.
By Lemma~\ref{lem:killing}, we have $\nabla h=0$.
Connectedness and the decay of $h$ then imply $h=0$.
\end{proof}

We use the growth convention for weighted spaces recalled in
Section~\ref{sec:preliminary}. More precisely, fix $0<\alpha<1$
and a compact set $K_*\supset\{r\leq2R_0\}$. For an integer $m\geq0$,
the space $C^{m,\alpha}_\lambda(X)$ consists of functions with finite norm
\[
 \|v\|_{C^{m,\alpha}_\lambda}
 =\|v\|_{C^{m,\alpha}(K_*)}
 +\sup_{R\geq R_0}R^{-\lambda}
   \|\widetilde v(R\,\cdot)\|_{C^{m,\alpha}(\{1<|y|<2\})},
\]
where $\widetilde v$ is the lift to the ALE end and the compact
norm is defined in a fixed finite coordinate cover. In particular,
$|\nabla^k v|=O(r^{\lambda-k})$ for $0\leq k\leq m$.

We solve the equation with weighted estimates, allowing controlled
growth when $n=2$, and then cut off the solution. The following lemma
gives the compactly supported approximations needed for uniqueness.
The full description of solvability at other weights is left to
Appendix~\ref{app:weighted-lichnerowicz}.

\begin{lemma}
\label{prop:weighted-density}
Under the assumptions of Lemma~\ref{lem:decaying-kernel}, let
$q\in C_c^\infty(X,\RR)$. There is a smooth real solution of
$L_\omega v=-q$ with
\begin{equation}\label{eq:compact-source-solution}
 \begin{aligned}
 &|\nabla^jv|\leq C_qr^{4-2n-j},
       &&0\leq j\leq4,\quad n\geq3,\\
 &|v|\leq C_q(1+\log r),
       &&n=2,\\
 &|\nabla^jv|\leq C_qr^{-j},
       &&1\leq j\leq4,\quad n=2.
 \end{aligned}
\end{equation}
The constant $C_q$ depends on $n,(X,\omega),r$, the finite ALE
bounds in \eqref{eq:finite-order-linear-ale}, and $q$.
Moreover, for every $\beta>0$ and every positive sequence
$\epsilon_j\to0$, there are $v_j\in C_c^\infty(X,\RR)$ such that
\[
 \int_Xr^{-\beta}|L_\omega v_j+q|\,\omega^n\leq\epsilon_j.
\]
\end{lemma}
\begin{proof}
Fix $\sigma=\tfrac12\min\{1,\delta_\omega\}$ and set
$\lambda=2-n+\sigma$. We first use the weighted Fredholm
alternative to solve $L_\omega v=-q$ at growth weight $\lambda$.
By the finite-order bounds \eqref{eq:finite-order-linear-ale}, we have,
on the ALE cover,
\[
 L_\omega=c\Delta_{\mathrm{euc}}^2+
       \sum_{1\leq|\kappa|\leq4}a_\kappa\partial^\kappa,
 \qquad c>0,
\]
with
\[
 |\partial^\zeta a_\kappa|
 \leq C r^{|\kappa|-4-\delta_\omega-|\zeta|},
 \qquad |\zeta|\leq1.
\]
After rescaling an annulus of radius $R$ and multiplying the
operator by $R^4$, its coefficients converge in $C^1$ to those
of $c\Delta_{\mathrm{euc}}^2$. The constant $C$ depends on
$n,\delta_\omega$, the bounds in \eqref{eq:finite-order-linear-ale},
and the uniform comparison of $g_\omega$ with $g_{\mathrm{euc}}$,
but not on $R$. The weighted Fredholm theorem of
\cite[Theorem~6.2]{LockhartMcOwen} therefore applies. There are
no exceptional weights in $(4-2n,0)$ for $n\geq3$, or in $(0,1)$
for $n=2$, and $\lambda$ lies in the corresponding interval.

Formal self-adjointness places the obstruction at complementary
weight $4-2n-\lambda=2-n-\sigma$
\cite[Section~7]{LockhartMcOwen}. By rescaled elliptic estimates, we have
$h=O(r^{2-n-\sigma+\eta})$ for every sufficiently small $\eta>0$.
Taking $\eta<\sigma$, Lemma~\ref{lem:decaying-kernel} makes the
obstruction space zero. The Fredholm alternative gives a real
solution $v$, and local elliptic regularity makes it smooth.
By weighted Sobolev embedding and rescaled Schauder estimates, we have
$v\in C^{4,\alpha}_\lambda$ for a fixed $0<\alpha<1$.

We now use the compact support of $q$ to improve the decay. On
the end, $L_\omega v=0$, so the coefficient bounds imply
\begin{equation}\label{eq:compact-source-model-error}
 c\Delta_{\mathrm{euc}}^2v
 =(c\Delta_{\mathrm{euc}}^2-L_\omega)v
       \in C^{0,\alpha}_{\lambda-4-\delta_\omega}.
\end{equation}
For any nonexceptional weight $\lambda-\delta_\omega<\nu<\lambda$,
extend the right-hand side of \eqref{eq:compact-source-model-error},
after multiplication by $r^4$, to the full cylinder in the variable
$\log r$. The weighted inverse of the model operator
$r^4c\Delta_{\mathrm{euc}}^2$ gives a particular solution at weight
$\nu$ on the end. Subtracting it from $v$ leaves a biharmonic
function there.
Its spherical harmonic expansion contains terms $r^\rho Y_d$
with
\[
 \rho\in\{d,d+2,2-2n-d,4-2n-d\},
\]
where $Y_d$ is a $\Gamma$-invariant spherical harmonic of degree
$d$; repeated roots also give logarithmic terms. Keeping the
terms with $\nu<\rho<\lambda$ leaves a remainder in
$C^{4,\alpha}_\nu$. This is the model weighted-inverse and
expansion argument of \cite[Sections~2 and~5]{LockhartMcOwen},
applied in the variable $\log r$.

For $n\geq3$, there are no roots in $(4-2n,0)$. Repeating the
preceding improvement reaches a weight
$4-2n<\lambda<\min\{0,4-2n+\delta_\omega/2\}$.
We may then cross the simple root $4-2n$, obtaining, for some
$0<\eta<\min\{1,\delta_\omega/2\}$,
\[
 v=c_0r^{4-2n}+O(r^{4-2n-\eta}).
\]
For $n=2$, the initial choice $\sigma<\delta_\omega$ permits
crossing the double root zero, and we obtain
\[
 v=c_0\log r+c_1+O(r^{-\eta}),
 \qquad 0<\eta<\min\{1,\delta_\omega-\sigma\}.
\]
The remainder estimates hold with derivatives through order four.
Thus \eqref{eq:compact-source-solution} follows, with bounds on
the compact part supplied by smoothness. The improvement uses
only the finite coefficient bounds above.

For the compactly supported approximation, let $\chi_R$ be the
usual cutoff with $|\nabla^a\chi_R|\leq C_aR^{-a}$ for
$1\leq a\leq4$. Once $\chi_R=1$ on $\supp q$,
\[
 L_\omega(\chi_Rv)+q
       =L_\omega(\chi_Rv)-\chi_RL_\omega v,
\]
which is supported on $\{R<r<2R\}$. By the product rule and
\eqref{eq:compact-source-solution}, the right-hand side is bounded
by $C_qR^{-2n}$ for $n\geq3$, and by
$C_qR^{-4}(1+\log R)$ for $n=2$. The lower-order terms satisfy
the same bounds. Thus, by annular volume growth, we have
\begin{equation}\label{eq:compact-source-cutoff}
 \int_Xr^{-\beta}|L_\omega(\chi_Rv)+q|\,\omega^n
 \leq
 \begin{cases}
 C_qR^{-\beta},&n\geq3,\\
 C_q(1+\log R)R^{-\beta},&n=2.
 \end{cases}
\end{equation}
Here $C_q$ also depends on $\beta$ and the cutoff profile, but
not on $R$. Both bounds tend to zero for every $\beta>0$.
Choose $R_j\to\infty$ so that the bound is at most $\epsilon_j$,
and set $v_j=\chi_{R_j}v$.
\end{proof}

\noindent\textit{Remark.}
Suppose $q\in C_c^\infty(X,\RR)$ and
$\int_Xq\,\omega^n\ne0$. A solution of $L_\omega v=q$ cannot
decay faster than $r^{4-2n}$ in the weighted spaces above, and
in particular cannot be compactly supported.
Indeed, if $v\in C^{4,\alpha}_\lambda$ with $\lambda<4-2n$,
then $L_\omega v\in L^1$, and by formal self-adjointness, we have
\begin{equation}\label{eq:zero-flux}
 \int_XL_\omega v\,\omega^n
 =\lim_{R\to\infty}\int_X\chi_RL_\omega v\,\omega^n
 =\lim_{R\to\infty}\int_XvL_\omega\chi_R\,\omega^n=0.
\end{equation}
Here $L_\omega1=0$, so $L_\omega\chi_R$ is supported on
$\{R<r<2R\}$ and bounded by $CR^{-4}$.
The final integral is bounded in absolute value by
$C\|v\|_{C^{4,\alpha}_\lambda}R^{2n+\lambda-4}\to0$,
contradicting $\int_Xq\,\omega^n\ne0$.
Thus the exact solution in Lemma~\ref{prop:weighted-density}
cannot in general be compactly supported, which is why we use
compactly supported approximate solutions.
The corresponding Fredholm obstruction is also described in
Proposition~\ref{prop:fredholm} in
Appendix~\ref{app:weighted-lichnerowicz}.

\subsection{Perturbation of the Mabuchi functional near scalar-flat metrics}
We return to the reference form $\omega$ of
Section~\ref{sec:preliminary}. We perturb the Mabuchi functional
by a sequence of compactly supported densities, chosen to be
positive on successively larger compact subsets of $X$. More precisely,
choose relatively compact domains 
\[
K_1\Subset K_2\Subset\cdots \Subset \cdots \Subset X
\]
exhausting $X$, and smooth nonnegative densities $d\mu_i$ supported
in $K_{i+1}$ and strictly positive near $\overline K_i$. Set
\[
 \mathcal F_{\mu_i}(v)=\int_Xv\,d\mu_i,\qquad
 \MM_{s,\mu_i}(v)=\MM(v)+s\mathcal F_{\mu_i}(v),\qquad s>0.
\]
The added functional is finite and convex along weak geodesics,
since $t\mapsto u_t(x)$ is convex for each $x$. We work with one domain and density at a
time, denoted by $K$ and $\mu$; no limit of the densities is needed in this subsection.
For a smooth potential $v$ and a smooth, compactly supported variation $\xi$,
the first variation is
\[
 \delta_\xi\MM_{s,\mu}(v)
 =\int_X\xi\bigl(-S_{\omega_v}\omega_v^n+s\,d\mu\bigr).
\]

We perturb the Mabuchi functional near the scalar-flat endpoint potentials.
Since $d\mu_i$ is strictly positive near $\overline K_i$, the
increase of the first derivative of $\mathcal F_{\mu_i}$ along a
weak geodesic controls the spatial gradient of the difference
of its endpoint potentials on $K_i$. We construct nearby endpoint
potentials for which the perturbed first variation is small
relative to the perturbation parameter. Letting that parameter
tend to zero then forces the original metrics to agree on $K_i$;
uniqueness on $X$ follows by exhaustion. The main additional work
on an ALE manifold is to control the errors from the compactly
supported approximations relative to the small perturbation.

Let $u_0,u_1\in\HH_{-\gamma}(X,\omega)$ be smooth, with
$\omega_i=\omega+dd^cu_i$ scalar-flat for $i=0,1$, and fix a
smooth nonnegative compactly supported density $d\mu$.
At $u_i$, the perturbed first variation is
$\delta_\xi\MM_{s,\mu}(u_i)=s\int_X\xi\,d\mu$;
the error need not tend to zero after division by $s$.
We therefore also perturb the endpoint potentials. To construct
the perturbation terms $v_{i,j}$, set $\beta=\gamma/2>n-2$ and
\[
 L_i=L_{\omega_i},\qquad q_i=\frac{d\mu}{\omega_i^n}
       \in C_c^\infty(X,\RR),\qquad i=0,1.
\]
Lemma~\ref{prop:weighted-density}, applied to each $\omega_i$,
gives a smooth real solution $v_i$ of $L_iv_i=-q_i$ with the
weighted derivative bounds in \eqref{eq:compact-source-solution}.
Given a positive sequence $\epsilon_j\downarrow0$, use the cutoff
$\chi_R$ from that lemma and take $R_j\uparrow\infty$ with
$\operatorname{supp}\mu\subset\{r<R_1\}$ and set
\[
 v_{i,j}=\chi_{R_j}v_i\in C_c^\infty(X,\RR),\qquad
 u_i^j=u_i+s_jv_{i,j}.
\]
By \eqref{eq:compact-source-cutoff}, the radii can be chosen for
both $i=0,1$ so that
$\int_Xr^{-\beta}|L_iv_{i,j}+q_i|\,\omega_i^n\leq\epsilon_j$.
We will choose $s_j>0$ so that the first-variation error is small
compared with $s_j$. When $u_i^j\in\HH_{-\gamma}(X,\omega)$,
let $U^j$ denote the weak geodesic joining $u_0^j$ to $u_1^j$.
Write $U_t^j(x)=U^j(x,t)$ for its time slices; dots denote time
derivatives, taken one-sided at $t=0,1$.

\begin{proposition}\label{prop:perturbed-mabuchi}
For every positive sequence $\epsilon_j\downarrow0$, the
perturbation terms $v_{i,j}$ and parameters $s_j\downarrow0$ can be
chosen so that the perturbed forms have a common positive
lower bound relative to $\omega$, and
\begin{equation}\label{eq:almost-critical}
 \left|\int_X\dot U_i^j
       \bigl(-S_{\omega_{u_i^j}}\omega_{u_i^j}^n+s_j\,d\mu\bigr)\right|
       \leq Cs_j\epsilon_j,\qquad i=0,1.
\end{equation}
Consequently,
\[
 \bigl(\MM_{s_j,\mu}(U^j)\bigr)'_+(0)\geq-Cs_j\epsilon_j,
 \qquad
 \bigl(\MM_{s_j,\mu}(U^j)\bigr)'_-(1)\leq Cs_j\epsilon_j.
\]
Here $C$ depends on $n,\gamma,k,(X,\omega),r,\mu$, the weighted
$C^k$ norms of $u_0,u_1$, and a common positive lower bound for
$\omega_0,\omega_1$ relative to $\omega$; it is independent of $j$.
\end{proposition}

The proof has two parts. First, the compactly supported perturbation terms
from Lemma~\ref{prop:weighted-density} cancel the first variation
up to a small weighted integral error. The uniform bound for
$\dot U^j$ in Section~\ref{subsec:weighted-time-derivatives} then permits
pairing that error with the endpoint velocities. Second, we need
to relate those pairings to the one-sided derivatives of $\MM$.
The next two lemmas justify changing the reference form and give
the required endpoint inequalities. We then complete the
perturbation argument.

\begin{lemma}\label{prop:cocycle}
Let $v\in\HH_{-\gamma}(X,\omega)$ be smooth, let
$\omega_v=\omega+dd^cv$, and let $u$ be a slice of a weak geodesic with
\[
 0\leq\omega_u\leq\Lambda\omega,\quad
 |u|\leq Cr^{-\gamma},\quad
 dd^cu\in L^2,\quad \Ent_{\mathrm{rel},\omega}(u)<\infty.
\]
Then all terms below are finite and
\begin{equation}\label{eq:cocycle}
 \MM_\omega(u)-\MM_\omega(v)=\MM_{\omega_v}(u-v).
\end{equation}
The subscript specifies the reference form in the definition of $\MM$.
\end{lemma}
\begin{proof}
Put $\vartheta=\Ric(\omega)$ and $\vartheta_v=\Ric(\omega_v)$;
the superscript on $E$ specifies its reference form.
Expanding the relative entropies and integrating by parts, we have
\begin{equation}\label{eq:entropy-cocycle}
 \begin{aligned}
 H_\omega(u)-H_\omega(v)-H_{\omega_v}(u-v)
 &=\int_X\log\frac{\omega_v^n}{\omega^n}
                 (\omega_u^n-\omega_v^n)\\
 &=nE^{\omega_v}_{\vartheta-\vartheta_v}(u-v),
 \end{aligned}
\end{equation}
since $dd^c\log(\omega_v^n/\omega^n)=\vartheta-\vartheta_v$.
For $u_s=v+s(u-v)$, differentiating \eqref{eq:twisted-energy} yields
\[
 \frac{d}{ds}E^\omega_\vartheta(u_s)
 =\int_X(u-v)\vartheta\wedge\omega_{u_s}^{n-1}
 =\frac{d}{ds}E^{\omega_v}_\vartheta\bigl(s(u-v)\bigr).
\]
Integrating over $0\leq s\leq1$, we obtain
\begin{equation}\label{eq:energy-transgression}
 E^\omega_\vartheta(u)-E^\omega_\vartheta(v)
       =E^{\omega_v}_\vartheta(u-v).
\end{equation}

Perform the integrations by parts on $\{r<R\}$, whose boundary
$S_R=\{r=R\}$ has area $O(R^{2n-1})$. Local smooth approximation
justifies the calculation for the weak slice for almost every
sufficiently large $R$. By the weighted gradient estimate
\eqref{eq:weighted-gradient-bound} and the decay of $v$, we have
$|\nabla u|+|\nabla v|=O(r^{-\gamma-1})$; moreover,
\[
 \log\frac{\omega_v^n}{\omega^n}=O(r^{-\gamma-2}),\qquad
 \left|\nabla\log\frac{\omega_v^n}{\omega^n}\right|
       =O(r^{-\gamma-3}).
\]
The energy boundary integrands are bounded by
$C\bigl(|v||\nabla(u-v)|+|u-v||\nabla v|\bigr)|\vartheta|$.
For the entropy, the corresponding terms pair the logarithm
with $\nabla(u-v)$, or its gradient with $u-v$.
Since the remaining mixed form factors are bounded, the two
boundary integrals are bounded, respectively, by
\begin{align}
 CR^{2n-1}R^{-\gamma}R^{-\gamma-1}R^{-\delta-2}
 &=CR^{2n-2\gamma-\delta-4}\longrightarrow0,
       \label{eq:transgression-error}\\
 CR^{2n-1}\bigl(R^{-\gamma-2}R^{-\gamma-1}
               +R^{-\gamma}R^{-\gamma-3}\bigr)
 &\leq CR^{2n-2\gamma-4}\longrightarrow0,
       \label{eq:entropy-transgression-error}
\end{align}
with $C$ independent of $R$. We may therefore let $R\to\infty$.
Absolute convergence follows from $S_\omega\in L^1$, the decay
of $v$, and the stated weighted and $L^2$ bounds.
Thus the identities hold for the weak slice $u$.
Subtracting $n$ times \eqref{eq:energy-transgression} from
\eqref{eq:entropy-cocycle} proves \eqref{eq:cocycle}.
\end{proof}

\begin{lemma}\label{prop:endpoint-slopes}
Let $U(x,t)=u_t(x)$ be the weak geodesic joining smooth endpoint
potentials $u_0,u_1\in\HH_{-\gamma}(X,\omega)$.
If $S_{\omega_{u_0}}$ and $S_{\omega_{u_1}}$ have compact support, then
\begin{align}
 (\MM(U))'_+(0)&\geq
       -\int_X\dot u_0S_{\omega_{u_0}}\,\omega_{u_0}^n,
       \label{eq:initial-slope}\\
 (\MM(U))'_-(1)&\leq
       -\int_X\dot u_1S_{\omega_{u_1}}\,\omega_{u_1}^n.
       \label{eq:final-slope}
\end{align}
The derivatives on the left are one-sided convex derivatives.
\end{lemma}
\begin{proof}
Use $\omega_{u_0}$ as the reference form and put
$v_t=u_t-u_0$ and $\vartheta_0=\Ric(\omega_{u_0})$.
By \eqref{eq:energy-first-derivative} with this reference form, we have
\[
 \left|\frac{nE^{\omega_{u_0}}_{\vartheta_0}(v_t)}{t}
       -\int_X\frac{v_t}{t}S_{\omega_{u_0}}\,\omega_{u_0}^n\right|
 \leq C\int_X\frac{|\partial v_t|^2}{t}
                  |\vartheta_0|\,\omega_{u_0}^n.
\]
Joint $C^{1,1}$ regularity and \eqref{eq:weighted-gradient-bound}
imply $|\partial v_t|\leq C\min\{t,r^{-\gamma-1}\}$ for $0<t\leq1$.
The mixed form factors are bounded, and
\begin{equation}\label{eq:slope-tail}
 \frac{|\partial v_t|^2}{t}|\vartheta_0|
 \leq Cr^{-\gamma-1}|\vartheta_0|
       \in L^1(X,\omega_{u_0}^n),
\end{equation}
with $C$ independent of $t$. The integrability follows from
$\gamma+\min\{\delta,\gamma+2\}+3>2n$.
The gradient terms tend pointwise to zero, while
$v_t/t\to\dot u_0$ uniformly on the compact support of
$S_{\omega_{u_0}}$. By dominated convergence, the
boundary derivative is
\[
 n\left.\frac{d}{dt}\right|_{0+}
 E^{\omega_{u_0}}_{\vartheta_0}(v_t)
 =\lim_{t\downarrow0}\frac{nE^{\omega_{u_0}}_{\vartheta_0}(v_t)}{t}
 =\int_X\dot u_0S_{\omega_{u_0}}\,\omega_{u_0}^n.
\]

Since $v_0=0$, the relative entropy satisfies
$H_{\omega_{u_0}}(v_0)=0$ and $H_{\omega_{u_0}}(v_t)\geq0$.
Thus, by Lemma~\ref{prop:cocycle}, we have
\begin{equation}\label{eq:slope-lower-bound}
 \frac{\MM_\omega(u_t)-\MM_\omega(u_0)}{t}
 \geq-\frac{nE^{\omega_{u_0}}_{\vartheta_0}(v_t)}{t}.
\end{equation}
Letting $t\downarrow0$ proves \eqref{eq:initial-slope}.
At $t=1$, use $\omega_{u_1}$ as reference and $v_t=u_t-u_1$.
The same computation based at $t=1$ reverses the inequality and proves \eqref{eq:final-slope}.
\end{proof}

The endpoint inequalities are the relative ALE version of
\cite[Lemma~3.5]{BB}. The estimates
\eqref{eq:transgression-error}, \eqref{eq:entropy-transgression-error},
and \eqref{eq:slope-tail} justify the additional limits at infinity.

\begin{proof}[Proof of Proposition~\ref{prop:perturbed-mabuchi}]
Set $\beta=\gamma/2$, $q_i=d\mu/\omega_i^n$, and $L_i=L_{\omega_i}$.
Apply Lemma~\ref{prop:weighted-density} at each scalar-flat metric.
The cutoff construction in its proof gives radii $R_j\uparrow\infty$
and real functions $v_{i,j}\in C_c^\infty(X)$ supported in
$\{r\leq2R_j\}$ such that
\begin{equation}\label{eq:approximate-corrections}
 \int_Xr^{-\beta}|L_iv_{i,j}+q_i|\,\omega_i^n\leq\epsilon_j,
 \qquad i=0,1.
\end{equation}
We choose $R_1$ so that $\operatorname{supp}\mu\subset\{r<R_1\}$.
Thus $v_{i,j}$ are supported in an exhaustion of compact
subsets, while the density $\mu$ remains fixed.

For each fixed $j$, taking the Taylor expansion of the scalar-curvature
measure at $u_i$ and using \eqref{eq:scalar-density-linearization}, we have
\begin{equation}\label{eq:almost-critical-expansion}
 -S_{\omega_{u_i+sv_{i,j}}}\omega_{u_i+sv_{i,j}}^n+s\,d\mu
 =s(L_iv_{i,j}+q_i)\omega_i^n
       +s^2e_{i,j}(s)\omega_i^n.
\end{equation}
The remainder coefficient $e_{i,j}(s)=e_{i,j}(s,x)$ extends
smoothly to $s=0$. Indeed, choose $0<s_*\leq1$, depending on $j$,
so that
\[
 \omega_{u_i+s v_{i,j}}\geq\tfrac12\omega_i,
 \qquad |s|\leq s_*,\qquad i=0,1.
\]
The scalar-curvature density is then smooth jointly in $(x,s)$,
and Taylor's integral remainder formula defines the smooth extension
of $e_{i,j}$. On $\{r>2R_j\}$, the perturbation term $v_{i,j}$
and all its derivatives vanish. Hence
$\omega_{u_i+s v_{i,j}}=\omega_i$, whose scalar curvature is zero,
and $L_iv_{i,j}=0$ there. It follows from
\eqref{eq:almost-critical-expansion} that $e_{i,j}(s)=0$ on
$\{r>2R_j\}$, also at $s=0$ by continuity.
For $|s|\leq s_*$, we therefore have
\[
 \begin{aligned}
 \int_Xr^{-\beta}|e_{i,j}(s)|\,\omega_i^n
 &=\int_{\{r\leq2R_j\}}r^{-\beta}|e_{i,j}(s)|\,\omega_i^n\\
 &\leq
 \left(\sup_{\substack{|\sigma|\leq s_*\\r(x)\leq2R_j}}
       |e_{i,j}(\sigma,x)|\right)
 \int_{\{r\leq2R_j\}}r^{-\beta}\,\omega_i^n
 \leq C_j.
 \end{aligned}
\]
Joint continuity on the compact product
$\{r\leq2R_j\}\times[-s_*,s_*]$ bounds the supremum, and the
weighted volume is finite since $r$ is positive.
The constant $C_j$ depends on $n,\beta,(X,\omega),r,R_j$, the
$C^2$ bounds and positive lower bounds of $\omega_0,\omega_1$
on $\{r\leq2R_j\}$, and the $C^4$ norms of $v_{0,j},v_{1,j}$,
with all norms taken relative to $\omega$.
It is independent of $i$ and $s$ for $|s|\leq s_*$;
no bound uniform in $j$ is needed.

Choose $s_j\downarrow0$ within those intervals so that
\begin{equation}\label{eq:diagonal}
 s_jC_j\leq\epsilon_j,\qquad
 s_j\left(\|v_{0,j}\|_{C^k_{-\gamma}} + \|v_{1,j}\|_{C^k_{-\gamma}} \right)\leq j^{-1}.
\end{equation}
The weighted norms are finite because $v_{i,j}$ are compactly
supported. Decrease $s_j$ further, if necessary, so that
$u_i^j=u_i+s_jv_{i,j}$ satisfy
$\omega_{u_i^j}\geq\frac12\kappa\omega$, where
$\omega_i\geq\kappa\omega$ for $i=0,1$.
The perturbed endpoint potentials then have uniformly bounded
weighted $C^k$ norms and a common positive metric lower bound.
By Proposition~\ref{prop:weighted-time-derivative} and the joint
$C^{1,1}$ estimate in \cite[Theorem~A]{YaoGeo}, we have
\[
 |\dot U_t^j|\leq Cr^{-\gamma}\leq Cr^{-\beta},\qquad
 0\leq\omega_{U_t^j}\leq\Lambda\omega,
\]
with $C,\Lambda$ independent of $j$ and the supports of $v_{i,j}$.
Their dependence is through $n,\gamma,k,(X,\omega),r$, the common
endpoint bounds, and $\kappa$. The weak geodesics are the solutions
specified in Section~\ref{sec:preliminary}.

Apply \eqref{eq:almost-critical-expansion}, with $s=s_j$, and integrate against
$\dot U_i^j$. By \eqref{eq:approximate-corrections} and
\eqref{eq:diagonal},
\begin{equation}\label{eq:perturbed-variation-error}
 \begin{aligned}
 &\left|\int_X\dot U_i^j
       \bigl(-S_{\omega_{u_i^j}}\omega_{u_i^j}^n+s_j\,d\mu\bigr)\right|\\
 &\qquad\leq Cs_j\int_Xr^{-\beta}|L_iv_{i,j}+q_i|\,\omega_i^n
       +Cs_j^2\int_Xr^{-\beta}|e_{i,j}(s_j)|\,\omega_i^n\\
 &\qquad\leq C(s_j\epsilon_j+s_j^2C_j)
       \leq Cs_j\epsilon_j.
 \end{aligned}
\end{equation}
This proves \eqref{eq:almost-critical}. The perturbed forms are
scalar-flat outside compact sets, so Lemma~\ref{prop:endpoint-slopes}
applies to $U^j$. By compact support of $\mu$ and joint $C^{1,1}$
regularity, we have
\[
 \frac{d}{dt}\mathcal F_\mu(U_t^j)=\int_X\dot U_t^j\,d\mu
\]
also at the endpoints, with one-sided derivatives there.
Adding $s_j$ times this identity to \eqref{eq:initial-slope} and
\eqref{eq:final-slope}, then using \eqref{eq:almost-critical},
gives the two asserted bounds for the perturbed functional.
\end{proof}

\noindent\textit{Remark.}
On a compact K\"ahler manifold, Berman--Berndtsson
\cite[Section~4]{BB} use the perturbing functional
$\int_Xv\,d\mu-E(v)$, where $E$ is the Aubin--Mabuchi energy and
$\int_Xd\mu=\int_X\omega^n$. Its first variation at $u_i$ is
$\int_X\xi(d\mu-\omega_i^n)$, which vanishes on constant variations.
After choosing representatives under the automorphism group so
that the first variation also annihilates the remaining holomorphy
potentials, they solve
\[
 L_i v_i\,\omega_i^n=\omega_i^n-d\mu.
\]
The perturbation term $v_i$ cancels the term of order $s$. For each fixed $v_i$,
the first variation of the compact perturbed functional at
$u_i+sv_i$, evaluated on $\xi$, is bounded by
$Cs^2\|\xi\|_{C^0}$ for sufficiently small $s$.

In the compact case, the first-variation error is $O(s_j^2)$;
on the ALE manifold, it is
$O(C_js_j^2+s_j\epsilon_j)$ by \eqref{eq:perturbed-variation-error}.
For $v_{i,j}=\chi_{R_j}v_i$, the constant $C_j$ may grow with
$R_j$ and the $C^4$ norms of the perturbation terms.
The additional term $s_j\epsilon_j$ comes from the weighted
approximation of $L_iv_i=-q_i$ by $L_iv_{i,j}$ in
\eqref{eq:approximate-corrections}, whose error tends to zero
as $R_j\to\infty$.

For the uniqueness proof, take the perturbed endpoint potentials
and weak geodesics supplied by Proposition~\ref{prop:perturbed-mabuchi}
for the fixed domain $K$ and density $\mu$ chosen above.

\subsection{A proof of uniqueness}
\begin{proof}[Proof of Theorem~\ref{thm:uniqueness}]
By \cite[Theorem~D(i)]{YaoMass}, we may take $\omega=\omega_0$
and write $\omega_1=\omega+dd^cu_1$ with
$u_1\in\HH_{-\gamma}(X,\omega)$ for some $\gamma>2n-4$; set $u_0=0$.
Set $m_j(t)=\MM(U_t^j)$ and $p_j(t)=\mathcal F_\mu(U_t^j)$.
Both functions are convex. The joint $C^{1,1}$ regularity of $U^j$
and compact support of $\mu$ imply $p_j\in C^1([0,1])$, with
one-sided derivatives at $0$ and $1$.
By Lemma~\ref{prop:endpoint-slopes} and
Proposition~\ref{prop:perturbed-mabuchi}, we therefore have
\begin{align}
 0\leq s_j\bigl(p_j'(1)-p_j'(0)\bigr)
 &\leq(m_j+s_jp_j)'_-(1)-(m_j+s_jp_j)'_+(0)\notag\\
 &\leq\int_X\dot U_1^j
       \bigl(-S_{\omega_{u_1^j}}\omega_{u_1^j}^n+s_j\,d\mu\bigr)
       \notag\\
 &\quad-\int_X\dot U_0^j
       \bigl(-S_{\omega_{u_0^j}}\omega_{u_0^j}^n+s_j\,d\mu\bigr)
 \leq Cs_j\epsilon_j .\label{eq:slope-comparison}
\end{align}
We now use the bound in \eqref{eq:slope-comparison} to control
the gradient of the difference of the endpoint potentials on $K$.
In the local logarithmic coordinate $\tau=t+\sqrt{-1}s$, positivity of
$\Omega_{U^j}=\Omega+\ddc U^j$ means that its
block matrix is positive semidefinite, with spatial block $g_{U_t^j}$, mixed block
$\frac12\partial_X\dot U_t^j$, and time entry
$\frac14\ddot U_t^j$. Add $\epsilon g$, with $\epsilon>0$,
to the spatial block and take the Schur complement.
Since $g_{U_t^j}\leq\Lambda g$, with $\Lambda$ independent of $j$,
letting $\epsilon\downarrow0$, we have, almost everywhere,
\begin{equation}\label{eq:local-strictness}
 \ddot U_t^j\geq\Lambda^{-1}
              |\partial_X\dot U_t^j|_\omega^2 .
\end{equation}
The estimate remains valid when the spatial form is degenerate.

Choose $a_K>0$ with $d\mu\geq a_K\omega^n$ on $K$.
Integrate \eqref{eq:local-strictness} in space and time and use
Cauchy--Schwarz in time:
\begin{align*}
 p_j'(1)-p_j'(0)
 &=\int_0^1\int_X\ddot U_t^j\,d\mu\,dt\\
 &\geq\frac{a_K}{\Lambda}
       \int_0^1\int_K|\partial_X\dot U_t^j|_\omega^2\,\omega^n\,dt\\
 &\geq\frac{a_K}{\Lambda}
       \int_K|\partial_X(u_1^j-u_0^j)|_\omega^2\,\omega^n .
\end{align*}
Joint $C^{1,1}$ regularity on
$\operatorname{supp}\mu\times[0,1]$ justifies Fubini's theorem
and the fundamental theorem of calculus.
Combine the preceding estimate with \eqref{eq:slope-comparison}
and divide by $s_j>0$.
By \eqref{eq:diagonal}, $u_1^j-u_0^j\to u_1-u_0$
locally in $C^1$. Letting $j\to\infty$, we obtain
\[
 \int_K|\partial_X(u_1-u_0)|_\omega^2\,\omega^n=0 .
\]
Since $K$ was arbitrary and $X$ is connected, $u_1-u_0$ is
constant on $X$. Both endpoint potentials tend to zero at infinity,
so the constant is zero. Thus $u_1=u_0$ and $\omega_1=\omega_0$.
\end{proof}

\subsection{Proof of Corollary~\ref{cor:order-zero}}
\label{subsec:extensions}
A scalar-flat K\"ahler metric satisfies
$\Delta_g\Ric(g)=\mathrm{Rm}(g)*\Ric(g)$, where the right-hand side
is a fixed sum of curvature contractions. Hence
\cite[Theorem~1.4]{AcheViaclovsky}, with equation parameter $1$, gives
coordinates on each end in which each $g_i$ has decay $2-2n$ with differentiated
estimates. Only finitely many input derivatives are required; see
\cite[Remark~5.9]{AcheViaclovsky}. This improved decay satisfies the
hypotheses of the one-end results of Hein--LeBrun \cite{HeinLeBrun},
so $X$ has only one end. We must check that both metrics
have this decay in a common chart.

Let $x$ and $y$ be the improved coordinates for $g_0$ and $g_1$,
respectively, lifted to the Euclidean covers of the ends, and put
$F=y\circ x^{-1}$. We regard $g_0$ and $g_1$ as metrics on their
respective coordinate domains. Thus $F^*g_1$ is the second metric
expressed in the $x$ chart: for tangent vectors $v,w$ at $x$,
\[
 (F^*g_1)_x(v,w)=(g_1)_{F(x)}(dF_xv,dF_xw).
\]
In particular, $(F^*g_{\mathrm{euc}})_x(v,w)$ is the Euclidean
inner product of $dF_xv$ and $dF_xw$. Unless otherwise indicated,
tensor norms below are induced by the Euclidean metric on the
corresponding coordinate domain. The operator norm of $dF_x$ uses
the Euclidean metrics on the source $T_x\RR^{2n}$ and the target
$T_{F(x)}\RR^{2n}$; the norm of $(dF_x)^{-1}$ uses the same two
metrics in reverse order.

By \eqref{eq:order-zero-pair},
$|F^*g_1-g_0|_{g_0}=o(1)$. The two metrics on $X$ are therefore
uniformly equivalent. Each ALE radius is comparable to distance
from a fixed compact set in the corresponding metric, so
$c|x|\leq|F(x)|\leq C|x|$ for uniform positive constants $c,C$.
Moreover, for a function $\epsilon(x)\to0$,
\[
 (1-\epsilon)g_0\leq F^*g_1\leq(1+\epsilon)g_0.
\]
Both metrics are uniformly equivalent to the Euclidean metrics in
their respective charts. The comparison above therefore bounds the
Euclidean length of $dF_xv$ above and below by fixed positive
multiples of the Euclidean length of $v$. The upper bound controls
$dF_x$, and the lower bound controls its inverse, uniformly on the
end. To compare the Euclidean metrics in the two charts, write
\[
 F^*g_{\mathrm{euc}}-g_{\mathrm{euc}}
 =(F^*g_1-g_0)+(g_0-g_{\mathrm{euc}})
       -F^*(g_1-g_{\mathrm{euc}}).
\]
The first two terms tend to zero in the Euclidean tensor norm on
the $x$ chart, and the last term satisfies
\[
 |F^*(g_1-g_{\mathrm{euc}})|_{g_{\mathrm{euc}}}(x)
 \leq |dF_x|_{\mathrm{op}}^2
       |g_1-g_{\mathrm{euc}}|_{g_{\mathrm{euc}}}(F(x))=o(1).
\]
In conclusion, we have therefore proved
\begin{equation}\label{eq:common-chart-basic-bounds}
 \begin{gathered}
 c|x|\leq|F(x)|\leq C|x|,\qquad
 |dF_x|_{\mathrm{op}}+|(dF_x)^{-1}|_{\mathrm{op}}\leq C,\\
 |F^*g_{\mathrm{euc}}-g_{\mathrm{euc}}|_{g_{\mathrm{euc}}}=o(1).
 \end{gathered}
\end{equation}
We put $r=|x|$, and fix a sufficiently large finite $k\geq3$
up to which the improved metric estimates hold.

\begin{lemma}\label{lem:transition-laplacian}
The transition map $F$ satisfies $\Delta_{\mathrm{euc}}F=h$, where
\begin{equation}\label{eq:transition-laplacian-decay}
 |D^jh|\leq C_jr^{1-2n-j},\qquad 0\leq j\leq k-2.
\end{equation}
Here $D$ denotes Euclidean coordinate differentiation, and $C_j$
depends only on $n,k$, the metric bounds above, and the bounds in
\eqref{eq:common-chart-basic-bounds}.
\end{lemma}

\begin{proof}
The identity map from $(X,g_0,J)$ to $(X,g_1,J)$ is harmonic.
Indeed, in common holomorphic coordinates the mixed Christoffel
symbols of both K\"ahler metrics vanish, so the $g_0$-trace of the
difference of their connections is zero. Expressing the identity
map in the real ALE coordinates, we have
\begin{equation}\label{eq:transition-harmonic-map}
 (g_0)^{ab}\bigl(\partial_{ab}F^\alpha
 -\Gamma(g_0)^c_{ab}\partial_cF^\alpha
 +\Gamma(g_1)^\alpha_{\mu\nu}(F)
          \partial_aF^\mu\partial_bF^\nu\bigr)=0.
\end{equation}
The ALE coordinates need not be holomorphic for this equation to
hold. The improved metric
estimates and $cr\leq|F(x)|\leq Cr$ imply
\[
 \begin{gathered}
 |g_0^{-1}-g_{\mathrm{euc}}^{-1}|_{g_{\mathrm{euc}}}
       =O(r^{2-2n}),\\
 \Gamma(g_0)(x)=O(r^{1-2n}),\qquad
 \Gamma(g_1)(F(x))=O(r^{1-2n}).
 \end{gathered}
\]

We first prove $D^2F=O(r^{-1})$ by rescaling. Set
$F_R(z)=R^{-1}F(Rz)$ on $\{1/2<|z|<2\}$. Then
\[
 |F_R|+|D_zF_R|\leq C,\qquad
 D_z^2F_R(z)=R D_x^2F(Rz).
\]
Equation \eqref{eq:transition-harmonic-map} becomes
\[
 \begin{aligned}
 (g_0)^{ab}(Rz)\partial_{ab}F_R^\alpha
 ={}&R(g_0)^{ab}(Rz)\Gamma(g_0)^c_{ab}(Rz)
          \partial_cF_R^\alpha\\
 &-R(g_0)^{ab}(Rz)\Gamma(g_1)^\alpha_{\mu\nu}(RF_R(z))
          \partial_aF_R^\mu\partial_bF_R^\nu.
 \end{aligned}
\]
The rescaled equation is uniformly elliptic with uniform
coefficient bounds. For fixed $0<\alpha<1$, the uniform $C^1$
bound for $F_R$ bounds the H\"older norm of the right-hand side by
$C(1+[D_zF_R]_\alpha)$. By interior Schauder estimates and
interpolation, absorbing the first-derivative H\"older term, we obtain
$\|F_R\|_{C^{2,\alpha}(\{3/4<|z|<3/2\})}\leq C$.
Rescaling back yields $|D_x^2F|\leq Cr^{-1}$. Differentiating the
equation and repeating the interior estimates, we obtain
\[
 |D^jF|\leq C_jr^{1-j},\qquad 1\leq j\leq k,
\]
with the corresponding uniform bounds in rescaled H\"older norms.

We now estimate $h=\Delta_{\mathrm{euc}}F$. Replacing the leading
coefficients in \eqref{eq:transition-harmonic-map} by the Euclidean
ones, we have
\[
 \begin{aligned}
 h^\alpha
 ={}&\bigl(\delta^{ab}-(g_0)^{ab}\bigr)\partial_{ab}F^\alpha\\
 &+(g_0)^{ab}\Gamma(g_0)^c_{ab}\partial_cF^\alpha\\
 &-(g_0)^{ab}\Gamma(g_1)^\alpha_{\mu\nu}(F)
          \partial_aF^\mu\partial_bF^\nu.
 \end{aligned}
\]
Consequently,
\[
 |h|\leq Cr^{2-2n}|D^2F|
       +Cr^{1-2n}\bigl(|DF|+|DF|^2\bigr)
 \leq Cr^{1-2n}.
\]
Differentiating the expression for $h$ gives
\eqref{eq:transition-laplacian-decay}. For example, the terms
containing differentiated leading coefficients satisfy
\[
 |D^a(g_0^{-1}-g_{\mathrm{euc}}^{-1})|\,|D^{j-a+2}F|
 \leq Cr^{2-2n-a}r^{-1-j+a}
 =Cr^{1-2n-j},\qquad 0\leq a\leq j.
\]
The differentiated Christoffel terms have the same order, and the
bounds also hold in the corresponding rescaled H\"older norms.
\end{proof}

\begin{proof}[Proof of Corollary~\ref{cor:order-zero}]
Extend $F$ smoothly across the bounded missing region to
$\widetilde F$ on $\RR^{2n}$. By
Lemma~\ref{lem:transition-laplacian} and the weighted Laplace
solvability theorem \cite[Theorem~2.11]{ConlonHein}, applied
componentwise at the weight $3-2n\in(2-2n,0)$, there is a solution
\[
 \Delta_{\mathrm{euc}}w=\Delta_{\mathrm{euc}}\widetilde F,
 \qquad D^jw=O(r^{3-2n-j}),\quad 0\leq j\leq k.
\]
The function $\widetilde F-w$ is harmonic on all of $\RR^{2n}$ and
has at most linear growth, so each component is affine:
\begin{equation}\label{eq:common-improved-chart}
 F(x)=Ax+b+O(r^{3-2n}),\qquad
 DF=A+O(r^{2-2n}).
\end{equation}
Here $A$ is a constant linear map. Taking the limit in
\eqref{eq:common-chart-basic-bounds}, we have
$A^*g_{\mathrm{euc}}=g_{\mathrm{euc}}$, so $A$ is a Euclidean
isometry. By the pullback estimate above and
\eqref{eq:common-improved-chart}, we now have, in the $x$ chart,
\[
 \begin{aligned}
 |F^*g_1-g_{\mathrm{euc}}|_{g_{\mathrm{euc}}}
 &\leq |F^*(g_1-g_{\mathrm{euc}})|_{g_{\mathrm{euc}}}
       +|F^*g_{\mathrm{euc}}-g_{\mathrm{euc}}|_{g_{\mathrm{euc}}}\\
 &=O(r^{2-2n}),
 \end{aligned}
\]
with differentiated bounds through order $k-1$. Both metrics
therefore approach the same Euclidean metric in this chart.
Since $2n-2>n-1$ and the initial $k$ can be taken sufficiently
large, Theorem~\ref{thm:uniqueness} applies.
\end{proof}

\subsection{Proof of Corollary~\ref{cor:ac}}
\label{subsec:ac-extension}
The AC argument uses the same ingredients and follows the same
main steps as the ALE proof of uniqueness. We use the preceding
proofs wherever the arguments are unchanged and explain the
modifications needed on the cone end.

We work in the setting of Corollary~\ref{cor:ac}. Let $(L,g_L)$
be a smooth compact Riemannian manifold of dimension $2n-1$.
Its metric cone is
\[
 C(L)=(0,\infty)\times L,\qquad g_C=dr^2+r^2g_L.
\]
Assume that $g_C$ is a Ricci-flat K\"ahler metric.
A complete K\"ahler manifold $(X,J,\omega)$, with associated metric
$g$, is asymptotic to this cone with rate $-\delta$ if there are
a compact set $K\Subset X$ and a diffeomorphism
$\Phi:(R_0,\infty)\times L\longrightarrow X\setminus K$ such that
\begin{equation}\label{eq:ac-metric-decay}
 |\nabla_C^j(\Phi^*g-g_C)|_{g_C}
       =O(r^{-\delta-j}),\qquad 0\leq j\leq k.
\end{equation}
Here $\nabla_C$ is the Levi-Civita connection of $g_C$, and $k$
is sufficiently large. We fix the cone and $\Phi$, assume
$\delta>n-1$, and extend $r$ smoothly to $X$ with $r\geq1$.
For $\gamma>0$, the space $\HH_{-\gamma}(X,\omega)$ consists of
smooth real functions $u$ satisfying
\[
 \omega_u=\omega+\ddc u>0,\qquad
 |\nabla^j u|_g=O(r^{-\gamma-j}),\qquad 0\leq j\leq k,
\]
where $\nabla$ is the connection of $g$. For the Mabuchi functional
we retain $\gamma>2n-4$ and $S_\omega\in L^1(X,\omega^n)$,
as in Corollary~\ref{cor:ac}.

Let $A=\{\tau\in\CC:1<|\tau|<e\}$, let
$p:X\times A\longrightarrow X$ be the projection, and set
$\Omega=p^*\omega$. As in the ALE case, a weak geodesic joining
$u_0,u_1\in\HH_{-\gamma}(X,\omega)$ is a bounded,
$S^1$-invariant function $U(x,e^{t+\sqrt{-1}s})=u_t(x)$ solving
\begin{equation}\label{eq:ac-hcma}
 \begin{cases}
 (\Omega+\ddc U)^{n+1}=0, \qquad
       &\text{on }X\times A,\\
  \Omega+\ddc U\geq0, \qquad
       &\text{on }X\times A,\\
 U(x,\tau)=u_j(x)&\text{when }|\tau|=e^j,\quad j=0,1,
 \end{cases}
\end{equation}
in the sense of Bedford--Taylor. We use the same Perron class
$\mathcal B(u_0,u_1)$ as in \eqref{eq:perron}: its functions are
bounded, $S^1$-invariant, and $\Omega$-psh, with joint boundary
upper limits at most $u_j$ on $|\tau|=e^j$.

\begin{subequations}\label{eq:ac-geometric-bounds}
The AC assumption supplies the uniform local geometry needed for
the geodesic equation. After dilating a cone annulus to unit scale,
the rescaled metrics converge to $g_C$ on a fixed annulus. Since
$L$ is compact, the metric $g$ has a positive injectivity radius
and bounded curvature derivatives of the required finite orders;
the same bounds hold on the compact part of $X$. We therefore
have holomorphic coordinate balls with uniform metric bounds,
as in \cite[Theorem~9]{WuYau20}. More precisely, the curvature of
the cone has order $r^{-2}$, while its Ricci curvature vanishes,
so by \eqref{eq:ac-metric-decay}, we have
\begin{equation}
 |\nabla^j\mathrm{Rm}(g)|\leq C_jr^{-2-j},\qquad
 |\nabla^j\Ric(\omega)|\leq C_jr^{-\delta-2-j},
 \qquad 0\leq j\leq k-2.
\end{equation}
Using the cone volume form and radial coordinate, we also have
\begin{equation}
 \operatorname{Vol}_\omega\{R<r<2R\}\leq CR^{2n},\qquad
 |dr|\leq C,\qquad |\nabla^2r|\leq Cr^{-1}
\end{equation}
on the end. In particular, the extended radius is a proper
function with bounded gradient and Hessian. The constants depend
on $n,(X,\omega),r$ and the finite AC bounds, but not on $R$.
\end{subequations}

\begin{lemma}\label{lem:ac-geodesic}
For $u_0,u_1\in\HH_{-\gamma}(X,\omega)$, the weak geodesic $U$
is unique and continuous on $X\times\overline A$ and has joint
$C^{1,1}$ regularity. For $\widehat g=g+dt^2+ds^2$, where
$\tau=e^{t+\sqrt{-1}s}$,
\begin{equation}\label{eq:ac-joint-bound}
 \|U\|_{L^\infty}
 +\|\widehat\nabla U\|_{L^\infty}
 +\|\widehat\nabla^2U\|_{L^\infty}\leq C.
\end{equation}
The norms are on $X\times\overline A$, and the Hessian is understood
almost everywhere. The constant depends on $n,(X,\omega)$,
a common positive lower bound $\kappa\omega$ for
$\omega_{u_0},\omega_{u_1}$, and a bound for
$\|u_0\|_{C^3(X,g)}+\|u_1\|_{C^3(X,g)}$.
\end{lemma}

\begin{proof}
We follow the $\epsilon$-geodesic construction of Chen
\cite{Chen00} and its noncompact version in
\cite[Sections~2--5]{YaoGeo}. The Wu--Yau generalized
maximum principle used in \cite[Lemma~2.2 and Theorem~3.2]{YaoGeo}
therefore applies by the same argument, subtracting a small multiple
of $r$. Consequently, the linear Dirichlet problems and continuity
method for the positive $\epsilon$-geodesic equations in
\cite[Section~2]{YaoGeo} apply on the AC manifold.

To make the dependence of the estimates explicit, use the fixed
product form $\Theta=\Omega+dd^c(t(t-1))$ and solve
\[
 (\Omega+dd^cU_\epsilon)^{n+1}
 =\epsilon\Theta^{n+1},\qquad 0<\epsilon\leq1,
\]
with the prescribed boundary values. Relative to $\Theta$, the
potential $U_\epsilon-t(t-1)$ has constant logarithmic
Monge--Amp\`ere density $\log\epsilon$, whose derivatives vanish.
The curvature bounds of $\Theta$ are independent of the endpoint
potentials. The interpolation $(1-t)u_0+tu_1+C t(t-1)$, with
$C$ sufficiently large, supplies a subsolution and a starting point
for the continuity method; the boundary estimates require endpoint
derivatives only through order three.

The boundary and Laplacian estimates of Chen \cite{Chen00},
the gradient estimate of B{\l}ocki \cite{Blocki09,Blocki12}, and
the joint Hessian estimate of Chu--Tosatti--Weinkove \cite{CTW}
apply as in \cite[Sections~3--5]{YaoGeo}. Their local calculations
are unchanged; the preceding maximum principle justifies their
use on the noncompact manifold. For the $\epsilon$-geodesics
$U_\epsilon$, we therefore have
\[
 \sup_{0<\epsilon\leq1}
 \left(\|U_\epsilon\|_{L^\infty}
 +\|\widehat\nabla U_\epsilon\|_{L^\infty}
 +\|\widehat\nabla^2U_\epsilon\|_{L^\infty}\right)\leq C,
\]
with the same dependence as in \eqref{eq:ac-joint-bound}.
A subsequence converges locally in $C^{1,\alpha}$, $0<\alpha<1$,
and Bedford--Taylor continuity gives a solution of
\eqref{eq:ac-hcma} satisfying \eqref{eq:ac-joint-bound}.

The solution can also be constructed by the Perron envelope,
following \cite[Section~3 and Theorem~5.7]{YaoAsymptotic}.
The envelope and comparison arguments there use uniform local
coordinate bounds, boundary barriers on the annulus, and the
proper-radius maximum principle, all available above. The endpoint
potentials are bounded and uniformly continuous, and the annulus
has strongly regular boundary. By the same arguments, we therefore have
\[
 U=\sup_{V\in\mathcal B(u_0,u_1)}V
\]
as the unique bounded continuous solution of \eqref{eq:ac-hcma}.
In particular, the Perron solution agrees with the
$\epsilon$-geodesic limit.
\end{proof}

The relative Mabuchi functional has the same definition as in
the ALE case. With $\vartheta=\Ric(\omega)$ and
$h(s)=s\log s-s+1$, where $h(0)=1$, we set
\[
 \begin{aligned}
 \MM(u)&=\Ent_{\mathrm{rel}}(u)-nE_\vartheta(u)\\
 &=\int_X h\!\left(\frac{\omega_u^n}{\omega^n}\right)\omega^n
   -\sum_{j=0}^{n-1}\int_X
       u\,\vartheta\wedge\omega_u^{n-1-j}\wedge\omega^j.
 \end{aligned}
\]
For $u\in\HH_{-\gamma}(X,\omega)$, the same decay estimates as
in \eqref{eq:smooth-entropy-tail} imply that the entropy is finite. Finiteness
of the functional along weak geodesics requires the approximation
argument below.

We use the same compactly supported approximation as in
Sections~\ref{sec:ale-finiteness} and~\ref{sec:ale}.
The radial-barrier proof of Lemma~\ref{lem:support} applies by
\eqref{eq:ac-geometric-bounds}, so each geodesic with compactly supported endpoint
potentials has compact spatial support independent of time.
By Propositions~\ref{prop:compact-convexity}
and~\ref{prop:compact-continuity}, the Mabuchi functional is then
continuous and convex along those geodesics.

For general endpoint potentials, the same radial barriers establish
the estimates of Propositions~\ref{prop:weighted-geodesic-estimates}
and~\ref{prop:weighted-time-derivative}. Rescaling cone annuli
provides the uniform local metric and curvature bounds used in
the spatial-flow argument of Proposition~\ref{prop:weighted-gradient},
so the weighted spatial gradient estimate also holds.

In particular, the AC weak geodesics satisfy
\[
 |u_t|+r|\nabla_Xu_t|+|\dot U|\leq Cr^{-\gamma},\qquad
 0\leq t\leq1,
\]
and Corollary~\ref{cor:integrated-time-derivatives} gives
\[
 \int_0^1\bigl(\ddot U+|\ddot U|^2+
                  |\nabla_X\dot U|^2\bigr)\,dt
       \leq Cr^{-2\gamma-2}
 \quad\text{for almost every }x\in X.
\]
The constants have the dependencies stated in the cited results,
with the AC bounds in \eqref{eq:ac-geometric-bounds}; they are
uniform for the cutoff endpoint potentials and the small
perturbations with common endpoint bounds used below.

Cutting off the endpoint potentials, we obtain compactly supported
weak geodesics $u_{t,l}$ converging uniformly to $u_t$ by
\eqref{eq:perron-stability}. Following the proofs of
Lemmas~\ref{lem:endpoint-energy-convergence}
and~\ref{lem:twisted-energy-convergence}, we obtain convergence of the
endpoint Mabuchi energies and uniform convergence in $t$ of
$E_\vartheta(u_{t,l})$. Convexity therefore bounds the relative
entropies uniformly. Lower semicontinuity in
Lemma~\ref{lem:entropy-lsc}, followed by the entropy--Hessian
estimate in Lemma~\ref{lem:entropy-hessian}, implies finiteness of
$\MM(u_t)$ exactly as in the proof of Theorem~\ref{thm:main-a}.

The rooftop and contact-set argument of
Lemma~\ref{lem:compact-density-convergence} proves entropy convergence
on compact subsets without change. The volume and Ricci bounds
in \eqref{eq:ac-geometric-bounds}, together with the integrated
time estimates above, also imply the estimates of
Lemma~\ref{lem:quantitative-entropy-tails} with the same decay rates.
As in Proposition~\ref{prop:interior-entropy-convergence}, the
difference between a subsequential limit of $\MM(u_{t,l})$ and
$\MM(u_t)$ is nonnegative and convex, with zero endpoint limits,
and hence vanishes. Thus the Mabuchi energies converge at every
interior time. Convexity passes to the limit, and the endpoint
argument of Theorem~\ref{thm:main-b} proves continuity on $[0,1]$.

For uniqueness, we need the following solvability result for
compactly supported sources. The Fredholm argument is the same
as in Lemma~\ref{prop:weighted-density}; we verify the vanishing
of decaying Killing fields and the required weights on the cone.

\begin{lemma}\label{lem:ac-lichnerowicz}
Under the geometric assumptions of Corollary~\ref{cor:ac}, assume
$n\geq3$ and $S_\omega=0$. For every
$q\in C_c^\infty(X,\RR)$, there is a smooth real function $v$
satisfying
\begin{equation}\label{eq:ac-compact-source-solution}
 L_\omega v=-q,\qquad
 |\nabla^jv|\leq C_qr^{4-2n-j},\qquad 0\leq j\leq4.
\end{equation}
The constant $C_q$ depends on $n,(X,\omega),r$, the finite AC
bounds, and $q$, and is independent of any later cutoff radius.
\end{lemma}

\begin{proof}
We first show that a Killing field tending to zero on the AC end
vanishes. By the metric decay, we have
\begin{equation}\label{eq:ac-radius-hessian}
 |\nabla r|_g^2=1+O(r^{-\delta}),\qquad
 \nabla_g^2r=r^{-1}(g-dr\otimes dr)+O(r^{-1-\delta}).
\end{equation}
Choose a unit-speed geodesic $\zeta$ initially directed along
$\nabla r$ from a sufficiently distant point, and set
$a=(r\circ\zeta)'$. Then
\[
 a'=\frac{1-a^2}{r}+O(r^{-1-\delta}),\qquad
 a^2\leq1+Cr^{-\delta}.
\]
While $a\geq1/2$, we have, by integration,
$a(t)\geq a(0)-Cr(\zeta(0))^{-\delta}>1/2$ if the initial
point is sufficiently far out. Thus $\zeta$ escapes and
$c(1+t)\leq r(\zeta(t))\leq C(1+t)$. Moreover,
$(r^2(1-a^2))'=O(r^{1-\delta})$, so integrating, we obtain
\[
 |\dot\zeta-a\partial_r|_{g_C}^2
 =1-a^2+O(r^{-\delta})\leq Cr^{-\epsilon},
 \qquad \epsilon=\min\{\delta,1\}>0.
\]
Since the cone curvature vanishes in radial directions,
\[
 |\mathrm{Rm}_g(\dot\zeta,\cdot)\dot\zeta|
 \leq Cr^{-2}|\dot\zeta-a\partial_r|_{g_C}^2+Cr^{-2-\delta}
 \leq C(1+t)^{-2-\epsilon}.
\]
The constants may depend on the initial point and the AC bounds,
but not on $t$. The Jacobi-field argument of
Lemma~\ref{lem:killing} now makes every decaying Killing field
vanish along $\zeta$. Since the initial point was arbitrary
sufficiently far out, the field vanishes on the end, and hence
on connected $X$.

If a real distributional solution of $L_\omega h=0$ satisfies
$h=O(r^{-\beta})$ for $\beta>n-2$,
we may apply the rescaled estimates and integration by parts in
Lemma~\ref{lem:decaying-kernel} on cone annuli to obtain
\[
 0\leq\int_X\chi_R|\mathcal D_\omega h|^2\,\omega^n
       \leq CR^{2n-2\beta-4}\longrightarrow0.
\]
Thus $\mathcal D_\omega h=0$, so $J\nabla h$ is a decaying
Killing field. The preceding vanishing and the decay of $h$
imply $h=0$.

Since the cone is Ricci-flat, the model operator is a positive
multiple of $\Delta_C^2$. For $-\Delta_LY=\mu Y$, $\mu\geq0$,
\begin{equation}\label{eq:ac-indicial}
 \begin{aligned}
 \Delta_C^2(r^\lambda Y)
 ={}&[\lambda(\lambda+2n-2)-\mu]\\
 &\cdot[(\lambda-2)(\lambda+2n-4)-\mu]
                 r^{\lambda-4}Y.
 \end{aligned}
\end{equation}
The roots are real, and both factors are strictly negative for
$4-2n<\lambda<0$. Hence there are no exceptional weights in
that interval.

Fix $\sigma=\tfrac12$ and set $\lambda=2-n+\sigma$.
The rescaled AC coefficient bounds allow the weighted Fredholm
theorem of \cite[Theorem~6.2 and Section~7]{LockhartMcOwen} to
be applied as in Lemma~\ref{prop:weighted-density}. Formal
self-adjointness places the obstruction at complementary weight
$4-2n-\lambda=2-n-\sigma$; by rescaled elliptic estimates, we have
$h=O(r^{2-n-\sigma+\eta})$ for every sufficiently small $\eta>0$.
Choosing $\eta<\sigma$, the vanishing above eliminates the
obstruction. By the Fredholm alternative and elliptic regularity,
there is a smooth real solution with
$|\nabla^jv|=O(r^{2-n+\sigma-j})$, $0\leq j\leq4$.

Outside $\supp q$, comparing with the cone operator, we have
\[
 c\Delta_C^2v=(c\Delta_C^2-L_\omega)v
       =O(r^{-2n-\epsilon}),\qquad
 \epsilon=\delta-(n-2+\sigma)>0.
\]
The model expansion argument in the proof of
Lemma~\ref{prop:weighted-density} applies with the eigenfunctions
of the cone link. Formula~\eqref{eq:ac-indicial} shows that the
largest negative root is $4-2n$, a simple root corresponding to
the constant link eigenfunction. All other negative roots are
strictly smaller. Choose $0<\eta<\epsilon$ so that
$[4-2n-\eta,4-2n)$ contains no indicial root. By the expansion, we have
\[
 v=c_0r^{4-2n}+O(r^{4-2n-\eta}),
\]
with differentiated remainder estimates through order four.
This proves \eqref{eq:ac-compact-source-solution}.
\end{proof}

\begin{proof}[Proof of Corollary~\ref{cor:ac}]
When $n=2$, the three-dimensional Einstein link of the Ricci-flat
cone has constant sectional curvature one. The cone is a flat
quotient, so the assertion is already the ALE result. Assume
$n\geq3$ below.

It remains to obtain the compactly supported perturbation terms
for uniqueness. For a scalar-flat AC reference form and a real
$q\in C_c^\infty(X)$, take the solution $v$ from
Lemma~\ref{lem:ac-lichnerowicz}. If $\chi_R=1$ on $\supp q$, then
\[
 L_\omega(\chi_Rv)+q
       =L_\omega(\chi_Rv)-\chi_RL_\omega v.
\]
The product estimate in the proof of
Lemma~\ref{prop:weighted-density} uses only the differentiated
bounds for $v$, the cutoff, and the rescaled operator coefficients.
Those bounds hold on cone annuli. Using their volume growth, we have
\begin{equation}\label{eq:ac-compact-source-cutoff}
 \int_Xr^{-\beta}|L_\omega(\chi_Rv)+q|\,\omega^n
       \leq C_qR^{-\beta}.
\end{equation}
Here $C_q$ also depends on $\beta$ and the cutoff profile, but not
on $R$. The right-hand side tends to zero for every $\beta>0$.
Taking $\beta=\gamma/2>n-2$, we obtain the approximation required in
Proposition~\ref{prop:perturbed-mabuchi}.

Finally, by \cite[Theorem~D(i)]{YaoMass}, which is stated for AC
K\"ahler manifolds, we have $\omega_1=\omega_0+\ddc u_1$ with
$u_1\in\HH_{-\gamma}(X,\omega_0)$ for some $\gamma>2n-4$.
Take $\omega=\omega_0$ and $u_0=0$. The proofs of the
change-of-reference identity and endpoint slope inequalities in
Lemmas~\ref{prop:cocycle} and~\ref{prop:endpoint-slopes} use the
integrable bounds verified above and therefore apply without change.

Fix $K\Subset X$ and a smooth nonnegative compactly supported
density $d\mu$ that is strictly positive near $K$. As in the ALE
argument, set
\[
 \mathcal F_\mu(v)=\int_Xv\,d\mu,\qquad
 \MM_{s,\mu}(v)=\MM(v)+s\mathcal F_\mu(v).
\]
We also perturb the endpoint potentials to make the first
variation of $\MM_{s,\mu}$ small compared with $s$.
For $i=0,1$, set $L_i=L_{\omega_i}$ and
$q_i=d\mu/\omega_i^n$. By Lemma~\ref{lem:ac-lichnerowicz}, we can solve
$L_iv_i=-q_i$, and \eqref{eq:ac-compact-source-cutoff} allows us
to choose $v_{i,j}=\chi_{R_j}v_i$, with $R_j\to\infty$, so that
\[
 \int_Xr^{-\beta}|L_iv_{i,j}+q_i|\,\omega_i^n\leq\epsilon_j,
 \qquad \beta=\gamma/2,\qquad \epsilon_j\downarrow0.
\]
Set $u_i^j=u_i+s_jv_{i,j}$ and choose $s_j\downarrow0$ as in
\eqref{eq:diagonal}, preserving a common positive lower bound
for the perturbed forms and ensuring $u_i^j\to u_i$ in
$C^k_{-\gamma}$. The scalar-curvature expansion and weighted
time-derivative estimate used in the proof of
Proposition~\ref{prop:perturbed-mabuchi} apply on the AC manifold.
Choosing $s_j$ small enough also controls the quadratic remainder,
so the perturbed first-variation bound \eqref{eq:almost-critical}
holds for the weak geodesic $U^j$ joining $u_0^j$ and $u_1^j$.

The remaining argument is the proof of
Theorem~\ref{thm:uniqueness}. Convexity and the endpoint slope
inequalities imply \eqref{eq:slope-comparison}. Since $d\mu$ is
positive near $K$, the positivity estimate
\eqref{eq:local-strictness} implies
\[
 s_jc_K\int_K|\partial_X(u_1^j-u_0^j)|_\omega^2\,\omega^n
       \leq Cs_j\epsilon_j,
\]
where $c_K>0$ and $C$ are independent of $j$; their dependence
is through $K,\mu$, the reference manifold $(X,\omega)$, and the
common endpoint and weak-geodesic bounds.
Divide by $s_j$ and let $j\to\infty$. Local $C^1$ convergence of
the endpoint potentials implies $\partial_X(u_1-u_0)=0$ on $K$.
Since $K$ is arbitrary, connectedness and decay at infinity imply
$u_1=u_0$, hence $\omega_1=\omega_0$.
\end{proof}

For completeness, the obstruction mentioned in
Remark~\ref{rem:ac-decay} can be seen directly from periods.
The cone K\"ahler form $\omega_C$ is exact. If the restriction of
$[\omega]$ to $L$ is nonzero, choose a two-cycle $\Sigma\subset L$
with nonzero period, and denote its copy in $\{r=R\}$ by
$\Sigma_R$. Since $\operatorname{Area}_{g_C}(\Sigma_R)=O(R^2)$,
\[
 0\ne\int_{\Sigma_R}\Phi^*\omega
    =\int_{\Sigma_R}(\Phi^*\omega-\omega_C)
\]
is independent of $R$ and cannot tend to zero. Thus
$|\Phi^*\omega-\omega_C|_{g_C}=o(R^{-2})$ is impossible in that
class. For example, the K\"ahler classes on the canonical bundle of
$\PP^1\times\PP^1$ whose restrictions to the zero section are not
proportional to $c_1(\PP^1\times\PP^1)$ in
\cite[Theorem~1.1 and Section~4.2]{VanCoevering} have this property.
They provide Ricci-flat examples with genuine order $r^{-2}$
terms; the broader scalar-flat line-bundle setting is discussed in
\cite[Section~4.2.2]{YaoInvariant}.

\appendix
\section{The Lichnerowicz equation in weighted spaces}
\label{app:weighted-lichnerowicz}

Let $(X,J,\omega)$ be a connected complete K\"ahler
manifold of complex dimension $n\geq2$, with one end asymptotic
to a Ricci-flat K\"ahler cone
\[
 C(L)=(0,\infty)\times L,\qquad g_C=dr^2+r^2g_L,
\]
where $L$ is smooth, compact, and connected. Fix an identification
$\Phi:(R_0,\infty)\times L\longrightarrow X\setminus K$ and a
smooth extension $r\geq1$ to $X$. For the full weighted theory,
assume differentiated decay of every order:
\begin{equation}\label{eq:all-order-ac}
 |\nabla_C^j(\Phi^*g_\omega-g_C)|_{g_C}
       \leq C_jr^{-\delta-j},\qquad j=0,1,\ldots,\qquad\delta>0.
\end{equation}
In this appendix, $L_\omega$ denotes the real self-adjoint
Lichnerowicz operator
\[
 \begin{aligned}
 L_\omega v
 &:=\operatorname{Re}(\mathcal D_\omega^*\mathcal D_\omega v)\\
 &=\Delta_\omega^2v+\langle\Ric(\omega),dd^cv\rangle_\omega
   +\operatorname{Re}\!\left(
       \omega^{a\bar b}(\partial_aS_\omega)\partial_{\bar b}v
     \right),\qquad v\in C^\infty(X,\RR),
 \end{aligned}
\]
where $\mathcal D_\omega=\bar\partial\nabla^{1,0}_\omega$ and
$\mathcal D_\omega^*$ is its Hermitian formal adjoint. Equivalently,
$L_\omega$ is obtained by taking the adjoint with respect to the
real $L^2$ pairing. For $v\in C_c^\infty(X,\RR)$,
\[
 \int_XvL_\omega v\,\omega^n
   =\int_X|\mathcal D_\omega v|^2\,\omega^n.
\]
When $S_\omega=0$, this operator agrees with
\eqref{eq:lich-formula}. Since the cone is Ricci-flat, its model
operator is $c\Delta_C^2$ for a positive constant $c$. The ALE
case corresponds to the round quotient $L=S^{2n-1}/\Gamma$.

For an integer $m\geq0$, $0<\alpha<1$, and an open set
$D\subset\RR^{2n}$, write
\[
 \|v\|_{C^{m,\alpha}(D)}
 =\sum_{|\beta|\leq m}\sup_D|D^\beta v|
 +\sum_{|\beta|=m}\sup_{\substack{x,y\in D\\x\ne y}}
       \frac{|D^\beta v(x)-D^\beta v(y)|}{|x-y|^\alpha}.
\]
Here $\beta$ is a multi-index of real coordinate derivatives.
Fix a compact set $K_*$ containing the compact part and
$\{R_0<r\leq2R_0\}$. Choose fixed finite coordinate covers of $K_*$ and
$[1,2]\times L$ to define the corresponding H\"older norms. Set
$\widetilde v=v\circ\Phi$ and
\[
 \|v\|_{C^{m,\alpha}_\lambda}
 =\|v\|_{C^{m,\alpha}(K_*)}
 +\sup_{R\geq R_0}R^{-\lambda}
   \|\widetilde v(R\,\cdot,\cdot)\|_{C^{m,\alpha}((1,2)\times L)}.
\]
The space $C^{m,\alpha}_\lambda(X)$ consists of real-valued,
locally $C^{m,\alpha}$ functions with finite norm above; in particular,
$|\nabla^jv|=O(r^{\lambda-j})$ for $0\leq j\leq m$.

The exceptional weights are determined by the eigenvalues of
$-\Delta_L$. For $-\Delta_LY=\mu Y$,
\begin{equation}\label{eq:indicial}
 \begin{aligned}
 \Delta_C^2(r^\lambda Y)
 ={}&[\lambda(\lambda+2n-2)-\mu]\\
 &\cdot[(\lambda-2)(\lambda+2n-4)-\mu]
                 r^{\lambda-4}Y.
 \end{aligned}
\end{equation}
Thus the exceptional set is
\[
 \mathcal I_L=
 \bigcup_{\mu\in\operatorname{Spec}(-\Delta_L)}
 \left\{-(n-1)\pm\sqrt{(n-1)^2+\mu},\,
              3-n\pm\sqrt{(n-1)^2+\mu}\right\}.
\]
It is invariant under $\rho\mapsto4-2n-\rho$ and has no element
in $(4-2n,0)$ when $n>2$. Write
$\rho_+=\min(\mathcal I_L\cap(0,\infty))$.
For each $\rho\in\mathcal I_L$, let $m_L(\rho)$ be the dimension
of the space of biharmonic solutions on $C(L)$ with radial exponent $\rho$,
including logarithmic solutions associated with repeated roots
in \eqref{eq:indicial}. In particular, $m_L(0)=1$ for $n>2$;
for $n=2$, the two radial solutions are $1$ and $\log r$, so
$m_L(0)=2$.

For the ALE link $S^{2n-1}/\Gamma$, the eigenvalues are
$d(d+2n-2)$ for degrees $d$ admitting nonzero $\Gamma$-invariant
spherical harmonics. With $\HHH_d^\Gamma$ denoting their space,
\[
 \mathcal I_\Gamma
 =\bigcup_{d:\,\HHH_d^\Gamma\ne0}
       \{d,d+2,2-2n-d,4-2n-d\}.
\]
The exceptional weights are therefore integers in the ALE case;
for a general cone, they depend on the spectrum of its link.

\begin{proposition}[The Lichnerowicz equation in weighted spaces]
\label{prop:fredholm}
Let $(X,J,\omega)$ be a connected complete AC
K\"ahler manifold with one end asymptotic to a Ricci-flat
K\"ahler cone as above, satisfying \eqref{eq:all-order-ac}.
Fix an integer $m\geq0$ and $0<\alpha<1$, and put
\[
 \mathcal K_\nu
 =\{h\in C^{4,\alpha}_\nu(X):L_\omega h=0\}.
\]
For $\lambda\notin\mathcal I_L$, the operator
\[
 L_{\omega,\lambda}:C^{m+4,\alpha}_\lambda(X)
              \longrightarrow C^{m,\alpha}_{\lambda-4}(X)
\]
is Fredholm, with kernel $\mathcal K_\lambda$ and range
\begin{equation}\label{eq:weighted-obstruction}
 \left\{q\in C^{m,\alpha}_{\lambda-4}(X):
       \int_Xqh\,\omega^n=0
       \text{ for all }h\in\mathcal K_{4-2n-\lambda}\right\}.
\end{equation}
The pairings converge absolutely. Solutions are unique modulo
$\mathcal K_\lambda$, and
\begin{equation}\label{eq:weighted-quotient-estimate}
 \inf_{h\in\mathcal K_\lambda}
       \|v-h\|_{C^{m+4,\alpha}_\lambda}
 \leq C\|L_\omega v\|_{C^{m,\alpha}_{\lambda-4}},
\end{equation}
where $C$ depends on $n,m,\alpha,\lambda$, $(X,\omega)$, and the
fixed weighted norms, but not on $v$. For nonexceptional weights
$\lambda_1<\lambda_2$,
\begin{equation}\label{eq:lich-index-jump}
 \operatorname{ind}L_{\omega,\lambda_2}
 -\operatorname{ind}L_{\omega,\lambda_1}
 =\sum_{\rho\in\mathcal I_L\cap(\lambda_1,\lambda_2)}m_L(\rho).
\end{equation}
Moreover,
\begin{equation}\label{eq:all-negative-kernel}
 \mathcal K_\nu=\{0\}\quad(\nu<0),\qquad
 \mathcal K_\nu=\RR\quad(0\leq\nu<\rho_+).
\end{equation}
In particular, $L_{\omega,\lambda}$ is surjective at every
nonexceptional weight $\lambda>4-2n$, and the following maps
are isomorphisms:
\begin{align}
 &L_\omega:C^{m+4,\alpha}_\lambda
   \longrightarrow C^{m,\alpha}_{\lambda-4},
 &4-2n<\lambda<0,\quad n>2;
       \label{eq:inverse-negative}\\
 &L_\omega:C^{m+4,\alpha}_\lambda/\RR
   \longrightarrow C^{m,\alpha}_{\lambda-4},
 &0<\lambda<\rho_+;
       \label{eq:inverse-positive}\\
 &L_\omega:C^{m+4,\alpha}_\lambda
   \longrightarrow
   \left\{q\in C^{m,\alpha}_{\lambda-4}:
                   \int_Xq\,\omega^n=0\right\},
 &4-2n-\rho_+<\lambda<4-2n.
       \label{eq:inverse-charge}
\end{align}
\end{proposition}

\begin{proof}
The metric estimates \eqref{eq:all-order-ac} and Ricci-flatness
of the cone imply $|\nabla^jS_\omega|=O(r^{-\delta-2-j})$.
The additional first-order term in $L_\omega$ therefore has
coefficient $O(r^{-\delta-3})$, with the corresponding differentiated
bounds. Thus the coefficient comparison in the proof of
Lemma~\ref{lem:ac-lichnerowicz}, now with bounds of every order,
still identifies $c\Delta_C^2$ as the limiting operator.
The weighted Fredholm theorem of
\cite[Theorem~6.2]{LockhartMcOwen} and its H\"older formulation
\cite[Theorem~6.10]{Marshall02} give Fredholmness and the range
formula, using formal self-adjointness and the indicial polynomial
\eqref{eq:indicial}. By the bounded inverse theorem on the quotient,
we have \eqref{eq:weighted-quotient-estimate}. The index formula
\eqref{eq:lich-index-jump} follows from
\cite[Theorem~6.2]{LockhartMcOwen}; the H\"older and Sobolev
kernels and indices agree by
\cite[Corollary~4.27 and Equation~(6.14)]{Marshall02}.

Kernels are unchanged when the weight varies in an interval
without exceptional weights \cite[Lemma~7.1]{LockhartMcOwen}.
For $\nu=4-2n-\lambda\notin\mathcal I_L$, choose $\epsilon>0$
so that $\mathcal K_\nu=\mathcal K_{\nu-\epsilon}$.
For $q\in C^{m,\alpha}_{\lambda-4}$ and $h\in\mathcal K_\nu$,
the integral of $|qh|$ on the end is then bounded by a constant
times
\[
 \int_{R_0}^{\infty}
       r^{\lambda-4+\nu-\epsilon+2n-1}\,dr
       =\int_{R_0}^{\infty}r^{-1-\epsilon}\,dr<\infty.
\]
This proves absolute convergence in
\eqref{eq:weighted-obstruction}.

For $h\in\mathcal K_\nu$ with $\nu<2-n$, test $L_\omega h=0$
against $\chi_Rh$ using the real adjoint identity above.
The integration-by-parts argument in the proof of
Lemma~\ref{lem:ac-lichnerowicz} bounds the resulting error by
$CR^{2n+2\nu-4}\to0$, so $\mathcal D_\omega h=0$.
Therefore $J\nabla h$ is a decaying Killing field. The Killing-field
argument in that proof uses only $\delta>0$, so $J\nabla h=0$;
since $h$ decays, $h=0$. Thus $\mathcal K_\nu=0$ for $\nu<2-n$.
For $n=2$, this covers all negative weights. For $n>2$, the absence of exceptional
weights in $(4-2n,0)$, together with inclusion of weighted
spaces, extends the vanishing to every $\nu<0$.

By formal self-adjointness, we have
\[
 \operatorname{ind}L_{\omega,\lambda}
 =\dim\mathcal K_\lambda-\dim\mathcal K_{4-2n-\lambda}
 =-\operatorname{ind}L_{\omega,4-2n-\lambda}.
\]
For $n>2$, the index is zero on $(4-2n,0)$. Crossing zero adds
$m_L(0)=1$, so for $0<\lambda<\rho_+$ the cokernel is zero
and the kernel consists precisely of constants. For $n=2$,
choose $0<\epsilon<\rho_+$. By the index jump and symmetry, we have
\[
 2\dim\mathcal K_\epsilon
 =\operatorname{ind}L_{\omega,\epsilon}
       -\operatorname{ind}L_{\omega,-\epsilon}
 =m_L(0)=2,
\]
so the same conclusion holds. Inclusion proves
$\mathcal K_0=\RR$. The stated surjectivity and isomorphisms now
follow from the range formula. In the last interval, the
complementary kernel consists of constants, and
$\lambda-4<-2n$ ensures that $\int_Xq\,\omega^n$ is finite.
\end{proof}

\end{document}